\documentclass[12pt]{article}

\usepackage{latexsym}
\usepackage{amssymb}
\usepackage{amsmath}
\usepackage{enumerate}

\newtheorem{theorem}{Theorem}[section]
\newtheorem{lemma}[theorem]{Lemma}
\newtheorem{proposition}[theorem]{Proposition}
\newtheorem{corollary}[theorem]{Corollary}
\newtheorem{definition}[theorem]{Definition}

\newtheorem{remark}[theorem]{Remark}

\newcommand{\proof}{\noindent {\it {Proof. }}}

\newcommand{\prend}{ $\diamondsuit $\hfill \bigskip}

\newcommand\supp{\mathop{\rm supp}}

\newcommand\esssup{\mathop{\rm esssup}}

\newcommand\id{\mathop{\rm id}}

\newcommand\tr{\mathop{\rm tr}}

\newcommand\rank  {\mathop{\rm rank}}

\newcommand\op{\mathop{\rm op}}

\newcommand\univ{\mathop{\rm univ}}
\newcommand\dd{\mathop{\rm d}}
\newcommand\hh{\mathop{\rm h}}
\newcommand\ph{\mathop{\rm ph}}
\newcommand\ttt{\mathop{\rm t}}
\newcommand\mm{\mathop{\rm m}}

\newcommand\cb{\mathop{\rm cb}}

\newcommand\nph{\varphi}

\newcommand{\cl}[1]{\mathcal{#1}}
\newcommand{\bb}[1]{\mathbb{#1}}

\begin{document}

\title{Multidimensional operator multipliers}

\author{K. Juschenko, I. G. Todorov and L. Turowska}

\date{}

\maketitle

\begin{abstract}

We introduce multidimensional Schur multipliers and characterise
them generalising well known results by Grothendieck and Peller.
We define a multidimensional version of the two dimensional
operator multipliers studied recently by Kissin and Shulman. The
multidimensional operator multipliers are defined as elements of
the minimal tensor product of several C*-algebras satisfying
certain boundedness conditions. In the case of commutative
C*-algebras, the multidimensional operator multipliers reduce to
continuous multidimensional Schur multipliers. We show that the
multipliers with respect to some given representations of the
corresponding C*-algebras do not change if the representations are
replaced by approximately equivalent ones. We establish a
non-commutative and multidimensional version of the
characterisations by Grothendieck and Peller which shows that
universal operator multipliers can be obtained as certain weak
limits of elements of the algebraic tensor product of the
corresponding C*-algebras. 
\end{abstract}

\section{Introduction}\label{intro}

A bounded function
$\nph:\mathbb{N}\times\mathbb{N}\rightarrow\mathbb{C}$
is called a Schur multiplier if $(\nph(i,j)a_{ij})$ is the matrix
of a bounded linear operator on $\ell^2$ whenever $(a_{ij})$ is
such. The study of Schur multipliers was initiated by Schur in the
early 20th century. A characterisation of these objects was given
by A. Grothendieck in his {\it R$\acute{e}$sum$\acute{e}$}
\cite{Gro}, where he showed that Schur multipliers are precisely
the functions $\nph$ of the form $\nph(i,j) =
\sum_{k=1}^{\infty}a_k(i)b_k(j)$, where $a_k,b_k :
\bb{N}\rightarrow\bb{C}$ are such that $\sup_i \sum_{k=1}^{\infty}
|a_k(i)|^2 < \infty$ and $\sup_j \sum_{k=1}^{\infty} |b_k(j)|^2 <
\infty$. Schur multipliers have had many important applications in
Analysis, see e.g. \cite{bp}, \cite{dp} and \cite{Pi}. One of the
forms of the celebrated Grothendieck inequality can be given in
terms of these objects \cite{Pi}.

One of the most important developments in Analysis in recent years
has been \lq\lq quantisation'' \cite{effros}, starting with the
advent of the theory of operator spaces in the 1980's in the work
of Blecher, Effros, Haagerup, Paulsen, Pisier, Ruan, Sinclair and
many others, and based on Arveson's pioneering work in the 1970's.
Operator space (or non-commutative) versions are presently being
found for many results in classical Banach space theory \cite{blm,
paulsen, pisier_intr}. A construction underlying many of the
developments in Operator Space Theory is the Haagerup tensor
product, as well as its weak counterpart, the weak* Haagerup
tensor product \cite{bs} and its generalisation, the extended
Haagerup tensor product \cite{effros_ruan}. Grothendieck's
characterisation can be formulated by saying that the set of Schur
multipliers coincides with the extended (or the weak*) Haagerup
tensor product $\ell^{\infty}\otimes_{eh} \ell^{\infty}$ of the
space $\ell^{\infty}$ of all bounded complex sequences, with
itself.

Schur multipliers are elements of the commutative von Neumann
algebra $\ell^{\infty}(\bb{N}\times\bb{N})$, or equivalently of
the (von Neumann) tensor product of (the commutative von Neumann
algebra) $\ell^{\infty}$ with itself.  Subsequently, they form a
commutative algebra themselves. Their quantisation was initiated
by Kissin and Shulman in \cite{ks}. Suppose that $\cl A$ and $\cl
B$ are C*-algebras and $\pi$ and $\rho$ their representations on
$H$ and $K$, respectively. The Hilbert space tensor product
$H\otimes K$ can be naturally identified with the Hilbert space
$\cl C_2(H^{\dd},K)$ of Hilbert-Schmidt operators from the dual
$H^{\dd}$ of $H$ into $K$. It follows that $\pi$ and $\rho$ give
rise to a representation $\sigma_{\pi,\rho}$ of the minimal tensor
product $\cl A\otimes\cl B$ of $\cl A$ and $\cl B$ on $\cl
C_2(H^{\dd},K)$. Kissin and Shulman call an element $\nph\in\cl
A\otimes \cl B$ a $\pi,\rho$-multiplier if
$\sigma_{\pi,\rho}(\nph)$ is bounded in norm of $\cl
C_2(H^{\dd},K)$ induced by its inclusion into the algebra $\cl
B(H^{\dd},K)$ of all bounded operators from $H^{\dd}$ into $K$. In
\cite{ks}, they study two sets of problems: the dependence of
$\pi,\rho$-multipliers on $\pi$ and $\rho$ and the description of
the norm of an operator multiplier. Most of their results are
established in the more general setting of symmetrically normed
ideals.

Assume that $\cl A$ and $\cl B$ are commutative, say $\cl A =
C_0(X)$ and $\cl B = C_0(Y)$, for some locally compact Hausdorff
spaces $X$ and $Y$, and that the representations $\pi$ and $\rho$
arise from some spectral measures on $X$ and $Y$. The notion of a
$\pi,\rho$-multiplier is in this case closely related to double
operator integrals. The theory of these integrals was developed by
Birman and Solomyak \cite{BS1, BS2,BS3,BS4} in connection with
various problems of Mathematical Physics and in particular of
Perturbation Theory. If $(X,{\mathcal E})$ and $(Y,{\mathcal F})$
are spectral measures on  Hilbert spaces $H$ and $K$, they defined
the double operator integral
$$I_{\psi}(T)=\int_{{X}\times{Y}}\psi(x,y)\,d{\mathcal
E}(x)T\, d{\mathcal F}(y)$$ for every bounded measurable function
$\psi$ and every operator $T$ from the Hilbert-Schmidt class $\cl
C_2(H, K)$. A function $\psi$ is called a Schur multiplier with
respect to ${\mathcal E}$ and ${\mathcal F}$ if $I_{\psi}$ can be
extended to a bounded linear transformer on the space $(\cl
B(H,K), \|\cdot\|_{\op})$ of bounded operators from $H$ to $K$,
that is, if there exists $C>0$ such that
$\|I_{\psi}(T)\|_{\op}\leq C\|T\|_{\op}$ for all $T\in\cl
C_2(H,K)$. Peller \cite{peller_two_dim} (see also \cite{HK})
characterised Schur multipliers with respect to ${\mathcal E}$ and
${\mathcal F}$ in several ways. In particular, he showed that the
space of Schur multipliers with respect to ${\mathcal E}$ and
${\mathcal F}$ coincides with the extended Haagerup tensor product
$L^{\infty}({X})\otimes_{eh} L^{\infty}({Y})$ and the integral
projective tensor product $L^{\infty}({X})\hat\otimes_i
L^{\infty}({Y})$.

Several attempts were made to generalise the Birman-Solomyak
theory to the case of multiple operator integrals \cite{pavlov,
stenkin, solomyak_stenkin}. Such integrals appear, for instance,
in the study of differentiability of functions of operators
depending on a parameter. A recent definition of multiple operator
integrals by Peller in \cite{peller} is based on the integral
projective tensor product. For some fixed spectral measures $(X_1,
{\mathcal E}_1),\ldots,(X_n, {\mathcal E}_n)$ on  Hilbert spaces
$H_1,\ldots,H_n$, he defines
$$I_{\psi}(T_1,\ldots,T_{n-1})=\int_{X_1\times\ldots\times
X_n}\psi(x_1,\ldots,x_n)\, d{\mathcal E}_1(x_1)T_1\, d{\mathcal
E}_2(x_2)\ldots T_{n-1}\, d{\mathcal E}_n(x_n),$$ where $\psi\in
L^{\infty}({X}_1)\hat\otimes_i\ldots\hat\otimes_i
L^{\infty}({X}_n)$ and $T_1,\ldots,T_{n-1}$ are bounded linear
operators, and shows that
$$\|I_{\psi}(T_1,\ldots,T_{n-1})\|_{\op}\leq
\|\psi\|_{i}\|T_1\|_{\op}\ldots \|T_{n-1}\|_{\op},$$ where
$\|\psi\|_i$ denotes the integral projective tensor norm of
$\psi$. If the spectral measures are multiplicity free and
$T_1,\ldots,T_{n-1}$ are Hilbert-Schmidt operators with kernels
$f_1,\dots,f_{n-1}$, respectively, then $I_{\psi}(T_1,\ldots,$
$T_{n-1})$ is a Hilbert-Schmidt operator with kernel
$S_{\psi}(f_1,\ldots,f_{n-1})\in L^2({X}_1\times{X}_n)$ equal to
\begin{equation}\label{mult}\int_{{X}_2\times\ldots\times {X}_{n-1}}
\psi(x_1,\ldots,x_n)f_1(x_1,x_2)\ldots f_{n-1}(x_{n-1},x_n)\, d
{\mathcal E}_2(x_2)\ldots\, d{\mathcal E}_{n-1}(x_{n-1}).
\end{equation}

This was the starting point for our definition of multidimensional
Schur multipliers in Section~\ref{s_mm2}. Let $(X_i,\mu_i)$,
$i=1,\ldots,n$, be standard $\sigma$-finite measure spaces and
$\Gamma(X_1,\ldots,X_n)=L^2(X_1\times X_2)\odot L^2(X_2\times
X_3)\odot\ldots\odot L^2(X_{n-1} \times X_n)$ be the algebraic
tensor product of the corresponding $L^2$-spaces equipped with the
projective tensor norm, where each of the $L^2$-spaces is equipped
with its $L^2$-norm. An element $\psi\in
L^{\infty}(X_1\times\dots\times X_n)$ determines a bounded linear
map $S_{\psi}$ from $\Gamma(X_1,\ldots,X_n)$ to $L^2(X_1,X_n)$
given on elementary tensors $f_1\otimes\ldots\otimes
f_n\in\Gamma(X_1,\ldots,X_n)$ by (\ref{mult}) (where the
integration is now with respect to $\mu_i$ instead of ${\mathcal
E}_i$). On the other hand, for any measure spaces $(X,\mu)$ and
$(Y,\nu)$, the space $L^2(X\times Y)$ can be identified with the
class of all Hilbert-Schmidt operators from $L^2(X)$ to $L^2(Y)$;
to each $f\in L^2(X\times Y)$ there corresponds the operator $T_f$
given by $T_f\xi (y) = \int_{X}f(x,y)\xi(x)d\mu(x)$, $\xi\in
L^2(X)$. Using this identification, one can equip the space
$L^2(X\times Y)$ with the opposite operator space structure
arising from the inclusion of $L^2(X\times Y)$ into $\cl
B(L^2(X),L^2(Y))$. We further equip $\Gamma(X_1,\ldots, X_n)$ with
the Haagerup tensor norm $\|\cdot\|_{\hh}$, where the $L^2$-spaces
are given their opposite operator space structure described above,
and say that an element $\psi\in L^{\infty}(X_1\times \ldots\times
X_n)$ is a Schur multiplier (with respect to $\mu_1,\ldots,\mu_n$)
if there exists $C > 0$ such that
\begin{equation}\label{def}
\|S_{\psi}(\Phi)\|_{\op}\leq C\|\Phi\|_{\hh}, \text{ for all
}\Phi\in\Gamma(X_1,\ldots,X_n).
\end{equation}
Using a generalisation of a result of Smith \cite{smith} on the
complete boundedness of certain bounded bimodule maps to the case
of multilinear modular maps, we obtain a characterisation of
multidimensional Schur multipliers as elements of the extended
Haagerup tensor product
$L^{\infty}(X_1)\otimes_{eh}\ldots\otimes_{eh} L^{\infty}(X_n)$
(Theorem~\ref{th_g3}). This generalises Grothendieck's and
Peller's characterisations in the case $n=2$. We show that the
integral projective tensor product consists of multipliers and,
therefore, $L^{\infty}(X_1)\hat\otimes_i\ldots\hat\otimes_i
L^{\infty}(X_n)\subset
L^{\infty}(X_1)\otimes_{eh}\ldots\otimes_{eh}L^{\infty}(X_n)$. The
converse inclusion is true in the case $n=2$ \cite{peller_two_dim}
but remains an open problem for $n > 2$.

In Section~\ref{s_nonc} we consider a non-commutative version of
multidimensional multipliers following the Kissin-Shulman approach
in the two dimensional case. We replace the functions $\psi$ by
elements of the minimal tensor product ${\mathcal
A}_1\otimes\ldots\otimes {\mathcal A}_n$ of some given C*-algebras
${\mathcal A}_1,\ldots,{\mathcal A}_n$ and the measure $\mu_i$ by
a representation $\pi_i$ of ${\mathcal A}_i$. We thus obtain a
class of operator $\pi_1,\ldots,\pi_n$-multipliers. If each
${\mathcal A}_i$ is a commutative C*-algebra, say $\cl A_i =
C_0(X_i)$ for some locally compact Hausdorff space $X_i$, and
$\pi_i(f)$ is the operator of multiplication by $f\in C_0(X)$
acting on $L^2(X_i,\mu_i)$, then $\psi$ is a
$\pi_1,\ldots,\pi_n$-multiplier if and only if $\psi$ is a Schur
multiplier with respect to $\mu_1,\dots,\mu_n$
(Proposition~\ref{com+noncom}). As in the two-dimensional case, we
show that the set of $\pi_1,\ldots,\pi_n$-multipliers does not
change if we replace each $\pi_i$ by an  approximately equivalent
representation (Theorem~\ref{ap}). A consequence of this result is
the fact that the class of continuous (multidimensional) Schur
multipliers depends only on the supports of the measures $\mu_i$.

In Section~\ref{s_universal} we study universal mutlipliers, that
is, the elements of ${\mathcal A}_1\otimes\ldots\otimes{\mathcal
A}_n$ which are $\pi_1,\ldots,\pi_n$-multipliers for all
representations $\pi_i$ of ${\mathcal A}_i$, $i = 1,\dots,n$. We
characterise such multipliers as the elements of a certain weak
completion of the algebraic tensor product ${\mathcal
A}_1\odot\ldots\odot{\mathcal A}_n$ (Theorem~\ref{desc}). In the
case where the C*-algebras are commutative and $n=2$ this was proved
in \cite{ks}; the case of arbitrary C*-algebras was left as a
conjecture. Our result may be thought of as a non-commutative and
multidimensional version of Grothendieck's and Peller's
characterisations of Schur multipliers. The key ingredient in the
proof is the observation that a universal multiplier determines a
completely bounded multilinear modular map from the Cartesian
product of the C*-algebras of compact operators into the C*-algebra
of compact operators which allows us to use a result by Christensen
and Sinclair \cite{cs} providing a description of all such mappings.

\vspace{0.3cm}
\begin{center}
{\bf Acknowledgements}
\end{center}
We are grateful to V.S. Shulman for stimulating results, questions
and discussions. We thank T. Itoh and A. Katavolos for providing
the references \cite{itoh} and \cite{harris}, respectively, and R.
Levene for pointing out a mistake in an earlier version of the
paper. We are indebted to the referee for a number of suggestions
which helped improving the manuscript.

The work was partially written when the second author was visiting
Chalmers University of Technology in G\"oteborg, Sweden, supported
by funds from Queen' s University Belfast. The third author was
partially supported by the Swedish Reseach Council, while the
first and the third authors were partially supported by
Engineering and Physical Sciences Research Council grant
EP/D050677/1.

\section{Preliminaries}\label{prel}

In this section we collect some preliminary notions and results
which will be needed in the sequel.

Let $H$ be a Hilbert space. The dual space $H^{\dd}$ of $H$ is a
Hilbert space and there exists an anti-isometry $\partial :
H\rightarrow H^{\dd}$ given by $\partial(x)(y) = (y,x)$, $x,y\in
H$. We set $x^{\dd} =
\partial(x)$.

If $H$ and $K$ are Hilbert spaces, we let $\cl B(H,K)$ be the
space of all bounded linear operators from $H$ into $K$, and
$\|\cdot\|_{\op}$ be the usual operator norm on $\cl B(H,K)$. We
let $\cl K(H,K)$ be the subspace of all compact operators, and
$\cl C_2(H,K)$ be the subspace of all Hilbert-Schmidt operators,
from $H$ into $K$. For each $T\in \cl C_2(H,K)$, we denote by
$\|T\|_2$ the Hilbert-Schmidt norm of $T$. The space $\cl
C_2(H,K)$ is a Hilbert space with respect to the inner product
$(T,S) = \tr(TS^*)$, where $S^*$ denotes the adjoint of the
operator $S$. We let $\cl B(H) = \cl B(H,H)$, $\cl K(H) = \cl
K(H,H)$ and $\cl C_2(H) = \cl C_2(H,H)$.

If $T\in \cl B(H,K)$ we denote by $T^{\dd}\in B(K^{\dd},H^{\dd})$
the conjugate of $T$. We have that $\|T^{\dd}\|_{\op} =
\|T\|_{\op}$ and $T^{\dd}x^{\dd}=(T^*x)^{\dd}$, whenever $x\in
H_{2}$. Another way of expressing the last identity is
\begin{equation}
T^{\dd}={\partial} T^*{\partial}^{-1}. \label{partial}%
\end{equation}
We also have
\begin{equation}
(T^*)^{\dd}=(T^{\dd})^* \ \mbox{ and } \
({\lambda}T)^{\dd}={\lambda} T^{\dd},\ \ \ {\lambda}\in \bb{C}. \label{e1.5}%
\end{equation}

We let $H\otimes K$ be the Hilbert space tensor product of $H$ and
$K$. There exists a unitary operator $\theta : H\otimes
K\rightarrow \cl C_2(H^{\dd},K)$ given on elementary tensors
$x\otimes y\in H\otimes K$ by
$$\theta(x\otimes y)(z^{\dd}) =
(x,z)y, \ \ \ z^{\dd}\in H^{\dd}.$$ If $A\in\cl B(H)$, $B\in\cl
B(K)$, $x\in H$ and $y\in K$, we have that $\theta((A{\otimes}
B)(x{\otimes}y))$ $=$ $B{\theta}(x{\otimes} y)A^{\dd}$, and hence
\begin{equation}\label{e1.6}
{\theta}((A{\otimes} B){\xi})=B{\theta}({\xi})A^{\dd}\ \ \mbox{
for all }{\xi}\in H{\otimes} K.
\end{equation}

If $\nph\in\cl B(H\otimes K)$, let $\sigma(\nph)\in \cl B(\cl
C_2(H^{\dd},K))$ be given by the formula
$$\sigma(\nph)\theta(\xi) = \theta(\nph\xi), \ \
\xi\in H\otimes K.$$ Then $\sigma$ implements a unitary
equivalence between $\cl B(H\otimes K)$ and $\cl B(\cl C_2$
$(H^{\dd},$ $K))$. We will call an element $\nph\in \cl B(H\otimes
K)$ a concrete (operator) multiplier if there exists $C > 0$ such
that $\|\sigma(\nph)T\|_{\op}\leq C\|T\|_{\op}$, for each $T\in\cl
C_2(H^{\dd},K)$. Suppose that $H = l^2(X)$, $K = l^2(Y)$ for some
sets $X$ and $Y$ and $\nph$ is the operator on $H\otimes K =
\ell^2(X\times Y)$ of multiplication by a function $\phi\in
\ell^{\infty}(X\times Y)$. The concrete operator multipliers of
this form are precisely the classical Schur multipliers on
$X\times Y$ (see e.g. \cite{Pi}).

Let $\cl A$ and $\cl B$ be C*-algebras. We denote by $\cl
A\otimes\cl B$ the minimal tensor product of $\cl A$ and $\cl B$.
Let $\pi : \cl A\rightarrow\cl B(H)$ (resp. $\rho : \cl
B\rightarrow\cl B(K)$) be a representation of $\cl A$ (resp. $\cl
B$). Then $\pi\otimes\rho : \cl A\otimes\cl B\rightarrow \cl
B(H\otimes K)$, given on elementary tensors by
$(\pi\otimes\rho)(a\otimes b) = \pi(a)\otimes\rho(b)$, is a
representation of $\cl A\otimes \cl B$. Let $\sigma_{\pi,\rho} =
\sigma\circ (\pi\otimes\rho)$; clearly, $\sigma_{\pi,\rho}$ is a
representation of $\cl A\otimes\cl B$ on $\cl C_2(H^{\dd},K)$,
unitarily equivalent to $\pi\otimes\rho$. We moreover have
$$\sigma_{\pi,\rho}(a\otimes b)T = \rho(b)T\pi(a)^{\dd}, \ \
a\in \cl A, b \in \cl B, T\in\cl C_2(H^{\dd},K).$$ An element
$\nph\in\cl A\otimes\cl B$ is called a $\pi,\rho$-multiplier
\cite{ks} if there exists $C > 0$ such that
\begin{equation}\label{2dm}
\|\sigma_{\pi,\rho}(\nph)T\|_{\op}\leq C\|T\|_{\op}, \ \ \
\mbox{for each } T\in\cl C_2(H^{\dd},K),
\end{equation}
in other words, if $(\pi\otimes\rho)(\nph)$ is a concrete operator
multiplier. The set of all $\pi,\rho$-multipliers in $\cl
A\otimes\cl B$ is denoted by ${\bf M}_{\pi,\rho}(\cl A,\cl B)$,
and the smallest constant $C$ appearing in (\ref{2dm}) is denoted
by $\|\nph\|_{\pi,\rho}$. If $\nph$ is a $\pi,\rho$-multiplier for
all representations $\pi$ of $\cl A$ and $\rho$ of $\cl B$ then
$\nph$ is called a universal multiplier. The set of all universal
multipliers is denoted by ${\bf M}(\cl A,\cl B)$; if $\nph\in {\bf
M}(\cl A,\cl B)$ we let $\|\nph\|_{\univ} = \sup_{\pi,\rho}
\|\nph\|_{\pi,\rho}$. It is not difficult to see that in this case
$\|\nph\|_{\univ} < \infty$ \cite{ks}.

We now recall some notions from Operator Space Theory. We refer
the reader to \cite{blm}, \cite{er} and \cite{pisier_intr} for
more details. An operator space $\cl E$ is a closed subspace of
$\cl B(H,K)$, for some Hilbert spaces $H$ and $K$. If
$n,m\in\bb{N}$, by $M_{n,m}(\cl E)$ we will denote the space of
all $n$ by $m$ matrices with entries in $\cl E$ and let $M_n(\cl
E) = M_{n,n}(\cl E)$. Note that $M_{n,m}(\cl E)$ can be identified
in a natural way with a subspace of $\cl B(H^m, K^n)$ and hence
carries a natural operator norm. If $n = \infty$ or $m = \infty$,
we will denote by $M_{n,m}(\cl E)$ the space of all (singly or
doubly infinite) matrices with entries in $\cl E$ which represent
a bounded linear operator between the corresponding amplifications
of the Hilbert spaces and set $M_{\infty}(\cl E) =
M_{\infty,\infty}(\cl E)$. We also write $M_{n,m} =
M_{n,m}(\bb{C})$ and $M_{\infty} = M_{\infty,\infty}(\bb{C})$. If
$a=(a_{ij})\in M_{n,m}(\cl E)$, where $a_{ij}\in \cl E$, we let
$a^{\dd} = (a_{ij}^{\dd})$; thus $a^{\dd}\in \cl
B(K^{\dd,m},H^{\dd,n})$. We also let $a^{\ttt} = (a_{ji})\in
M_{m,n}(\cl E)$; thus $a^{\ttt}\in \cl B(H^{n},K^{m})$. We have
$\|a^{\dd}\|_{\op}=\|a^{\ttt}\|_{\op}$ and $\|a^{\dd,\ttt}\|_{\op}
= \|a\|_{\op}$. The opposite operator space $\cl E^o$ of the
operator space $\cl E$ is defined as follows: if $\cl
E\subseteq\cl B(H,K)$ then $\cl E^o = \{x^{\dd} : x\in\cl
E\}\subseteq\cl B(K^{\dd},H^{\dd})$.

If $\cl E$ and $\cl F$ are operator spaces, a linear map $\Phi :
\cl E\rightarrow \cl F$ is called completely bounded if the map
$\Phi^{(k)} : M_k(\cl E)\rightarrow M_k(\cl F)$, given by
$\Phi^{(k)}((a_{ij})) = (\Phi(a_{ij}))$, is bounded for each
$k\in\bb{N}$ and $\|\Phi\|_{\cb} \stackrel{def}{=} \sup_k
\|\Phi^{(k)}\| < \infty$.

Let $\cl E, \cl E_1,\dots,\cl E_n$ be operator spaces. We denote
by $\cl E_1\odot\dots\odot \cl E_n$ the algebraic tensor product
of $\cl E_1,\dots,\cl E_n$. Let $a_k = (a_{ij}^k)\in
M_{m_k,m_{k+1}}(\cl E_k)$, $k = 1,\dots,n$. We denote by
\begin{equation}\label{odot}
a^1\odot \dots\odot a^n\in M_{m_1,m_{n+1}}(\cl E_1\odot\dots\odot
\cl E_n)
\end{equation}
the matrix whose $i,j$-entry is
\begin{equation}\label{odot2}
\sum_{i_2,\dots,i_n} a^1_{i,i_2}\otimes
a^2_{i_2,i_3}\otimes\dots\otimes a^n_{i_n,j}.
\end{equation}
Let $\Phi : \cl E_1\times\dots\times \cl E_n \rightarrow\cl E$ be
a  multilinear map and
$$\Phi^{(m)} : M_m(\cl E_1)\times M_m(\cl E_2)\times\dots\times M_m(\cl E_n)
\rightarrow M_m(\cl E)$$ be the multilinear map given by
\begin{equation}\label{rfl}
\Phi^{(m)}(a^1,\dots,a^n)_{ij} = \sum_{i_2,\dots,i_n}
\Phi(a^1_{i,i_2},a^2_{i_2,i_3},\dots, a^n_{i_n,j}),
\end{equation}
where $a^k = (a^k_{ij})\in M_m(\cl E_k)$, $1\leq i,j\leq m$. The
map $\Phi$ is called completely bounded if there exists $C
> 0$ such that for all $m\in\bb{N}$ and all elements $a^k\in
M_{m}(\cl E_k)$, $k = 1,\dots,n$, we have
$$\|\Phi^{(m)}(a^1,\dots,a^n)\|\leq C \|a^1\|\dots\|a^n\|.$$

Every completely bounded multilinear map $\Phi : \cl
E_1\times\dots\times \cl E_n \rightarrow\cl E$ gives rise to a
completely bounded linear map from the Haagerup tensor product
$\cl E_1\otimes_{\hh}\dots\otimes_{\hh} \cl E_n$ into $\cl E$. For
details on the Haagerup tensor product we refer the reader to
\cite{er}.

If $R_1,\dots,R_{n+1}$ are rings, $M_i$ is a $R_i$-left and
$R_{i+1}$-right module for each $i = 1,\dots,n$, and $M$ is an
$R_1,R_{n+1}$-module, a multilinear map $\Phi : M_1\times\dots
\times M_n\rightarrow M$ will be called
$R_1,\dots,R_{n+1}$-modular (or simply modular if
$R_1,\dots,R_{n+1}$ are clear from the context) if
$$\Phi(a_1 m_1 a_2,m_2 a_3,m_3 a_4,\dots,m_n a_{n+1})
= a_1\Phi(m_1, a_2 m_2, a_3 m_3,\dots,a_n m_n) a_{n+1},$$ for all
$m_i\in M_i$ ($i = 1,\dots,n$) and $a_j\in R_j$ ($j =
1,\dots,n+1$). If $R_i = \cl A_i$ are C*-algebras and $M_i = \cl
E_i$ are operator spaces, we let $\cl
B_{\mathcal{A}_1,\ldots,\mathcal{A}_{n+1}}
(\mathcal{E}_1,\ldots,\mathcal{E}_n; \mathcal{E})$ (resp.
$CB_{\mathcal{A}_1,\ldots,\mathcal{A}_{n+1}}(\mathcal{E}_1,\ldots,\mathcal{E}_n;
\mathcal{E})$) denote the spaces of all bounded (resp. completely
bounded) $\cl A_1,\dots,\cl A_{n+1}$-modular maps from $\cl
E_1\times\dots\times\cl E_n$ into $\cl E$.

\section{Multidimensional Schur multipliers}\label{s_mm2}

In this section, we define multidimensional Schur multipliers on
the direct product of finitely many measure spaces. The main
result of the section is Theorem \ref{th_g3} which characterises
multidimensional Schur multipliers generalising the results of
Peller \cite{peller_two_dim} and Spronk \cite{spronk}.

Let $(X_i,\mu_i)$, $i = 1,2,\dots,n$, be standard $\sigma$-finite
measure spaces. For notational convenience, integration with
respect to $\mu_i$ will be denoted by $dx_i$. Direct products of
the form $X_{i_1}\times\dots\times X_{i_k}$ will be equipped with
the corresponding product measure. We equip the space
$L^2(X_1\times X_2)$ with an
$L^{\infty}(X_1),L^{\infty}(X_2)$-module  action by letting $(a
\xi b)(x,y) = a(x)\xi(x,y)b(y)$. We will denote by $M_a$ the
operator of multiplication by the essentially bounded function $a$
acting on the corresponding $L^2$-space.

\begin{theorem}\label{th_mod}
A multilinear map
$$S : L^2(X_1\times X_2)\times L^2(X_2\times
X_3)\times\dots\times L^2(X_{n-1}\times X_n)\rightarrow
L^2(X_1\times X_n)$$ is a bounded modular map if and only if there
exists $\nph\in L^{\infty}(X_1\times\dots\times X_n)$ such that $S
= S_{\nph}$ where $S_{\nph}(f_1,\dots,f_{n-1})(x_1,x_n)$ is
defined as
$$\int_{X_2\times\dots\times X_{n-1}}
\nph(x_1,\dots,x_n)f_1(x_1,x_2)f_2(x_2,x_3)\dots
f_{n-1}(x_{n-1},x_n) dx_2\dots dx_{n-1}.$$ Moreover, $\|S_{\nph}\|
= \|\nph\|_{\infty}$.
\end{theorem}
\proof We first show that for each $\nph$, the map $S_{\nph}$ is a
bounded modular map with norm not exceeding $\|\nph\|_{\infty}$.
For simplicity, we will assume in this part of the proof that $n =
3$. Fix $\nph$, $f_1$ and $f_2$. We have
\begin{eqnarray*}
& & \|S_{\nph}(f_1,f_{2})\|_2^2\\ & \leq & \int_{X_1\times X_3}
\left(\int |\nph(x_1,x_2,x_3) f_1(x_1,x_2)f_{2}(x_2,x_3)|dx_2\right)^2 dx_1 dx_3\\
& \leq & \|\nph\|_{\infty}^2 \int_{X_1\times X_3} \left(\int |
f_1(x_1,x_2)f_{2}(x_2,x_3)|dx_2\right)^2 dx_1 dx_3\\
& \leq & \|\nph\|_{\infty}^2 \int_{X_1\times X_3} \left(\int |
f_1(x_1,x_2)|^2dx_2\right) \left(\int |
f_{2}(x_2,x_2)|^2dx_2\right) dx_1 dx_3\\
& = & \|\nph\|_{\infty}^2 \|f_1\|_2^2\|f_2\|_2^2.
\end{eqnarray*}
Thus, $\nph$ is bounded with $\|S_{\nph}\| \leq
\|\nph\|_{\infty}$; the modularity of $S_{\nph}$ is obvious.

Conversely, let
$$S : L^2(X_1\times X_2)\times L^2(X_2\times
X_3)\times\dots\times L^2(X_{n-1}\times X_n)\rightarrow
L^2(X_1\times X_n)$$ be a bounded modular map. We first assume
that the measures $\mu_i$ are finite. Write $K_1 = L^2(X_1\times
X_n)$ and let
$$S_1 : L^2(X_2)\times L^2(X_2)\times L^2(X_3)\times
L^2(X_3)\times\dots\times L^2(X_{n-1})\times L^2(X_{n-1})
\rightarrow K_1$$ be given by
$$S_1(\xi_2,\eta_2,\xi_3,\eta_3,\dots,\xi_{n-1},\eta_{n-1}) =
S(1\otimes\xi_2,\eta_2\otimes\xi_3,\dots,\eta_{n-1}\otimes 1)$$
(here and in the sequel we denote by $1$ the constant function
taking value one). The fact that $S$ is modular implies that
$$S_1(\xi_2a_2,\eta_2,\xi_3a_3,\dots,\xi_{n-1}a_{n-1},\eta_{n-1})
= S_1(\xi_2,a_2\eta_2,\xi_3,\dots,a_{n-1}\eta_{n-1}),$$ whenever
$a_i\in L^{\infty}(X_i)$, $i = 2,\dots,n-1$. For fixed
$\xi_3,\eta_3,\dots,\xi_{n-1},\eta_{n-1}$, let $S_2 :
L^2(X_2)\times L^2(X_2) \rightarrow K_1$ be given by
$$S_2(\xi_2,\eta_2) =
S_1(\xi_2,\eta_2,\xi_3,\eta_3,\dots,\xi_{n-1},\eta_{n-1}).$$ For
$h\in K_1$, let $S_2^h : L^2(X_2)\times L^2(X_2) \rightarrow
\bb{C}$ be defined by $S_2^h(\xi_2,\eta_2) =
(S_2(\xi_2,\eta_2),h)$. Clearly,
$$|S_2^h(\xi_2,\eta_2)|\leq
\|h\|\|S\|\prod_{i=2}^{n-1}\|\xi_i\|\|\eta_i\|.$$ Hence there
exists a bounded operator $T_2^h : L^2(X_2)\rightarrow L^2(X_2)$
such that $S_2^h(\xi_2,\eta_2) = (T_2^h\xi_2,\overline{\eta_2})$,
for all $\xi_2,\eta_2\in L^2(X_2)$ and $\|T_2^h\|$ $\leq$
$\|h\|\|S\|\prod_{i=3}^{n-1}$ $\|\xi_i\|$ $\|\eta_i\|$. For each
$a\in L^{\infty}(X_2)$ and $\xi_2,\eta_2\in L^2(X_2)$, we have
that
\begin{eqnarray*}
(T_2^hM_{a}\xi_2,\overline{\eta_2}) & = & S_2^h(a\xi_2,\eta_2) =
S_2^h(\xi_2,a\eta_2)\\ & = & (T_2^h\xi_2, \overline{a\eta_2}) =
(T_2^h\xi_2, M_{\overline{a}}\overline{\eta_2}) =
(M_{a}T_2^h\xi_2,\overline{\eta_2}).
\end{eqnarray*} Thus, there exists $\nph_2^h \in L^{\infty}(X_2)$
such that $T_2^h = M_{\nph_2^h}$. Moreover,
$$\|\nph_2^h\|_{\infty} \leq
\|h\|\|S\|\prod_{i=3}^{n-1}\|\xi_i\|\|\eta_i\|.$$ For each $f\in
L^1(X_2)$, the functional on $K_1$ given by $h \rightarrow
\int_{X_2} f(x_2)\nph_2^h(x_2) dx_2$ is conjugate linear and
bounded with norm not exceeding $\|f\|_1$ $\|S\|\prod_{i=3}^{n-1}$
$\|\xi_i\|\|\eta_i\|.$ Hence, there exists $\Phi_2(f)\in K_1$ such
that
$$(\Phi_2(f),h) = \int_{X_2} f(x_2)\nph_2^h(x_2) dx_2,$$ and
$\|\Phi_2(f)\|_{K_1}\leq
\|f\|_1\|S\|\prod_{i=3}^{n-1}\|\xi_i\|\|\eta_i\|.$ Thus, the
mapping $\Phi_2 : L^1(X_2)\rightarrow K_1$ is bounded and
$\|\Phi_2\|$ $ \leq$ $\|S\|$ $\prod_{i=3}^{n-1}$
$\|\xi_i\|\|\eta_i\|.$ Since Hilbert spaces possess Radon-Nikodym
property, the vector valued Riesz Representation Theorem
\cite[Theorem 5, p. 63]{du} implies that there exists $\nph_2\in
L^{\infty}(X_2,K_1)$ ($L^{\infty}(X_2,K_1)$ being the space of
essentially bounded $K_1$-valued measurable functions on $X_2$)
such that
$$\Phi_2(f) = \int_{X_2} f(x_2)\nph_2(x_2)dx_2,$$ where the
integral is in Bochner's sense. Moreover,
$$\|\nph_2\|_{L^{\infty}(X_2,K_1)} = \esssup_{x_2\in X_2} \|\nph_2(x_2)\|_{K_1}
= \|\Phi_2\| \leq \|S\|\prod_{i=3}^{n-1}\|\xi_i\|\|\eta_i\|.$$

For $\xi_2,\eta_2\in L^2(X_2)$, we have that
$\xi_2\overline{\eta_2}\in L^1(X_2)$ and hence
\begin{eqnarray*}
(S_2(\xi_2,\eta_2),h) & = & (T_2^h\xi_2,\overline{\eta_2}) =
\int_{X_2}\nph_2^h(x_2)\xi_2(x_2)\eta_2(x_2) dx_2\\ & = &
\left(\int_{X_2}\nph_2(x_2)\xi_2(x_2)\eta_2(x_2)dx_2, h\right);
\end{eqnarray*}
in other words,
$$S_2(\xi_2,\eta_2) =
\int_{X_2}\nph_2(x_2)\xi_2(x_2)\eta_2(x_2)dx_2,$$ where the
integral is in Bochner's sense.

We consider $\nph_2$ as a function on $X_1\times X_2\times X_n$ by
letting $\nph_2(x_1,x_2,x_n) = \nph_2(x_2)(x_1,x_n)$. Note that
$\nph_2$ depends on $\xi_3,\eta_3,\dots,\xi_{n-1},\eta_{n-1}$; we
denote this dependence by $\nph_2 =
\nph_{2,\xi_3,\eta_3,\dots,\xi_{n-1},\eta_{n-1}}$.

Let $K_2 = L^2(X_1\times X_2\times X_n)$. We have
\begin{eqnarray*}
\|\nph_2\|_{K_2} & = & \int_{X_2}\int_{X_1\times
X_n}|\nph_2(x_2)(x_1,x_n)|^2dx_1 dx_n dx_2 =
\int_{X_2}\|\nph_2(x_2)\|_{K_1}^2dx_2\\ & \leq &
\mu_2(X_2)\|\nph_2\|_{L^{\infty}(X_2,K_1)}.
\end{eqnarray*} It follows that the mapping $S_3 :
L^2(X_3)\times L^2(X_3)\rightarrow K_2$ given by
$$S_3(\xi_3,\eta_3) =
\nph_{2,\xi_3,\eta_3,\dots,\xi_{n-1},\eta_{n-1}}$$ is well-defined
and
$$\|S_3(\xi_3,\eta_3)\|_{K_2}\leq
\mu_2(X_2)\|S\|\prod_{i=3}^{n-1}\|\xi_i\|\|\eta_i\|.$$ Hence,
$S_3$ is bounded and $\|S_3\|\leq
\mu_2(X_2)\|S\|\prod_{i=4}^{n-1}\|\xi_i\|\|\eta_i\|.$ An argument
similar to the above implies the existence of $\nph_3\in
L^{\infty}(X_3,K_2)$ with
$$\|\nph_3\|_{L^{\infty}(X_3,K_2)} \leq
\mu_2(X_2)\|S\|\prod_{i=4}^{n-1}\|\xi_i\|\|\eta_i\|$$ such that
$$S_3(\xi_3,\eta_3) = \int_{X_3}\nph_3(x_3)\xi_3(x_3)\eta_3(x_3)dx_3,$$ where the
integral is in Bochner's sense. We may consider $\nph_3$ as a
function on $X_1\times X_2\times X_3\times X_n$ by letting
$\nph_3(x_1,x_2,x_3,x_n) = \nph_3(x_3)(x_1,x_2,x_n)$. We express
the dependence of $\nph_3$ on $\xi_4,\dots,\eta_{n-1}$ by writing
$\nph_3 = \nph_{3,\xi_4,\dots,\eta_{n-1}}$. We have that
\begin{eqnarray*}
& & S_1(\xi_2,\eta_2,\dots,\xi_{n-1},\eta_{n-1}) =\\
& &
\int_{X_2}\int_{X_3}\nph_{3,\xi_4,\dots,\eta_{n-1}}(x_1,x_2,x_3,x_n)\xi_2(x_2)\eta_2(x_2)\xi_3(x_3)\eta_3(x_3)
dx_3 dx_2,
\end{eqnarray*}
where both integrals are in Bochner's sense.

Continuing inductively, we obtain $\nph\in L^{\infty}(X_{n-1},
K_{n-2})$, where $K_{n-2} = L^2(X_1\times\dots\times X_{n-2}\times
X_n)$, such that
\begin{eqnarray*}
& & S_1(\xi_2,\eta_2,\dots,\xi_{n-1},\eta_{n-1}) =\\
& &
\int_{X_2}\dots\int_{X_{n-1}}\nph(x_1,\dots,x_n)\xi_2\eta_2\dots\xi_{n-1}\eta_{n-1}
dx_{n-1}\dots dx_2,
\end{eqnarray*}
where the integrals are understood in Bochner's sense and $\nph$
is viewed as a function on $X_1\times\dots\times X_n$ by letting
$\nph(x_1,\dots,x_n) = \nph(x_{n-1})(x_1,\dots,x_{n-2},x_n)$.

It is easy to see that if $\psi\in L^1(Y,L^2(Z))$, where $Y$ and
$Z$ are finite measure spaces, then $\int_{Y\times
Z}|\psi(y)(z)|dydz$ is finite and $\left(\int_Y\psi(y)
dy\right)(z) = \int_Y\psi(y)(z) dy$, for almost all $z\in Z$ (the
first integral is in Bochner's sense, while the second one is a
Lebesgue integral with respect to the variable $y$). It now
follows that the last equality holds when the integrals are
interpreted in the sense of Lebesgue.

The modularity of $S$ implies
\begin{eqnarray*}
& & S(a\otimes\xi_2,\eta_2\otimes\xi_3,\dots,\eta_{n-1}\otimes b) =\\
& &
\int_{X_2}\int_{X_3}\dots\int_{X_{n-1}}\nph(x_1,\dots,x_n)a\xi_2\eta_2\dots\xi_{n-1}\eta_{n-1}b
dx_{n-1}\dots dx_2,
\end{eqnarray*}
for all $a\in L^{\infty}(X_1)$, $b\in L^{\infty}(X_n)$ and
$\xi_i,\eta_i\in L^2(X_i)$, $i = 2,\dots,n-1$. Letting $a =
\chi_{\alpha_1}$, $b = \chi_{\alpha_n}$ and $\xi_i = \eta_i =
\chi_{\alpha_i}$, $i = 2,\dots,n-1$, the boundedness of $S$
implies
$$\int_{\alpha_1\times\dots\times \alpha_n}|\nph(x_1,\dots,x_n)|dx_1\dots
dx_n \leq \|S\|\mu_1(\alpha_1)\dots\mu_n(\alpha_n).$$ It follows
that the mapping
$$f = \sum_{i=1}^N
\lambda_i\chi_{\alpha_1^i\times\dots\times\alpha_n^i}\longrightarrow
\int_{X_1\times\dots\times X_n} \nph f, $$ where
$\{\alpha_1^i\times\dots\times\alpha_n^i\}$ is a finite family of
disjoint Borel rectangles, is a linear functional on a dense
subspace of $L^1(X_1\times\dots\times X_n)$ of norm not exceeding
$\|S\|$. Therefore, $\nph\in L^{\infty}(X_1\times\dots\times X_n)$
and $\|\nph\|_{\infty}\leq \|S\|$.

We have that the mappings $S$ and $S_{\nph}$ coincide on the
tuples of the form
$a\otimes\xi_2,\eta_2\otimes\xi_3,\dots,\eta_{n-1}\otimes b$; by
linearity and continuity, they are equal. By the first part of the
proof, $\|S\|\leq \|\nph\|_{\infty}$ and hence $\|\nph\|_{\infty}
= \|S\|$.

Now relax the assumption on the finiteness of $\mu_i$, and let
$X_i^k$, $k\in\bb{N}$, be a measurable subset of $X_i$ such that
$\mu_i(X_i^k) < \infty$, $X_i^{k}\subseteq X_i^{k+1}$ and $X_i =
\cup_{k=1}^{\infty}X_i^k$, $i = 1,\dots,n$. For each $k\in\bb{N}$,
let
$$S_k : L^2(X_1^k\times X_2^k)\times L^2(X_2^k\times
X_3^k)\times\dots\times L^2(X_{n-1}^k\times X_n^k)\rightarrow
L^2(X_1^k\times X_n^k)$$ be the map given by
$S_k(f_1,\dots,f_{n-1}) = S(\tilde{f}_1,\dots,\tilde{f}_{n-1})$,
where $\tilde{f}_i$ coincides with $f_i$ on $X_i^k$ and is equal
to zero on the complement of $X_i^k$. Since
\begin{eqnarray*}
S_k(f_1,\dots,f_{n-1}) & = &
S(\chi_{X_1^k}\tilde{f}_1,\dots,\tilde{f}_{n-1}\chi_{X_n^k})\\ & =
& \chi_{X_1^k}S(\tilde{f}_1,\dots,\tilde{f}_{n-1})\chi_{X_n^k},
\end{eqnarray*}
the map $S_k$ is well-defined and $\|S_k\|\leq \|S\|$. Since $S_k$
is obviously $L^{\infty}(X_n^k)$, $\dots,$ $L^{\infty}$
$(X_1^k)$-modular, the above paragraphs imply that there exists
$\nph_k\in L^{\infty}(X_1^k\times\dots\times X_n^k)$ such that
$S_k = S_{\nph_k}$, for each $k\in\bb{N}$. The space
$L^2(X_i^k\times X_{i+1}^k)$ can be considered as a subspace of
$L^2(X_i^{k+1}\times X_{i+1}^{k+1})$ in a natural way. We have
that the restriction of $S_{k+1}$ to $L^2(X_1^k\times X_2^k)$
$\times$ $L^2(X_2^k\times X_3^k)$ $\times$ $\dots$ $\times$
$L^2(X_{n-1}^k\times X_n^k)$ coincides with $S_k$. This implies
that the restriction of $\nph_{k+1}$ to $X_1^k\times\dots\times
X_n^k$ coincides (almost everywhere) with $\nph_k$. Hence, there
exists a function $\nph$ defined on $X_1\times\dots\times X_n$
which coincides with $\nph_k$ on $X_1^k\times\dots\times X_n^k$,
for each $k\in\bb{N}$. Since $\|\nph_k\|_{\infty} = \|S_k\|\leq
\|S\|$, we have that $\|\nph\|_{\infty}\leq \|S\|$. We have that
$S$ and $S_{\nph}$ coincide on the union of $L^2(X_1^k\times
X_2^k)$ $\times$ $L^2(X_2^k\times X_3^k)$ $\times$ $\dots$
$\times$ $L^2(X_{n-1}^k\times X_n^k)$, $k\in\bb{N}$, which is a
dense subset of $L^2(X_1\times X_2)$ $\times$ $L^2(X_2\times X_3)$
$\times$ $\dots$ $\times$ $L^2(X_{n-1}\times X_n)$. It follows
that $S = S_{\nph}$, and by the first part of the proof, $\|S\| =
\|\nph\|_{\infty}$. \prend

Let $(Y_1,\nu_1)$ and $(Y_2,\nu_2)$ be measure spaces. A subset
$E\subset Y_1\times Y_2$ is called marginally null \cite{a} if
$E\subset (A\times Y_2)\cup (Y_1\times B)$, $\nu_1(A)=\nu_2(B)=0$.
It is well-known that the projective tensor product
$L^2(Y_1)\hat{\otimes}L^2(Y_2)$ can be identified with a space of
complex-valued functions, defined marginally almost everywhere on
$Y_1\times Y_2$: the element $\sum_{i=1}^{\infty} f_i\otimes
g_i\in L^2(Y_1)\hat{\otimes}L^2(Y_2)$, where $f_i\in L^2(Y_1)$,
$g_i\in L^2(Y_2)$, $\sum_{i=1}^{\infty}\|f_i\|^2 < \infty$ and
$\sum_{i=1}^{\infty}\|g_i\|^2 < \infty$, is identified with the
function $h$ given by $h(x,y) = \sum_{i=1}^{\infty}f_i(x)g_i(y)$
(see e.g. \cite{a}).

Let
$$\Gamma(X_1,\dots,X_n) = L^2(X_1\times X_2)\odot\dots\odot L^2(X_{n-1}\times
X_n).$$ We identify the elements of $\Gamma(X_1,\dots,X_n)$ with
functions on $$X_1\times X_2\times X_2\times\dots\times
X_{n-1}\times X_{n-1}\times X_n$$ in the obvious fashion. We equip
$\Gamma(X_1,\dots,X_n)$ with two norms; one is the projective norm
$\|\cdot\|_{2,\wedge}$, where each of the $L^2$-spaces is equipped
with its $L^2$-norm, and the other is the Haagerup tensor norm
$\|\cdot\|_{\hh}$, where the $L^2$-spaces are given their opposite
operator space structure arising from the identification of
$L^2(X\times Y)$ with the class of Hilbert-Schmidt operators from
$L^2(X)$ into $L^2(Y)$ given by
\begin{equation}\label{hse}
(T_f\xi)(y) = \int_X f(x,y)\xi(x) dx, \ \ \ \ f\in L^2(X\times Y),
\xi\in L^2(X).
\end{equation}
For each $\nph\in L^{\infty}(X_1\times\dots\times
X_n)$, we consider the linearisation of the map $S_{\nph}$ from
Theorem \ref{th_mod} to a map defined on $\Gamma(X_1,\dots,X_n)$ and
taking values in $L^2(X_1\times X_n)$ and we denote it in the same
way. Thus, if $f_1\otimes\dots\otimes f_{n-1}$ is in
$\Gamma(X_1,\dots,X_n)$ then $S_{\nph}(f_1\otimes\dots\otimes
f_{n-1})(x_1,x_n)$ is equal to
$$\int_{X_2\times\dots\times X_{n-1}}
\nph(x_1,\dots,x_n)f_1(x_1,x_2)f_2(x_2,x_3)\dots
f_{n-1}(x_{n-1},x_n) dx_2\dots dx_{n-1}.$$ By Theorem
\ref{th_mod}, $S_{\nph}$ is bounded and $\|S_{\nph}\| =
\|\nph\|_{\infty}$. Hence it extends to a bounded map from
$(\Gamma(X_1,\dots,X_n),\|\cdot\|_{2,\wedge})$ into
$(L^2(X_1\times X_n),\|\cdot\|_2)$.

\begin{definition}\label{cmu}
Let $\nph\in L^{\infty}(X_1\times\dots\times X_n)$. We say that
$\nph$ is a Schur multiplier (relative to the measure spaces
$(X_1,\mu_1),\ldots (X_n,\mu_n)$) if there exists $C > 0$ such
that $\|S_{\nph}(\Phi)\|_{\op}\leq C \|\Phi\|_{\hh}$, for all
$\Phi\in\Gamma(X_1,\dots,X_n)$. The smallest constant $C$ with
this property will be denoted by $\|\nph\|_{\mm}$.
\end{definition}

Note that in the case where $n = 2$ and the measure spaces are
discrete, the definition above reduces to the definition of the
classical Schur multipliers. In the case of arbitrary measure
spaces and $n = 2$, we obtain the Schur multipliers studied by
Peller \cite{peller_two_dim} (see also \cite{spronk}).

We will present next a characterisation of the $n$-dimensional
Schur multipliers which generalises Grothendieck's and Peller's
characterisations. We will need the following generalisation of a
result of Smith \cite{smith}.

\begin{lemma}\label{l_gsmith}
Let $\cl E_i\subseteq B(H_i,H_{i+1})$, $i = 1,\dots,n$ be spaces
of operators and $\cl C\subseteq B(H_1)$, $\cl D\subseteq
B(H_{n+1})$ be C*-algebras with cyclic vectors. Assume that $\cl
E_1$ is a right $\cl C$-module and $\cl E_n$ is a left $\cl
D$-module. Let $\phi : \cl E_n\times\dots\times \cl E_1
\rightarrow B(H_1,H_{n+1})$ be a multilinear $\cl D, \cl C$-module
map (that is, $\phi(dy,\dots,xc) = d\phi(y,\dots,x)c$, whenever
$x\in\cl E_1$, $y\in\cl E_n$, $c\in\cl C$ and $d\in\cl D$) such
that the corresponding linear map from $\cl E_n\odot\dots\odot\cl
E_1$ into $B(H_1,H_{n+1})$ is bounded in the Haagerup norm. Then
$\phi$ is a completely bounded multilinear map.
\end{lemma}
\proof The proof is a straightforward generalisation of the
argument given by Smith \cite{smith}. We will denote by
$\tilde{\phi}$ the linear map from $\cl E_n\odot\dots\odot \cl
E_1$ into $\cl B(H_1,H_{n+1})$ defined by
$\tilde{\phi}(a_n\otimes\dots\otimes a_1) = \phi(a_n,\dots,a_1)$.
By the assumption of the lemma, it is bounded in the Haagerup norm
$\|\cdot\|_{\hh}$. Assume that $\|\tilde{\phi}\| = 1$. We will
show that $\|\tilde{\phi}\|_{\cb} = 1$. Suppose, to the contrary,
that $\|\tilde{\phi}\|_{\cb} > 1$. Then there exists $m\in\bb{N}$,
matrices $x^i = (x^i_{k j})\in M_m(\cl E_i)$, $i = 1,\dots,n$ and
column vectors $\xi_0 = (\xi_1,\dots,\xi_m)\in H_1^m$ and $\eta_0
= (\eta_1,\dots,\eta_m)\in H_{n+1}^m$ such that $\|\xi_0\| < 1$,
$\|\eta_0\| < 1$, all $\|x^i\| < 1$ and
\begin{equation}\label{eq_sm}
|(\phi^{(m)}(x^n,x^{n-1}\dots,x^1)\xi_0,\eta_0)|
> 1.
\end{equation}
If $\xi$ and $\eta$ are cyclic vectors for $\cl C$ and $\cl D$,
respectively, we may moreover assume that $\xi_i = a_i\xi$ and
$\eta_j = b_j\eta$, for some $a_i\in\cl C$ and $b_j\in \cl D$,
where $i,j=1,\dots,m$. Let $a = \sum_{i=1}^m a_i^*a_i$ and $b =
\sum_{j=1}^m b_j^*b_j$. Assume first that $a$ and $b$ are
invertible, and let $c_i = a_i a^{-1/2}$, $d_j = b_jb^{-1/2}$,
$\tilde{\xi} = a^{1/2}\xi$ and $\tilde{\eta} = b^{1/2}\eta$. Then
$\xi_i = c_i\tilde{\xi}$ and $\eta_j = d_j\tilde{\eta}$. Taking
into account (\ref{rfl}), the left hand side of (\ref{eq_sm})
becomes
\begin{eqnarray}\label{fromr}
& & \left|\sum_{i,j=1}^m(\phi^{(m)}(x^n,x^{n-1},
\dots,x^1)_{ji}c_i\tilde{\xi},d_j\tilde{\eta})\right|\nonumber\\
& = & \left|\sum_{k_1,\dots,k_{n-1} = 1}^m\sum_{i,j = 1}^m
(\phi(d_j^*x_{j k_{n-1}}^n,x_{k_{n-1} k_{n-2}}^{n-1},\dots,x_{k_1
i}^1c_i)\tilde{\xi},\tilde{\eta})\right|\nonumber\\
& = & \left|\sum_{k_1,\dots,k_{n-1}=1}^m
\left(\phi\left(\sum_{j=1}^m d_j^*x_{j k_{n-1}}^n, x^{n-1}_{k_{n-1}
k_{n-2}},\dots,\sum_{i=1}^m
x_{k_{1},i}^1c_i\right)\tilde{\xi},\tilde{\eta}\right)\right|\nonumber\\
& \leq & \left\|\sum_{k_1,\dots,k_{n-1}=1}^m \phi\left(\sum_{j=1}^m
d_j^*x_{j k_{n-1}}^n, x^{n-1}_{k_{n-1} k_{n-2}},\dots,\sum_{i=1}^m
x_{k_{1},i}^1c_i\right)\right\|\|\tilde{\xi}\|\|\tilde{\eta}\|.
\end{eqnarray}
We have that
$$\|\tilde{\xi}\| = (a^{1/2}\xi,a^{1/2}\xi) = (a\xi,\xi) =
\sum_{k=1}^n\|a_i\xi\|^2 = \sum_{k=1}^n\|\xi_i\|^2  = \|\xi_0\|
\leq 1,$$ and similarly $\|\tilde{\eta}\|\leq 1$. Set $d^* =
(d_j^*)\in M_{1,m}(\cl D)$, $c = (c_i)\in M_{m,1}(\cl C)$, $u =
d^*x^n \in M_{1,m}(\cl E_n)$ and $v = x^1c\in M_{m,1}(\cl E_1)$.
It follows from (\ref{odot}) and (\ref{odot2}) that
\begin{eqnarray}\label{onemo}
& & \left\|\sum_{k_1,\dots,k_{n-1}=1}^m \phi\left(\sum_{j=1}^m
d_j^*x_{j k_{n-1}}, x^{n-1}_{k_{n-1} k_{n-2}},\dots,\sum_{i=1}^m
x_{k_{1},i}c_i\right)\right\|\nonumber\\
& = & \left\|\sum_{k_1,\dots,k_{n-1}=1}^m
\phi\left(u_{k_{n-1}},x^{n-1}_{k_{n-1} k_{n-2}},\dots,v_{k_1}\right)\right\|\nonumber\\
& = & \left\| \tilde{\phi}\left( \sum_{k_1,\dots,k_{n-1}=1}^m
u_{k_{n-1}}\otimes x^{n-1}_{k_{n-1} k_{n-2}}\otimes
\dots\otimes v_{k_1}\right)\right\|\nonumber\\
& \leq & \left\|\sum_{k_1,\dots,k_{n-1}=1}^m u_{k_{n-1}}\otimes
x^{n-1}_{k_{n-1} k_{n-2}}\otimes
\dots\otimes v_{k_1}\right\|_{\hh}\nonumber\\
& = & \|u\odot x^{n-1}\odot\dots\odot x^2\odot
v\|_{\hh}\nonumber\\ & \leq &
\|d^*\|\|x^n\|\|x^{n-1}\|\dots\|x^2\|\|x^1\|\|c\|.
\end{eqnarray}
We have that
$$\|d^*\| = \left\|\sum_{j=1}^m d_j^*d_j\right\|^{1/2} = \|I\| = 1$$
and, similarly, $\|c\| = 1$. It follows from (\ref{fromr}) and
(\ref{onemo}) that
$$|(\phi^{(m)}(x^n,x^{n-1}\dots,x^1)\xi_0,\eta_0)| \leq 1,$$ which
contradicts (\ref{eq_sm}).

In the case $a$ or $b$ is not invertible, one can again follow
\cite{smith} and, for each $i$, consider the matrix $\hat{x}^i\in
M_{m+1}(\cl E_{i})$ which has the matrix $x^i$ in its upper left
corner and zeros in the last row and column. The vectors $\xi_0$
and $\eta_0$ are replaced with $\hat{\xi}_0 =
(\xi_1,\dots,\xi_m,\xi_{m+1})$ and $\hat{\eta}_0 =
(\eta_1,\dots,\eta_m,\eta_{m+1})$, where $\xi_{m+1} = \epsilon\xi$
and $\eta_{m+1} = \epsilon\eta$, respectively, for $\epsilon$
small enough so that the norms of these vectors remain less than
one. Letting $a_{n+1} = b_{n+1} = \epsilon I$, we have that
$a_i\xi = \xi_i$ and $b_i\eta = \eta_i$ for each $i =
1,\dots,m+1$. Finally,
$$(\phi^{(m)}(x^n,x^{n-1}\dots,x^1)\xi_0,\eta_0) =
(\phi^{(m+1)}(\hat{x}^n,\hat{x}^{n-1}\dots,\hat{x}^1)
\hat{\xi}_0,\hat{\eta}_0)$$ and the proof proceeds as before.
\prend

The main result of this section is the following

\begin{theorem}\label{th_g3}
Let $\nph\in L^{\infty}(X_1\times\dots\times X_n)$. The following
are equivalent:

(i) \ $\nph$ is a Schur multiplier and $\|\nph\|_{\mm} < 1$;

(ii) there exist essentially bounded functions $a_1 :
X_1\rightarrow M_{\infty,1}$, $a_n : X_n\rightarrow M_{1,\infty}$
and $a_i : X_i\rightarrow M_{\infty}$, $i = 2,\dots,n-1$, such
that, for almost all $x_1,\dots,x_n$ we have
$$\nph(x_1,\dots,x_n) = a_n(x_n)a_{n-1}(x_{n-1})\dots a_1(x_1) \ \mbox{ and } \
\esssup_{x_i\in X_i}\prod_{i=1}^n \|a_i(x_i)\| < 1.$$
\end{theorem}
\proof (i)$\Rightarrow$(ii) Let $\nph\in
L^{\infty}(X_1\times\dots\times X_n)$ be a Schur multiplier with
$\|\nph\|_{\mm} < 1$. Then the map $S_{\nph}$ induces a map,
denoted in the same way, from $L^2(X_1\times X_2)\times\dots\times
L^2(X_{n-1}\times X_n)$ into $L^2(X_1\times X_n)$. Let $H_i =
L^2(X_i)$, $\cl D_i = \{M_{\psi} : \psi \in L^{\infty}(X_i)\}$, $i
= 1,\dots,n$, and
$$\hat{S}_{\nph} : \cl C_2(H_1,H_2)\times\dots\times \cl C_2(H_{n-1},H_n)
\rightarrow \cl C_2(H_1,H_n)$$ be the map defined by
$\hat{S}_{\nph} (T_{f_1},\dots,T_{f_n}) =
T_{S_{\nph}(f_1,\dots,f_n)}$. Since $\nph$ is a Schur multiplier,
the linearisation of the map $\hat{S}_{\nph}$ from $(\cl
C_2(H_1,H_2)\odot\dots\odot \cl C_2(H_{n-1},H_n),
\|\cdot\|_{\hh})$ into $(\cl C_2(H_1,H_n),\|\cdot\|_{\op})$ is
bounded. (Here each of the operator spaces $\cl C_2(H_i,H_{i+1})$
is given its opposite operator space structure arising from the
inclusion $\cl C_2(H_i,H_{i+1})\subseteq\cl B(H_i,H_{i+1})$.) If
$a_i\in L^{\infty}(X_i)$, $i = 1,\dots,n$, then
\begin{eqnarray}\label{modu}
\hat{S}_{\nph}(T_{f_1}M_{a_1},T_{f_2}M_{a_2},\dots,M_{a_n}T_{f_n}M_{a_{n-1}})
& = & \hat{S}_{\nph}(T_{f_1a_1},T_{f_2a_2},\dots,T_{a_nf_na_{n-1}})\nonumber\\
& = & T_{S_{\nph}(f_1a_1,f_2a_2,\dots,a_nf_na_{n-1})}\\ & =
& T_{a_nS_{\nph}(a_2f_1,a_3f_2,\dots,a_{n-1}f_{n-2},f_n)a_1}\nonumber\\
& = &
M_{a_n}\hat{S}_{\nph}(M_{a_2}T_{f_1},\dots,T_{f_n})M_{a_1}\nonumber.
\end{eqnarray}

By continuity, the map $\hat{S}_{\nph}$ has an extension (denoted
in the same way)
$$\hat{S}_{\nph} : \cl K(H_1,H_2)\otimes_{\hh}\dots\otimes_{\hh} \cl K(H_{n-1},H_n)
\rightarrow \cl K(H_1,H_n)$$ to a  map with norm less than one,
where the spaces $\cl K(H_i,H_{i+1})$ are equipped with the operator
space structure opposite to their natural operator space structure.
It follows from (\ref{modu}) that the map
$$\check{S}_{\nph} :  \cl K(H_{n-1},H_n)\otimes_{\hh}\dots\otimes_{\hh} \cl K(H_1,H_2)
\rightarrow \cl K(H_1,H_n)$$ given by
$$\check{S}_{\nph}(T_{n-1}\otimes\dots\otimes T_1) =
\hat{S}_{\nph}(T_1\otimes\dots\otimes T_{n-1})$$ is modular and
bounded when the spaces $\cl K(H_i,H_{i+1})$ are given their
natural operator space structure. By Lemma \ref{l_gsmith},
$\check{S}_{\nph}$ is completely bounded. It follows that the
second dual
$$\check{S}_{\nph}^{**} : \cl B(H_{n-1},H_n)\otimes_{\sigma
h}\dots\otimes_{\sigma h}\cl B(H_1,H_2) \rightarrow \cl
B(H_1,H_n)$$ is a weak* continuous map with c.b. norm less than
one, which extends the map $\check{S}_{\nph}$. (Here
$\otimes_{\sigma h}$ denotes the normal Haagerup tensor product,
see e.g. \cite{blm}.)

Denote by $\tilde{S}_{\nph}$ the corresponding multilinear map
$$\tilde{S}_{\nph} : \cl B(H_{n-1},H_n)\times\dots\times \cl B(H_1,H_2)
\rightarrow \cl B(H_1,H_n).$$ The map $\tilde{S}_{\nph}$ is
separately weak* continuous and hence modular.

A modification of Corollary 5.9 of \cite{cs} now implies that
there exist bounded linear operators $V_1 : H_1\rightarrow
H_1^{\infty}$, $V_n : H_n^{\infty}\rightarrow H_n$ and $V_i :
H_i^{\infty}\rightarrow H_i^{\infty}$, $i = 2,\dots,n-1$, such
that the entries of $V_i$ belong to $\cl D_i$ and
$$\tilde{S}_{\nph}(T_{n-1},\dots,T_1) = V_n (T_{n-1}\otimes I) V_{n-1}
(T_{n-2}\otimes I) \dots (T_1\otimes I) V_1.$$ Moreover, the
operators $V_i$ can be chosen so that $\prod_{i=1}^n \|V_i\| < 1$.
Let $V_1 = (M_{a_1^1},M_{a_2^1},$ $\dots$ $)^{\ttt}$, $V_i =
(M_{a^i_{kl}})$ and $V_n = (M_{a^n_1},M_{a^n_2},\dots)$, for some
$a_1 = (a_1^1,a_2^1,$ $\dots$ $)^{\ttt}\in
L^{\infty}(X_1,M_{1,\infty})$, $a_n = (a^n_1,a^n_2,\dots)\in
L^{\infty}(X_n,M_{1,\infty})$ and $a_i = (a^i_{kl})\in
L^{\infty}(X_i,M_{\infty})$, $i = 2,\dots,n-1$. Moreover,
$$\esssup_{x_i\in X_i}\prod_{i=1}^n \|a_i(x_i)\| = \prod_{i=1}^n \|V_i\| <
1.$$

If $\xi\in L^2(X)$ and $\eta\in L^2(Y)$ denote by $\xi\otimes\eta$
the function on $X\times Y$ given by $(\xi\otimes\eta) (x,y) =
\xi(x)\eta(y)$; this function gives rise by (\ref{hse}) to a rank
one operator $T_{\xi\otimes\eta}$. Fix $\xi_i,\eta_i\in H_i$, $i =
1,\dots,n$. Then
\begin{eqnarray*}
& &
\tilde{S}_{\nph}(T_{\xi_{n-1}\otimes\eta_n},\dots,T_{\xi_1\otimes\eta_2})(\eta_1)
= V_n(T_{\xi_{n-1}\otimes\eta_n}\otimes I)\dots
(T_{\xi_{1}\otimes\eta_2}\otimes I) V_1 (\eta_1)\\
& = & V_n(T_{\xi_{n-1}\otimes\eta_n}\otimes I)\dots V_2
(T_{\xi_{1}\otimes\eta_2}\otimes I) (a^1_{k_1}\eta_1)_{k_1}\\
& = & V_n(T_{\xi_{n-1}\otimes\eta_n}\otimes I)\dots V_2
((\int_{X_1}a^1_{k_1}(x_1)\xi_1(x_1)\eta_1(x_1)dx_1)\eta_2)_{k_1}\\
& = & V_n\dots (T_{\xi_{2}\otimes\eta_3}\otimes I) ((\sum_{k_1=
1}^{\infty} \int_{X_1}a^1_{k_1}(x_1)\xi_1(x_1)\eta_1(x_1)dx_1)
a^2_{k_2,k_1}\eta_2)_{k_2}\\
& = & V_n\dots V_3((\sum_{k_1= 1}^{\infty} \int_{X_1\times
X_2}a^2_{k_2,k_1}(x_2)a^1_{k_1}(x_1)(\xi_1\eta_1)(x_1)(\xi_2\eta_2)(x_2)dx_1
dx_2)\eta_3)_{k_2}\\
\end{eqnarray*}
\begin{eqnarray*}
& = &
\dots\dots\dots\dots\dots\dots\dots\dots\dots\dots\dots\dots\dots\dots\dots\dots\dots\dots\dots\dots\dots\dots\\
& = & \sum_{k_n = 1}^{\infty}(\int_{X_1\times\dots\times
X_{n-1}}\sum_{k_1,\dots,k_{n-1} = 1}^{\infty}
a^{n-1}_{k_{n-1},k_{n-2}}(x_{n-1})\dots a^1_{k_1}(x_1)\times \\
& \times & \xi_1(x_1)\eta_1(x_1)\dots\xi_{n-1}(x_{n-1}))dx_1\dots
dx_{n-1}) M_{a^n_{k_n}}\eta_n.
\end{eqnarray*}
Thus,
\begin{eqnarray*}
& &
\tilde{S}_{\nph}(T_{\xi_{n-1}\otimes\eta_n},\dots,T_{\xi_1\otimes\eta_2})(\eta_1)(x_n)\\
& = & (\int_{X_1\times\dots\times X_{n-1}}\sum_{k_1,\dots,k_{n} =
1}^{\infty} a^n_{k_n}(x_n)a^{n-1}_{k_{n-1},k_{n-2}}(x_{n-1})\dots
a^1_{k_1}(x_1) \times \\ & \times &
\xi_1(x_1)\eta_1(x_1)\dots\xi_{n-1}(x_{n-1})dx_1\dots dx_{n-1})
\eta_n(x_n).
\end{eqnarray*}
On the other hand,
\begin{eqnarray*}
& &
\tilde{S}_{\nph}(T_{\xi_{n-1}\otimes\eta_n},\dots,T_{\xi_1\otimes\eta_2})(\eta_1)(x_n)
=
T_{S_{\nph}(\xi_1\otimes\eta_2,\dots,\xi_{n-1}\otimes\eta_n)}(\eta_1)(x_n)\\
& = & (\int_{X_1\times\dots\times
X_{n-1}}\nph(x_1,\dots,x_{n-1},x_n)\\
& \times & \xi_1(x_1)\eta_1(x_1)\dots\xi_{n-1}(x_{n-1})dx_1\dots
dx_{n-1}) \eta_n(x_n).
\end{eqnarray*} It follows that
$$\nph(x_1,\dots,x_n) = a_n(x_n)a_{n-1}(x_{n-1})\dots a_1(x_1),$$
for almost all $x_1,\dots,x_n$.

(ii)$\Rightarrow$(i) Assume that $\nph$ is given as in (ii), where
$a_1 = (a^1_1,a^1_2,\dots)^{\ttt}\in L^{\infty}(X_1,$
$M_{\infty,1})$, $a_n = (a^n_1,a^n_2,\dots) \in
L^{\infty}(X_n,M_{1,\infty})$ and $a_i = (a^i_{kl})\in
L^{\infty}(X_i,$ $M_{\infty})$, $i = 2,\dots,n-1$. Let $V_1 :
H_1\rightarrow H_1^{\infty}$ be the operator corresponding to the
column matrix $V_1 = (M_{a_1^{1}},M_{a_2^{1}},\dots)^{\ttt} :
H_1\rightarrow H_1^{\infty}$, $V_n : H_n^{\infty}\rightarrow H_n$ be
the operator corresponding to the row matrix $V_n =
(M_{a^n_{1}},M_{a^n_{2}},\dots)$ and $V_i : H_i^{\infty}\rightarrow
H_i^{\infty}$ be the operator corresponding to the matrix $V_i =
(M_{a^i_{kl}})$, $i = 2,\dots,n-1$. Then $\prod_{i=1}^n \|V_i\| <
1$. It follows from the first part of the proof that
$$\tilde{S}_{\nph}(T_{\xi_{n-1}\otimes\eta_n},\dots,T_{\xi_1\otimes\eta_2})
= V_n(T_{\xi_{n-1}\otimes\eta_n}\otimes I)\dots
(T_{\xi_{1}\otimes\eta_2}\otimes I) V_1,$$
for all $\xi_1\in H_1$, $\eta_n\in H_n$ and $\xi_i,\eta_i\in H_i$,
$i = 2,\dots,n-1$. Since the operator norm is dominated by the
Hilbert-Schmidt norm, we conclude that
$$\tilde{S}_{\nph}(T_{f_{n-1}},\dots,T_{f_1})\\
= V_n(T_{f_{n-1}}\otimes I)\dots (T_{f_{1}}\otimes I) V_1,$$
for all $f_i\in L^2(X_i\times X_{i+1})$, $i = 1,\dots,n-1$.

Let $$F = F_1\odot \dots\odot F_{n-1}\in L^2(X_1\times
X_2)\odot\dots\odot L^2(X_{n-1}\times X_n),$$ where $F_1\in
M_{1,\infty}(L^2(X_1\times X_2))$, $F_{n-1}\in
M_{\infty,1}(L^2(X_{n-1}\times X_n))$ and $F_i\in
M_{\infty}(L^2(X_i\times X_{i+1}))$, $i = 2,\dots,n-2$. Lemma
\ref{l_opmul} implies that
$$T_{S_{\nph}(F)} = V_n(T_{F_{n-1}}\otimes I)\dots (T_{F_{1}}\otimes I)
V_1,$$ where $T_{F_i} = (T_{f^i_{lk}})_{k,l}$ whenever $F_i =
(f^i_{kl})_{k,l}$. It follows that
$$\|T_{S_{\nph}(F)}\|_{\op}\leq
\prod_{i=1}^{n-1} \|F_i^t\|_{\op}\prod_{i=1}^n \|V_i\|.$$ Taking
infimum with respect to all representations of $F$, we conclude
that $\|T_{S_{\nph}(F)}\|_{\op}\leq  \|F\|_{\hh}\prod_{i=1}^n
\|V_i\|$ and so $\|\nph\|_{\mm} < 1$. \prend

\noindent {\bf Remark} The space of all functions
$\varphi(x_1,\ldots, x_n)$ satisfying condition (ii) of
Theorem~\ref{th_g3} can be identified with the extended Haagerup
tensor product $L^{\infty}(X_1)\otimes_{eh}$ $L^{\infty}(X_2)$
$\otimes_{eh}$ $\ldots$ $\otimes_{eh}$ $L^{\infty}(X_n)$.

\medskip

The next proposition relates our approach with a recent work of
Peller \cite{peller} on multiple operator integrals. For some
fixed spectral measures, Peller defines a multiple operator
integral $I_{\nph}(T_1,\dots,T_{n-1})$ of a function $\nph$ and
$(n-1)$-tuple of operators $(T_1,\dots,T_{n-1})$, and shows that
if $\nph$ belongs to the integral projective tensor product of the
corresponding $L^{\infty}$-spaces, then
$I_{\nph}(T_1,\dots,T_{n-1})$ is well-defined and, moreover,
$$\|I_{\nph}(T_1,\dots,T_{n-1})\|_{\op}\leq\|\nph\|_{i}\|T_1\|_{\op}\dots\|T_{n-1}\|_{\op}.$$
Recall that the integral projective tensor product
$L^{\infty}(X_1)\hat{\otimes}_{i}\dots\hat{\otimes}_{i}L^{\infty}(X_n)$
is the space of all functions $\nph$ for which there exists a
measure space $(\cl T,\nu)$ and measurable functions $g_i$ on
$X_i\times\cl T$ such that
\begin{equation}\label{iptp}
\nph(x_1,\dots,x_n) = \int_{\cl T}g_1(x_1,t)\dots
g_n(x_n,t)d\nu(t),
\end{equation}
for almost all $x_1\dots,x_n$, where
$$\int_{\cl T}\|g_1(\cdot,t)\|_{\infty}\dots\|g_n(\cdot,t)\|_{\infty} d\nu(t)
< \infty.$$ The integral projective norm $\|\nph\|_i$ of $\nph$ is
the infimum of the above expressions over all representations of
$\nph$ of the form (\ref{iptp}). It was proved by Peller in
\cite{peller_two_dim} that in the case where $n=2$ the integral
projective tensor product $L^{\infty}(X_1)\hat\otimes_i
L^{\infty}(X_2)$ coinsides with the set of all Schur mulipliers.
The next proposition shows that for $n>2$ the integral projective
tensor product consists of multipliers. We do not know whether it
coincides with the space of all Schur multipliers.

\begin{proposition}\label{p_peller}
Let $\nph\in L^{\infty}(X_1)\hat{\otimes}_{i}\dots
\hat{\otimes}_{i}L^{\infty}(X_n)$. Then $\nph$ is a Schur
multiplier and $\|\nph\|_{\mm}\leq\|\nph\|_{i}$.
\end{proposition}
\proof Suppose that
$$\nph(x_1,\dots,x_n) =
\int_{\cl T}g_1(x_1,t)\dots g_n(x_n,t)d\nu(t),$$ for almost all
$x_1\dots,x_n$, where $(\cl T,\nu)$ is a measure space, $g_i$ is a
measurable function on $X_i\times\cl T$, $i = 1,\dots,n$, such
that
$$\int_{\cl T}\|g_1(\cdot,t)\|_{\infty}\dots\|g_n(\cdot,t)\|_{\infty} d\nu(t)
< \infty.$$ Let $F = F_1\odot\dots\odot F_{n-1}$, where $F_1\in
M_{1,k_1}(L^2(X_1\times X_2))$, $F_{n-1}\in M_{k_{n-2},1}$
$(L^2(X_{n-1}$ $\times X_n))$ and $F_i\in
M_{k_{i-1},k_i}(L^2(X_i\times X_{i+1}))$, $i = 2,\dots,n-2$, and
$\tilde{F}(x_1,x_2\ldots,x_n)=F(x_1,x_2,x_2,x_3,\ldots,x_n)$.
Denoting by $M_{g_i(\cdot,t)}$ the multiplication operator by the
function $g_i(\cdot,t)$, and by $M_{g_i(\cdot,t)}\otimes I$ the
ampliation of $M_{g_i(\cdot,t)}$ of multiplicity $k_i$, we have
\begin{eqnarray*}
\|S_{\nph}(F)\|_{\op} & = & \|\int_{X_2\times\dots\times
X_{n-1}}\nph \tilde F dx_2\dots dx_{n-1}\|_{\op}\\ & = &
\|\int_{X_2\times\dots\times X_{n-1}} \left(\int_{\cl
T}g_1(x_1,t)\dots g_n(x_n,t)dt\right)\tilde F dx_2\dots
dx_{n-1}\|_{\op}\\
& = & \|\int_{\cl T} \left(\int_{X_2\times\dots\times
X_{n-1}}g_1(x_1,t)\dots g_n(x_n,t)dx_2\dots dx_{n-1}\right)\tilde F dt\|_{\op}\\
& = & \|\int_{\cl T} (\int_{X_2\times\dots\times X_{n-1}}
M_{g_1(\cdot,t)}F_1(M_{g_2(\cdot,t)}\otimes I) (x_1,x_2) \odot
\dots\\ & \odot & F_{n-1}M_{g_n(\cdot,t)}(x_{n-1},x_n)dx_2\dots
dx_{n-1})dt\|_{\op}\\
\end{eqnarray*}
\begin{eqnarray*}
& \leq & \int_{\cl T} \|\int_{X_2\times\dots\times X_{n-1}}
M_{g_1(\cdot,t)}F_1(M_{g_2(\cdot,t)}\otimes I) (x_1,x_2) \odot
\dots\\ & \odot & F_{n-1}M_{g_n(\cdot,t)}(x_{n-1},x_n)dx_2\dots
dx_{n-1}\|_{\op}dt\\
& \leq & \int_{\cl T}
\|M_{g_1(\cdot,t)}\|\|F_1\|_{\op}^{o}\|M_{g_2(\cdot,t)}\|
\dots\|F_{n-1}\|_{\op}^o\|M_{g_n(\cdot,t)}\|dt\\
& \leq & \|\nph\|_{i}\|F_1\|_{\op}^{o}\dots\|F_{n-1}\|_{\op}^{o}.
\end{eqnarray*}
where $\|\cdot\|_{\op}^o$ is the opposite operator norm (see Section
\ref{prel}). The claim follows by taking infimum over all
representations $F = F_1\odot\dots\odot F_{n-1}$. \prend

\begin{corollary}\label{ini}
$L^{\infty}(X_1)\hat\otimes_i\ldots\hat\otimes_i
L^{\infty}(X_n)\subseteq
L^{\infty}(X_1)\otimes_{eh}\ldots\otimes_{eh} L^{\infty}(X_n).$
\end{corollary}

In the case where $n = 2$, it follows by Peller's characterisation
of Schur multipliers \cite{peller_two_dim} that there is an
equality in the inclusion of Corollary \ref{ini}. We do not know
whether equality holds in the general case.

We finally point out another interesting open question, namely the
one of characterising the class of multipliers defined by using
the projective tensor norm instead of the Haagerup tensor norm in
(\ref{def}); equivalently, the class of multipliers obtained after
replacing (\ref{def}) with the weaker condition
$$\|S_{\psi}(f_1\otimes\ldots\otimes f_n)\|_{\op}\leq
C\|f_1\|_{\op}\ldots\|f_n\|_{\op}\text{ for all } f_i\in L^2(X_i),
i=1,\ldots,n.$$

\section{Multidimensional operator multipliers: the definition}\label{s_nonc}

In this section we generalise the notion of operator multipliers
given by Kissin and Shulman \cite{ks} to the multidimensional
case.

We recall the mapping $\theta_{K_1,K_2} : K_1\otimes
K_2\rightarrow \cl C_2(K_1^{\dd},K_2)$, where $K_1$ and $K_2$ are
Hilbert spaces, which is the unitary operator between the Hilbert
spaces $K_1\otimes K_2$ and $\cl C_2(K_1^{\dd},K_2)$ given on
elementary tensors by
$$\theta_{K_1,K_2}(\xi_1\otimes\xi_2)(\eta_1^{\dd}) =
(\xi_1,\eta_1)\xi_2.$$ Note that there is a natural identification
of $(K_1\otimes K_2)^{\dd}$ and $K_1^{\dd}\otimes K_2^{\dd}$. It
follows that $\cl C_2(K_1^{\dd},K_2)^{\dd}$ can be identified with
$\cl C_2(K_1,K_2^{\dd}) = \cl C_2((K_1^{\dd})^{\dd},K_2^{\dd})$;
we have that $\theta_{K_1^{\dd},K_2^{\dd}}(\xi^{\dd}) =
\theta_{K_1,K_2}(\xi)^{\dd}$.

Let $H_1,\dots,H_n$ be Hilbert spaces and $H = H_1\otimes\dots
H_n$. For any permutation $\pi$ of $\{1,\dots,n\}$, we will
identify $H$ with the tensor product $H_{\pi(1)}\otimes\dots
H_{\pi(n)}$ without explicitly mentioning this. The symbol
$\xi_{j_1,\dots,j_k}$ will denote an element of
$H_{j_1}\otimes\dots H_{j_k}$.

We define a Hilbert space $HS(H_1,\dots,H_n)$, isometrically
isomorphic to $H$. Let $HS(H_1,H_2) = \cl C_2(H_1^{\dd},H_2)$. In
the case where $n$ is even, we let by induction
$$HS(H_1,\dots,H_n) = \cl C_2(HS(H_2,H_3)^{\dd},
HS(H_1,H_4,\dots,H_n)),$$ and let
$$\theta_{H_1,\dots,H_n} : H \rightarrow HS(H_1,\dots,H_n)$$
be given by $$\theta_{H_1,\dots,H_n}(\xi_{2,3}\otimes
\xi)=\theta_{HS(H_2,H_3),HS(H_1,H_4,\dots,H_n)}(\theta_{H_2,H_3}(\xi_{2,3})\otimes
\theta_{H_1,H_4,\dots,H_n}(\xi)),$$ where $\xi\in H_1\otimes
H_4\otimes\dots\otimes H_n$. In particular, we have that
$$\theta_{H_1,\dots,H_n}(\xi_{2,3}\otimes
\xi)\theta_{H_2,H_3}(\eta_{2,3})^{\dd}=(\theta_{H_2,H_3}(\xi_{2,3}),
\theta_{H_2,H_3}(\eta_{2,3})) \theta_{H_1,H_4,\dots,H_n}(\xi).$$
In the case where $n$ is odd, we let
$$HS(H_1,\dots,H_n) = HS(\bb{C},H_1,\dots,H_n).$$
If $K$ is a Hilbert space, we will identify $\cl
C_2(\bb{C}^{\dd},K)$ with $K$ via the map $S\rightarrow
S(1^{\dd})$. Thus, $HS(H_1,\dots,H_n)$ can, in the case of odd
$n$, be defined inductively by letting $HS(H_1) = H_1$ and
$$HS(H_1,\dots,H_n) = \cl C_2(HS(H_1,H_2)^{\dd},HS(H_3,\dots,H_n)).$$
The isomorphism $\theta_{H_1,\dots,H_n}$ is in this case given by
$$\theta_{H_1,\dots,H_n}(\xi) =
\theta_{\bb{C},H_1,\dots,H_n}(1\otimes\xi).$$ We will usually omit
the subscripts and write simply $\theta$, when the corresponding
Hilbert spaces are understood.

\begin{lemma}\label{l_compos}
(i) Assume $n$ is even. Let $\xi\in H$ be of the form $\xi =
\xi_{1,2}\otimes\dots\otimes\xi_{n-1,n}$. If $\eta_{i,i+1}\in
H_i\otimes H_{i+1}$ ($i$ even) then
$$\theta(\xi)(\theta(\eta_{2,3}^{\dd}))\dots(\theta(\eta_{n-2,n-1}^{\dd}))
=\theta(\xi_{n-1,n})\theta(\eta_{n-2,n-1}^{\dd})\dots\theta(\xi_{3,4})\theta(\eta_{2,3}^{\dd})\theta(\xi_{1,2}).$$

(ii) Assume $n$ is odd. Let $\xi\in H$ be of the form $\xi =
\xi_{1}\otimes\xi_{2,3}\dots\otimes\xi_{n-1,n}$. If
$\eta_{i,i+1}\in H_i\otimes H_{i+1}$ ($i$ odd) then
$$\theta(\xi)(\theta(\eta_{1,2}^{\dd}))(\theta(\eta_{3,4}^{\dd}))\dots(\theta(\eta_{n-2,n-1}^{\dd}))
=\theta(\xi_{n-1,n})\theta(\eta_{n-2,n-1}^{\dd})\dots\theta(\eta_{1,2}^{\dd})(\xi_{1}).$$
\end{lemma}
\proof (i) Assume first that $\xi_{i-1,i} = \xi_{i-1}\otimes\xi_i$
and $\eta_{i,i+1} = \eta_i\otimes\eta_{i+1}$ ($i$ even). Fix
$\eta_1^{\dd}\in H_1^{\dd}$. The image of $\eta_1^{\dd}$ under the
operator on the right hand side of the identity in (i) is
$$(\xi_1,\eta_1)(\xi_2,\eta_2)\dots(\xi_{n-1},\eta_{n-1})\xi_n.$$
On the other hand, the image of $\eta_1^{\dd}$ under the operator
on the left hand side is
\begin{eqnarray*}
& & (\theta_{H_2,H_3}(\xi_2\otimes
\xi_3),\theta_{H_2,H_3}(\eta_2\otimes\eta_3))\\
& \times & \theta_{H_1,H_4,\dots,H_n}(\xi_1\otimes\xi_4
\otimes\dots\otimes\xi_n)(\theta(\eta_{4,5})^{\dd})\dots
(\theta(\eta_{n-2,n-1})^{\dd})(\eta_1^{\dd})\\
& = & (\xi_2,\eta_2)(\xi_3,\eta_3)\\ & \times &
\theta_{H_1,H_4,\dots,H_n}(\xi_1\otimes\xi_4
\otimes\dots\otimes\xi_n)(\theta(\eta_{4,5})^{\dd})\dots
(\theta(\eta_{n-2,n-1})^{\dd})(\eta_1^{\dd}).
\end{eqnarray*}
By induction, (i) holds in the case of elementary tensors.

By linearity, (i) holds for finite sums of elementary tensors.
Using continuity arguments and the fact that the operator norm is
dominated by the Hilbert-Schmidt norm, one can easily prove that
$(i)$ holds for genral $\xi$ and $\eta_{i,i+1}$. \prend

We define a representation $\sigma_{H}$ of $B(H)$ on
$HS(H_1,\dots,H_n)$ by letting
$$\sigma_{H}(A)\theta(\xi) = \theta(A\xi);$$ clearly,
$\sigma_{H}$ is unitarily equivalent to the identity
representation of $B(H)$. If $H_1,\dots,H_n$ are clear from the
context we will simply write $\sigma$ in the place of
$\sigma_{H}$. If $\cl A_1,\dots,\cl A_n$ are C*-algebras,
 $\pi_1,\dots,\pi_n$ corresponding representations on
$H_1,\dots,H_n$, and $\pi=\pi_1\otimes\dots\pi_n$ we let
$$\sigma_{\pi} = \sigma_{H}\circ
\pi \ ;$$ thus, $\sigma_{\pi}$ is a representation of $\cl
A_1\otimes\dots\otimes\cl A_n$ on $HS(H_1,\dots,H_n)$, unitarily
equivalent to $\pi$.

\begin{lemma}\label{l_compos2}
Let $A_i\in B(H_i)$, $i = 1,\dots,n$, and $A =
A_1\otimes\dots\otimes A_n$.

(i) \ Assume $n$ is even. Let $\xi_{i-1,i}\in H_{i-1}\otimes H_i$,
$\eta_{i,i+1}\in H_i\otimes H_{i+1}$ ($i$ even). If $\xi =
\xi_{1,2}\otimes\dots\xi_{n-1,n}$ then
\begin{eqnarray*}
& & \sigma(A)(\theta(\xi))(\theta(\eta_{2,3}^{\dd}))\dots
(\theta(\eta_{n-2,n-1}^{\dd}))\\
& = &
A_n\theta(\xi_{n-1,n})A_{n-1}^{\dd}\theta(\eta_{n-2,n-1})^{\dd}A_{n-2}\dots
A_2\theta(\xi_{1,2})A_1^{\dd}\\
& = & A_n\theta(\xi)(\theta((A_2^*\otimes
A_3^*(\eta_{2,3}))^{\dd}))\dots (\theta((A_{n-2}^*\otimes
A_{n-1}^*(\eta_{n-2,n-1}))^{\dd}))A_1^{\dd}.
\end{eqnarray*}

(ii) Assume $n$ is odd. Let $\xi_1\in H_1$, $\xi_{i-1,i}\in
H_{i-1}\otimes H_i$, $\eta_{i,i+1}\in H_i\otimes H_{i+1}$ ($i$
odd). If $\xi = \xi_1\otimes\xi_{2,3}\otimes\dots\xi_{n-1,n}$ then
\begin{eqnarray*}
& & \sigma(A)(\theta(\xi))(\theta(\eta_{1,2}^{\dd}))\dots
(\theta(\eta_{n-2,n-1}^{\dd}))\\
& = &
A_n\theta(\xi_{n-1,n})A_{n-1}^{\dd}\theta(\eta_{n-2,n-1})^{\dd}A_{n-2}\dots
A_2^{\dd}\theta(\eta_{1,2}^{\dd})(A_1\xi_1)\\
& = & A_n\theta(\xi)(\theta((A_1^*\otimes
A_2^*(\eta_{1,2}))^{\dd}))\dots (\theta((A_{n-2}^*\otimes
A_{n-1}^*(\eta_{n-2,n-1}))^{\dd})).
\end{eqnarray*}
\end{lemma}
\proof (i) Let first $n = 2$. If $\eta^{\dd}\in H_1^{\dd}$ and
$\xi = \xi_1\otimes\xi_2$ then
\begin{eqnarray*}
\sigma(A)(\theta(\xi))(\eta^{\dd}) & = & \theta(A_1\xi_1\otimes
A_2\xi_2)(\eta^{\dd}) = (A_1\xi_1,\eta)A_2\xi_2\\ & = &
(\xi_1,A_1^*\eta)A_2\xi_2 =
A_2\theta(\xi_1\otimes\xi_2)((A_1^*\eta)^{\dd})\\ & = &
A_2\theta(\xi_1\otimes\xi_2)A_1^{\dd}(\eta^{\dd}) =
A_2\theta(\xi)A_1^{\dd}(\eta^{\dd}).
\end{eqnarray*}
It follows by linearity and continuity that
$\sigma(A)(\theta(\xi)) = A_2\theta(\xi)A_1^{\dd}$, for every
$\xi\in H_1\otimes H_2$. Using Lemma \ref{l_compos} (i) we now
obtain
\begin{eqnarray*}
& & \sigma(A)(\theta(\xi))(\theta(\eta_{2,3})^{\dd})\dots
(\theta(\eta_{n-2,n-1}^{\dd}))\\
& = & \theta((A_1\otimes\dots
A_n)(\xi))(\theta(\eta_{2,3})^{\dd})\dots
(\theta(\eta_{n-2,n-1}^{\dd})) =  \theta((A_{n-1}\otimes
A_n)(\xi_{n-1,n}))\\&&\times\theta(\eta_{n-2,n-1}^{\dd})\dots\theta((A_3\otimes
A_4)(\xi_{3,4}))\theta(\eta_{2,3}^{\dd})\theta((A_1\otimes
A_2)(\xi_{1,2}))\\
& = &
A_n\theta(\xi_{n-1,n})A_{n-1}^{\dd}\theta(\eta_{n-2,n-1})^{\dd}A_{n-2}\dots
A_4\theta(\xi_{3,4})A_3^{\dd}\theta(\eta_{2,3})^{\dd}
A_2\theta(\xi_{1,2})A_1^{\dd}\\
& = & A_n\theta(\xi)(\theta((A_2^*\otimes
A_3^*)(\eta_{2,3}))^{\dd}))\dots (\theta((A_{n-2}^*\otimes
A_{n-1}^*)(\eta_{n-2,n-1}))^{\dd}))A_1^{\dd}.
\end{eqnarray*}

(ii) By Lemma \ref{l_compos} (ii),
\begin{eqnarray*}
& & \sigma(A)(\theta(\xi))(\theta(\eta_{1,2})^{\dd})\dots
(\theta(\eta_{n-2,n-1})^{\dd})\\
& = & \theta((A_1\otimes\dots
A_n)(\xi))(\theta(\eta_{1,2})^{\dd})\dots
(\theta(\eta_{n-2,n-1})^{\dd})\\
 &=&  \theta((A_{n-1}\otimes
A_n)(\xi_{n-1,n}))\theta(\eta_{n-2,n-1}^{\dd})\dots\theta(\eta_{1,2}^{\dd})(A_1\xi_1)\\
& = &
A_n\theta(\xi_{n-1,n})A_{n-1}^{\dd}\theta(\eta_{n-2,n-1})^{\dd}A_{n-2}\dots
A_2^{\dd}\theta(\eta_{1,2}^{\dd})(A_1\xi_1)\\
& = & A_n\theta(\xi)(\theta((A_1^*\otimes
A_2^*)(\eta_{1,2}))^{\dd}))\dots (\theta((A_{n-2}^*\otimes
A_{n-1}^*)(\eta_{n-2,n-1}))^{\dd})).
\end{eqnarray*}
\prend

Let $H_1,\dots,H_n$ be Hilbert spaces. If $n$ is even, we let
$$\Gamma(H_1,\dots,H_n) = (H_1\otimes H_2)\odot (H_2^{\dd}\otimes
H_3^{\dd})\odot (H_3\otimes H_4)\odot\dots\odot (H_{n-1}\otimes
H_n).$$ If $n$ is odd, we let
$$\Gamma(H_1,\dots,H_n) = (H_1^{\dd}\otimes H_2^{\dd})\odot (H_2\otimes
H_3)\odot (H_3^{\dd}\otimes H_4^{\dd})\odot\dots\odot
(H_{n-1}\otimes H_n).$$ After identifying $\bb{C}\otimes H_1$ with
$H_1$, for $n$ odd we have the identification
$$\Gamma(\bb{C},H_1,\dots,H_n) \equiv H_1\odot \Gamma(H_1,\dots,H_n).$$

Fix $\nph\in B(H)$. We define a mapping $S_{\nph}$ on
$\Gamma(H_1,\dots,H_n)$ taking values in $\cl B(H_1^{\dd},H_n)$ in
the case $n$ is even, and in $\cl B(H_1,H_n)$, in the case $n$ is
odd. Let first $n$ be even. On elementary tensors
$$\zeta =
\xi_{1,2}\otimes\eta_{2,3}^{\dd}\otimes\xi_{3,4}\otimes\dots\otimes\xi_{n-1,n}\in
\Gamma(H_1,\dots,H_n),$$ we let
$$S_{\nph}(\zeta) =
\sigma(\nph)\theta(\xi_{1,2}\otimes\xi_{3,4}\otimes\dots\otimes\xi_{n-1,n})
(\theta(\eta_{2,3}^{\dd}))\dots(\theta(\eta_{n-2,n-1}^{\dd}))$$
and extend $S_{\nph}$ on the whole of $\Gamma(H_1,\dots,H_n)$ by
linearity. Note that the values of $S_{\nph}$ are Hilbert-Schmidt
operators. Now assume $n$ is odd. Let
$\zeta\in\Gamma(H_1,\dots,H_n)$ and $\xi_1\in H_1$. Then
$$\xi_1\otimes\zeta\in H_1\odot\Gamma(H_1,\dots,H_n) =
\Gamma(\bb{C},H_1,\dots,H_n).$$ We let $S_{\nph}(\zeta)$ be the
operator defined on $H_1$ by
$$S_{\nph}(\zeta)(\xi_1) = S_{1\otimes\nph}(\xi_1\otimes\zeta).$$
Note that $S_{1\otimes\nph}(\xi_1\otimes\zeta)$ is an element of
$\cl C_2(\bb{C}^d,H_n)$, which can be identified with $H_n$ in a
natural way. In this way, $S_{\nph}(\zeta)(\xi_1)$ can be viewed
as an element of $H_n$. It is clear that the operator
$S_{\nph}(\zeta) : H_1\rightarrow H_n$ is linear. We moreover
claim that $S_{\nph}(\zeta)$ is bounded. Let $$\zeta =
\eta_{1,2}^{\dd}\otimes\dots\otimes\xi_{n-1,n}\in\Gamma(H_1,\dots,H_n)$$
and $\xi_1\in H_1$. Then $S_{\nph}(\zeta)$ is a bounded operator
and
\begin{equation}\label{pnes}
\|S_{\nph}(\zeta)\|_{\cl B(H_1,H_n)}\leq \|\nph\|_{\cl
B(H)}\|\eta_{1,2}\|\dots\|\eta_{n-2,n-1}\|\|\xi_{2,3}\|
\dots\|\xi_{n-1,n}\|.
\end{equation}
In fact, assuming for simplicity that $n=5$ we have
\begin{eqnarray*}
& & \|S_{\nph}(\zeta)(\xi_1)\| =
\|S_{1\otimes\nph}(\xi_1\otimes\zeta)\|\\ & = &\|\sigma(1\otimes
\nph)\theta((1\otimes\xi_1)\otimes\xi_{2,3}\otimes\xi_{4,5})
(\theta(\eta_{1,2}^{\dd}))(\theta(\eta_{3,4}^{\dd}))\|\\
& \leq &
\|\sigma(1\otimes\nph)\theta((1\otimes\xi_1)\otimes\xi_{2,3}\otimes\xi_{4,5})
(\theta(\eta_{1,2}^{\dd}))\|_{\op}\|(\theta(\eta_{3,4}^{\dd}))\|\\
& \leq &
\|\sigma(1\otimes\nph)\theta((1\otimes\xi_1)\otimes\xi_{2,3}\otimes\xi_{4,5})\|_{\op}
\|\eta_{1,2}\|\|\eta_{3,4}\|\\
& \leq & \|\nph\|_{\cl
B(H)}\|\xi_1\|\|\xi_{2,3}\|\|\xi_{4,5}\|\|\eta_{1,2}\|\|\eta_{3,4}\|\\
& = & \|\nph\|_{\cl B(H)}\|\zeta\|_{2,\wedge}\|\xi_1\|.
\end{eqnarray*}
Before proceeding, we identify two norms with which the space
$\Gamma(H_1,$ $\dots,$ $H_n)$ can be equipped. The first norm on
$\Gamma(H_1,\dots,H_n)$ is the projective tensor norm
$\|\cdot\|_{2,\wedge}$, where each of the terms $H_i\otimes
H_{i+1}$ (resp. $H_{i-1}^{\dd}\otimes H_i^{\dd}$) is given its
Hilbert space norm. In order to describe the second norm, note
that if $K_1$ and $K_2$ are Hilbert spaces then $K_1\otimes K_2$
can be endowed with an operator space structure by letting
$$\|(\xi_{ij})\| = \|\theta(\xi_{ji})\|_{M_m(\cl B(K_1^{\dd},K_2))}, \ \
(\xi_{ij})\in M_m(K_1\otimes K_2).$$ We write $(K_1\otimes
K_2)_{\op}^o$ for this operator space. Note that this is the
opposite operator space structure on $\cl C_2(K_1^{\dd},K_2)$,
after the identification of $K_1\otimes K_2$ and $\cl
C_2(K_1^{\dd},K_2)$. The norm $\|\cdot\|_{\hh}$ is the Haagerup
norm on $\Gamma(H_1,\dots,H_n)$ when $\Gamma(H_1,\dots,H_n)$ is
viewed as the algebraic tensor product of the operator spaces
$(H_i\otimes H_{i+1})_{\op}^o$ (resp. $(H_{i-1}^{\dd}\otimes
H_i^{\dd})_{\op}^o$). Thus, the norm $\|u\|_{\hh}$ of a finite sum
$u = \sum_i\xi_{1,2}^i\otimes\ldots\otimes\xi_{n-1,n}^i\in
\Gamma(H_1,\dots,H_n)$ of elementary tensors equals the Haagerup
norm of the element
$\sum_i\theta(\xi_{n-1,n}^i)\otimes\ldots\otimes\theta(\xi_{1,2}^i)$.

\begin{remark}\label{r_projnorm}
For each $\nph\in B(H)$ and each $\zeta\in\Gamma(H_1,\dots,H_n)$,
we have
$$\|S_{\nph}(\zeta)\|_{\op}\leq\|\nph\|_{\cl B(H)}\|\zeta\|_{2,\wedge}.$$
\end{remark}
\proof In the case where $n$ is odd and $\zeta$ is an elementary
tensor, the inequality coincides with (\ref{pnes}). In the case
$n$ is even and $\zeta$ is an elementary tensor, this is verified
similarly. The general case now follows by linearity. \prend

\begin{definition}\label{d_concmult}
An element $\nph\in B(H_1\otimes\dots\otimes H_n)$ is called a
concrete (operator) multiplier if there exists $C>0$ such that
$$\|S_{\nph}(\zeta)\|_{\op}\leq C\|\zeta\|_{\hh}, \ \mbox{ for each } \zeta\in
\Gamma(H_1,\dots,H_n).$$ The smallest such $C$ is denoted by
$\|\nph\|_{\mm}$.

Let $\cl A_1,\dots,\cl A_n$ be C*-algebras and $\pi_1,\dots,\pi_n$
be corresponding representations on Hilbert spaces
$H_1,\dots,H_n$. An element $\nph\in\cl A_1\otimes\dots\otimes\cl
A_n$ is called a $\pi_1,\dots,\pi_n$-multiplier if
$(\pi_1\otimes\dots\otimes\pi_n)(\nph)$ is a concrete multiplier.
We denote the set of all $\pi_1,\dots,\pi_n$-multipliers in $\cl
A_1\otimes\dots\otimes\cl A_n$ by ${\bf M}_{\pi_1,\dots,\pi_n}(\cl
A_1,\dots\cl A_n)$. If $\nph\in {\bf M}_{\pi_1,\dots,\pi_n}(\cl
A_1,\dots\cl A_n)$, we let $\|\nph\|_{\pi_1,\dots,\pi_n} =
\|(\pi_1\otimes\dots\otimes\pi_n)(\nph)\|_{\mm}$.

The element $\nph\in\cl A_1\otimes\dots\otimes\cl A_n$ is called a
universal multiplier if $\nph$ is a $\pi_1,\dots,\pi_n$-multiplier
for all representations $\pi_i$ of $\cl A_i$, $i = 1,\dots,n$. We
denote by ${\bf M}(\cl A_1,\dots\cl A_n)$ the set of all universal
multipliers in $\cl A_1\otimes\dots\otimes\cl A_n$.
\end{definition}

\begin{remark} In the case $n = 2$, Definition
\ref{d_concmult} reduces to the definition of $\cl
C_{\infty}$-multipliers studied in \cite{ks}.
\end{remark}

\bigskip

Next we show that an element $\varphi\in
L^{\infty}(X_1)\otimes\ldots\otimes L^{\infty}(X_n)$ $\subset$
$L^{\infty}(X_1\times\ldots\times X_n)$ is a Schur multiplier as
defined in Section~\ref{s_mm2} if and only if $\nph$ is a
$\pi_1,\ldots,\pi_n$-multiplier, where $\pi_i$ is the canonical
representation of $L^{\infty}(X_i)$ on $L^2(X_i)$ acting by
multiplication.

Let ${\cl A}$ be a commutative $C^*$-algebra with maximal ideal
space $X$, acting on a Hilbert space $H$. It is well-known that,
up to unitary equivalence, $H = \oplus_{\gamma\in \Gamma}
H_{\gamma}$, where $H_{\gamma} = L_2(X,\mu_{\gamma})$ is invariant
under ${\cl A}$ for each $\gamma\in\Gamma$, and an element
$f\in{\cl A}$ acts as on $H_{\gamma}$ by multiplication. Let $j:
H\to H$ be given by $\{\xi_{\gamma}(\lambda)\}\mapsto
\{\overline{\xi_{\gamma}(\lambda)}\}$. Then $V=\partial j$ is  a
unitary operator from $H$ to $H^{\dd}$ such that
$A^{\dd}=VAV^{-1}$ for all $A\in{\cl A}$. If $K$ is another
Hilbert space then $U(T)=TV$ (resp. $W(S)=V^{-1}S$) is an isometry
from ${\cl C}_2(H^{\dd}, K)$ to ${\cl C}_2(H, K)$ (resp. from
${\cl C}_2(K, H^{\dd})$ to ${\cl C}_2(K,H)$).

Let $\cl A_1,\ldots$, $\cl A_n$ be commutative $C^*$-algebras and
let $\pi_1,\ldots$, $\pi_n$ be corresponding representations on
$H_1,\ldots$, $H_n$ and $\pi=\pi_1\otimes\ldots\otimes\pi_n$. Let
$V_i: H_i\to H_i^{\dd}$ be unitary operator defined above with the
property $\pi_i(a_i)^{\dd} = V_i\pi_i(a_i)V_i^{-1}$ for each $a_i\in
\cl A_i$, $i = 1,\dots,n$. Define $U_{i,k}:{\cl C}_2(H_i^{\dd},
H_{k})\to {\cl C}_2(H_i, H_{k})$ and $W_{i,k}:{\cl C}_2(H_i,
H_{k}^{\dd})\to {\cl C}_2(H_i, H_k)$ to be $U_{i,k}(T)=TV_i$ and
$W_{i,k}(S)=V_{k}^{-1}S$. Then for $\varphi\in \cl
A_1\otimes\dots\otimes\cl A_n$, the mapping $S_{\pi(\varphi)}$ can
be identified with a mapping $\check S_{\pi(\varphi)}$ from ${\cl
C}_2(H_1,H_2)\odot{\cl C}_2(H_2,H_3)\odot\ldots\odot{\cl
C}_2(H_{n-1},H_n)$ into $\cl B(H_1,H_n)$ such that whenever
$\varphi=a_1\otimes\ldots\otimes a_n$ is an elementary tensor then
\begin{equation}\label{inten}
\check S_{\pi(\varphi)}(R_1\otimes\ldots\otimes
R_{n-1})=\pi_n(a_n)R_{n-1}\pi_{n-1}(a_{n-1})R_{n-2}\ldots
R_{1}\pi_1(a_1).
\end{equation}
In fact, let $\cl U=U_{1,2}\theta_{H_1,H_2}\otimes
W_{2,3}\theta_{H_2,H_3}\otimes\ldots\otimes
U_{n-1,n}\theta_{H_{n-1},H_n}$ if $n$ is even and $\cl
U=W_{1,2}\theta_{H_1,H_2}\otimes
U_{2,3}\theta_{H_2,H_3}\otimes\ldots\otimes
U_{n-1,n}\theta_{H_{n-1},H_n}$ if $n$ is odd. Then $\cl U$ maps
the space $\Gamma(H_1,H_2\ldots,H_n)$ onto ${\cl
C}_2(H_1,H_2)\odot{\cl C}_2(H_2,H_3)\odot\ldots\odot{\cl
C}_2(H_{n-1},H_n)$ and is an isometry with respect to the norm
$\|\cdot\|_{\hh}$ (this norm being defined on the algebraic tensor
product of the $\cl C_2$-spaces again as the Haagerup norm where
each of the $\cl C_2$-spaces is equipped with its opposite
operator space structure). Let
$$\check{S}_{\pi(\varphi)} = U_{1,n}S_{\pi(\varphi)}{\cl
U}^{-1}$$ in the case $n$ is even and
$$\check{S}_{\pi(\varphi)} = S_{\pi(\varphi)}{\cl
U}^{-1}$$ in the case $n$ is odd. Assume that
$\varphi=a_1\otimes\ldots\otimes a_n$. Then, in the case where $n$
is even, we have
\begin{eqnarray}\label{right}
&&\nonumber\check S_{\pi(\varphi)}(R_1\otimes\ldots\otimes
R_{n-1})\\
&&=\nonumber U_{1,n}S_{\pi(\varphi)}{\cl
U}^{-1}(R_1\otimes\ldots\otimes R_{n-1})\nonumber\\
&&=U_{1,n}(\pi_n(a_n)U_{n-1,n}^{-1}(R_{n-1})\pi_{n-1}(a_{n-1})^{\dd}W_{n-2,n-1}(R_{n-2})\ldots
\pi_1(a_1)^{\dd})\nonumber\\&&=
\pi_n(a_n)R_{n-1}V_{n-1}^{-1}\pi_{n-1}(a_{n-1})^{\dd}V_{n-1}R_{n-2}
\ldots R_{1}V_{1}^{-1}\pi_1(a_1)^{\dd}V_{1}\nonumber\\&&=\nonumber
\pi_n(a_n)R_{n-1}\pi_{n-1}(a_{n-1})R_{n-2}\ldots
R_{1}\pi_1(a_1).\nonumber
\end{eqnarray}
In the case where $n$ is odd one shows in a similar way that
(\ref{inten}) holds.

Now let ${\cl A}_i=L^{\infty}(X_i)$ and let $\pi_i$ be the
representation of ${\cl A}_i$ on $L^2(X_i)$ given by
$(\pi_i(f)\xi)(x)=f(x)\xi(x)$, $\xi \in L^2(X_i)$, $i=1,\ldots,n$.

Suppose $n$ is even. In this case
$\check{S}_{\pi(\nph)}(R_1\otimes\dots\otimes R_{n-1})$ is an
element of $\cl C_2(H_1,H_n)$. Using (\ref{right}) and the
identification $\psi_{k,l}:f\mapsto T_f$ of $L_2(X_k,X_l)$ with
the class of Hilbert-Schmidt operators from $L_2(X_k)$ to
$L_2(X_l)$, where
$$(T_f\xi)(y)=\int_{X_k}f(x,y)\xi(x) dx, \ \ f\in
L_2(X_k\times X_l), \xi\in L^2(X_k), y\in X_l,$$ we obtain that if
$f_1\otimes\ldots\otimes f_{n-1}\in \Gamma(X_1,\ldots,X_n)$ and
$\nph$ is an elementary tensor then
\begin{eqnarray}\label{number1}
&&\psi_{1,n}^{-1}(\check{S}
_{\pi(\varphi)}(\psi_{1,2}\otimes\ldots\otimes\psi_{n-1,n})(f_1\otimes\ldots\otimes
f_{n-1}))(x_1,x_n)\\
&&=\int_{X_2\times\ldots\times
X_{n-1}}\varphi(x_1,\ldots,x_n)f_1(x_1,x_2)\ldots
f_{n-1}(x_{n-1},x_n)dx_2\ldots
dx_{n-1}\nonumber\\&&=S_{\varphi}(f_1\otimes\ldots\otimes
f_{n-1})(x_1,x_n).\nonumber
\end{eqnarray}
By linearity and continuity, (\ref{number1}) holds for any
$\varphi\in L^{\infty}(X_1)\otimes\ldots\otimes L^{\infty}(X_n)$.

Now assume that $n$ is odd. Let $\xi\in H_1$, $\eta\in H_n$ and
$\psi_{0,1} : L^2(X_1)\rightarrow\cl C_2(\bb{C}, L^2(X_1))$ be the
natural identification. We have that
$(S_{\nph}(f_1\otimes\dots\otimes f_{n-1})\xi,\eta)$ coincides
with
$$(\check{S}_{(\id\otimes \pi)(1\otimes \nph)}(\psi_{0,1}\otimes\dots\otimes\psi_{n-1,n})((1\otimes\xi)\otimes
f_1\otimes\dots\otimes f_{n-1}),\eta)$$ whenever $\nph\in
L^{\infty}(X_1)\otimes\ldots\otimes L^{\infty}(X_n)$ is an
elementary tensor. By linearity and continuity, we have that
$\psi_{1,n}(S_{\nph}(f_1\otimes\dots\otimes f_{n-1}))$ is equal to
$$\check{S}_{\pi(\nph)}(\psi_{1,2}\otimes\dots\otimes\psi_{n-1,n})(f_1\otimes\dots\otimes
f_{n-1})$$ for all $\nph\in L^{\infty}(X_1)\otimes\ldots\otimes
L^{\infty}(X_n)$. In particular, $S_{\pi(\nph)}$ takes values in
$\cl C_2(H_1,H_n)$. As before, it follows that
\begin{eqnarray}\label{number2}
&&\psi_{1,n}^{-1}\check{S}_{\pi(\varphi)}
(\psi_{1,2}\otimes\ldots\otimes\psi_{n-1,n})(f_1\otimes\ldots\otimes
f_{n-1}))(x_1,x_n)\\&&= S_{\varphi}(f_1\otimes\ldots\otimes
f_{n-1})(x_1,x_n)\nonumber
\end{eqnarray}
for every $\nph\in L^{\infty}(X_1)\otimes\ldots\otimes
L^{\infty}(X_n)$. We have thus shown the following

\begin{proposition} \label{com+noncom}
An element $\varphi\in L^{\infty}(X_1)\otimes\ldots\otimes
L^{\infty}(X_n)$ is a Schur multiplier if and only if $\varphi\in
{\bf M}_{\pi_1,\ldots,\pi_n}(L^{\infty}(X_1),\ldots,
L^{\infty}(X_n))$.
\end{proposition}

Next we want to give a generalisation of Lemma \ref{l_compos2} for
the case where $\nph$ is a sum of elementary tensors. Let
$V$,$V_1,\dots,V_n$ be vector spaces, $L(V_1,V_2)$ be the space of
all linear mappings from $V_1$ into $V_2$ and $L(V) = L(V,V)$.
Recall that if $f : V_1\rightarrow V_2$ is a linear map, we let
$f_{k,l} : M_{k,l}(V_1)\rightarrow M_{k,l}(V_2)$ be the mapping
given by $f_{k,l}((v_{ij})) = (f(v_{ij}))$, for each $(v_{ij})\in
M_{k,l}(V_1)$. For an element $v = (v_{ij})\in M_{k,l}(V)$ we
denote by $v^{\ttt} = (v_{ji})\in M_{l,k}(V)$ the transpose of
$v$. Denote by $d : B(K)\rightarrow B(K^{\dd})$ the mapping
sending $A$ to its dual $A^{\dd}$. If $A=(A_{ij})\in
M_{k,l}(B(K))$ let $A^{\dd} = (A_{ij}^{\dd})$.

We will identify $M_{p,q}(\cl C_2(K_1,K_2))$ with $\cl
C_2(K_1^q,K_2^p)$. If $\xi\in M_{p,q}(K_1\otimes K_2)$ then
$\theta_{p,q}(\xi)\in M_{p,q}(\cl C_2(K_1^{\dd},K_2))$; using this
identification, we will be considering $\theta_{p,q}(\xi)$ as a
Hilbert-Schmidt operator from $K_1^q$ to $K_2^p$. If $A\in
B(K_1,K_2)$ then $A\otimes I_k \in B(K_1^k,K_2^k)$ is the $k$-fold
ampliation of $A$; under the identification $B(K_1^k,K_2^k)$ $=$
$M_k$ $(B(K_1,$ $K_2))$, the operator $A\otimes I_k$ has a $k$ by
$k$ diagonal matrix, whose every diagonal entry is $A$.

\begin{lemma}\label{l_opmul}
Let $V_1,\dots,V_n$ be vector spaces, $\cl L_i\subseteq
L(V_{i},V_{i+1})$ a subspace, $i = 1,\dots,n-1$, and
$$S : (L(V_n)\odot L(V_{n-1})\odot\dots\odot
L(V_1))\times (\cl L_{n-1}\odot\dots\odot \cl L_1)\rightarrow
L(V_1,V_n)$$ be a mapping satisfying
$$S(a_n\otimes\dots\otimes
a_1,\lambda_{n-1}\otimes\dots\otimes\lambda_1) =
a_n\lambda_{n-1}a_{n-1}\dots\lambda_1a_1.$$ If $A_1\in
M_{k_1,1}(L(V_1))$, $A_2\in M_{k_2,k_1}(L(V_2)),\dots$,$A_n\in
M_{1,k_{n-1}}(L(V_n))$ and $\Lambda_1\in M_{l_1,1}(\cl L_1)$,
$\Lambda_2\in M_{l_2,l_1}(\cl L_2),\dots$, $\Lambda_{n-1}$ $\in$
$M_{1,l_{n-2}}(\cl L_{n-1})$ then
$$S(A_n\odot\dots\odot
A_1,\Lambda_{n-1}\odot\dots\odot\Lambda_1) =
A_n\dots(\Lambda_2\otimes I_{k_2})(A_2\otimes
I_{l_1})(\Lambda_1\otimes I_{k_1})A_1 .$$
\end{lemma}
\proof \lq\lq A few moments' thought.'' \prend

\begin{lemma}\label{opmul2}
Let $A_1\in M_{1,k_1}(\cl B(H_1))$, $A_2\in M_{k_1,k_2}(\cl
B(H_2)),\dots,$ $A_n$ $\in$ $M_{k_{n-1},1}$ $(\cl B(H_n))$ and
$\nph = A_1\odot A_2\odot\dots\odot A_n$.

(i) Assume $n$ is even. Let $\xi_{1,2}\in M_{1,l_1}(H_1\otimes
H_2)$, $\eta_{2,3}\in M_{l_1,l_2}(H_2^{\dd}\otimes
H_3^{\dd}),\dots$, $\xi_{n-1,n}\in M_{l_{n-2},1}(H_{n-1}\otimes
H_n)$ and
$$\zeta =
\xi_{1,2}\odot\eta_{2,3}\odot\dots\odot\xi_{n-1,n}\in
\Gamma(H_1,\dots,H_n).$$ Then
$$S_{\nph}(\zeta) =
A_n^{\ttt}\dots(A_3^{\ttt,\dd}\otimes
I_{l_2})(\theta_{l_1,l_2}(\eta_{2,3})^{\ttt}\otimes
I_{k_2})(A_2^{\ttt,\dd}\otimes
I_{l_1})(\theta_{1,l_1}(\xi_{1,2})^{\ttt}\otimes
I_{k_1})A_1^{\ttt,\dd}.$$

(ii) Assume $n$ is odd. Let $\eta_{1,2}\in
M_{1,l_1}(H_1^{\dd}\otimes H_2^{\dd})$, $\xi_{2,3}\in
M_{l_1,l_2}(H_2\otimes H_3),\dots$, $\xi_{n-1,n}\in
M_{l_{n-2},1}(H_{n-1}\otimes H_n)$ and
$$\zeta =
\eta_{1,2}\odot\xi_{2,3}\odot\dots\odot\xi_{n-1,n}\in
\Gamma(H_1,\dots,H_n).$$ Then
$$S_{\nph}(\zeta) =
A_n^{\ttt}\dots(A_3^{\ttt}\otimes
I_{l_2})(\theta_{l_1,l_2}(\xi_{2,3})^{\ttt}\otimes
I_{k_2})(A_2^{\ttt,\dd}\otimes
I_{l_1})(\theta_{1,l_1}(\eta_{1,2})^{\ttt}\otimes
I_{k_1})A_1^{\ttt}.$$
\end{lemma}
\proof Let $f : V_1\odot\dots\odot V_n \rightarrow
V_n\odot\dots\odot V_1$ be the flip, namely the map given on
elementary tensors by $f(v_1\otimes\dots\otimes v_n) =
v_n\otimes\dots\otimes v_1$. Note that if $A_1\in M_{1,k_1}(V_1)$,
$A_2\in M_{k_1,k_2}(V_2),\dots$, $A_n\in M_{k_{n-1},1}(V_n)$ then
$$f(A_1\odot\dots\odot A_n) = A_n^{\ttt}\odot\dots\odot
A_1^{\ttt}.$$ Let
$$D : B(H_1)\odot B(H_2)\odot\dots\odot B(H_n) \longrightarrow B(H_n)\odot
B(H_{n-1}^{\dd})\odot\dots\odot B(H_1^{\dd})$$ be the map
$$D = f\circ (d\otimes\id\otimes d\otimes\dots\otimes \id).$$
We have that
$$D(A) = A_n^{\ttt}\odot A_{n-1}^{\ttt,\dd}\odot\dots\odot A_1^{\ttt,\dd}.$$

Define a mapping $S$ from
$$(B(H_n)\odot B(H_{n-1}^{\dd})\odot\dots\odot
B(H_1^{\dd}))\times (\cl C_2(H_{n-1}^{\dd},H_n)\odot\dots\odot\cl
C_2(H_1^{\dd},H_2))$$ into $\cl C_2(H_1^{\dd},H_n)$ by
$$S(\psi,\zeta') = S_{D^{-1}(\psi)}(\tilde{\theta}^{-1}(\zeta')),$$
where
$$\tilde{\theta} : \Gamma(H_1,\dots,H_n) \rightarrow \cl
C_2(H_{n-1}^{\dd},H_n)\odot\dots\odot\cl C_2(H_1^{\dd},H_2)$$ is
given on elementary tensors by
$$\tilde{\theta}(\xi_{1,2}\otimes\eta_{2,3}\otimes\dots\otimes\xi_{n-1,n})
=
\theta(\xi_{n-1,n})\otimes\dots\otimes\theta(\eta_{2,3})\otimes\theta(\xi_{1,2}).$$
By Lemma \ref{l_compos2} (i), the mapping $S$ satisfies the
requirements of Lemma \ref{l_opmul} and
$$S_{\nph}(\zeta) = S(A_n^{\ttt}\odot
A_{n-1}^{\ttt,\dd}\odot\dots\odot A_1^{\ttt,\dd},
\theta_{l_{n-2},1}(\xi_{n-1,n})^{\ttt}\odot\dots\odot\theta_{1,l_1}(\xi_{1,2})^{\ttt}).$$
The claim now follows from Lemma \ref{l_opmul}.

The proof of (ii) is similar. \prend

\section{Multipliers for tensor products of representations}

It was proved in \cite{ks} that the space of all
$(\pi,\rho)$-multipliers does not change if the representations
$\pi$ and $\rho$ are replaced by approximately equivalent
representations.  In this section we will prove a corresponding
result for multidimensional multipliers. We first recall the
notion of approximate equivalence and approximate suborditation
introduced by Voiculescu in \cite{voiculescu}.

Let $\pi$ and $\pi'$ be $*$-representations of a $C^*$-algebra
${\cl A}$ on Hilbert spaces $H$ and $H'$, respectively. We say
that $\pi'$ is {\it approximately subordinate} to $\pi$ and write
$\pi'\stackrel{a}\ll\pi$ if there is a net $\{U_{\lambda}\}$ of
isometries from $H'$ to $H$ such that
\begin{equation}\label{as}
\|\pi(a)U_{\lambda}-U_{\lambda}\pi'(a)\|\to 0\text{ for all
}a\in{\cl A}.
\end{equation}
The representations $\pi'$ and $\pi$ are said to be {\it
approximately equivalent} if the operators $U_{\lambda}$ can be
chosen to be unitary; in this case we write
$\pi'\stackrel{a}\sim\pi$.

For C*-algebras $\cl A_1,\ldots, \cl A_n$ and corresponding
representations $\pi_1,\ldots,\pi_n$, we will denote the
collection of all $\pi_1,\ldots,\pi_n$-multipliers in $\cl
A_1\otimes\dots\otimes\cl A_n$ simply by ${\bf
M}_{\pi_1,\ldots,\pi_n}$, in case there is no danger of confusion.

\begin{theorem}\label{ap}
Let $\cl A_1,\dots,\cl A_n$ be C*-algebras and $\pi_i$ and
$\pi_i'$ be representations of $\cl A_i$ on Hilbert spaces $H_i$
and $H_i'$, respectively, $i = 1,\dots,n$.

(i) If $\pi_i'\stackrel{a}\ll\pi_i$, $i=1,\ldots,n$, then
$${\bf M}_{\pi_1,\ldots,\pi_n}\subseteq {\bf M}_{\pi_1',\ldots,\pi_n'} \text{
and } \|\varphi\|_{\pi_1',\ldots,\pi_n'}\leq
\|\varphi\|_{\pi_1,\ldots,\pi_n},\text{ for } \varphi\in {\bf
M}_{\pi_1,\ldots,\pi_n}.$$

(ii) If $\pi_i'\stackrel{a}\sim\pi_i$, $i=1,\ldots,n$, then
$${\bf M}_{\pi_1,\ldots,\pi_n}= {\bf M}_{\pi_1',\ldots,\pi_n'} \text{ and }
\|\varphi\|_{\pi_1,\ldots,\pi_n}=
\|\varphi\|_{\pi_1',\ldots,\pi_n'},\text{ for } \varphi\in {\bf
M}_{\pi_1,\ldots,\pi_n}.$$
\end{theorem}
\proof (i) Let first $n$ be even and $\{U_{\lambda_i}\}$ be nets
of isometries from $H_i'$ into $H_i$ satisfying
$$\|\pi_i(a_i)U_{\lambda_i}-U_{\lambda_i}\pi_i'(a_i)\|\to 0, \text{ for all }
a_i\in{\cl A_i}.$$ Set $\pi=\otimes_{i=1}^n\pi_i$,
$\pi'=\otimes_{i=1}^n\pi_i'$,
$\lambda=(\lambda_1,\ldots,\lambda_n)$ and
$W_{\lambda}=U_{\lambda_1}\otimes\ldots\otimes U_{\lambda_n}$.
Then $W_{\lambda}$ are isometries from $\otimes_{i=1}^n H_i'$ to
$\otimes_{i=1}^n H_n$ and, for $x\in {\cl A_1}\odot\ldots\odot{\cl
A_n}$, we have
$$\|\pi(x)W_{\lambda}-W_{\lambda}\pi'(x)\|\longrightarrow
0.$$ As $\|W_{\lambda}\|=1$ for all $\lambda$, this holds for all
$x\in {\cl A_1}\otimes\ldots\otimes{\cl A_n}$. By Lemma
\ref{l_compos2} (i) we have that, for any $\xi\in \otimes_{i=1}^n
H_i$,
\begin{eqnarray*}
&&\theta(W^*_{\lambda}\xi)(\theta(\eta_{2,3}^{\dd}))\dots
(\theta(\eta_{n-2,n-1}^{\dd}))\\
&&=U_{\lambda_n}^*\theta(\xi)(\theta((W_{\lambda_2,\lambda_3}\eta_{2,3})^{\dd}))\dots
(\theta((W_{\lambda_{n-2},\lambda_{n-1}}\eta_{n-2,n-1})^{\dd}))(U_{\lambda_1}^*)^{\dd}
\end{eqnarray*}
where $W_{\lambda_k,\lambda_{k+1}}=U_{\lambda_k}\otimes
U_{\lambda_{k+1}}$. Therefore, if
$\displaystyle\zeta=\xi_{1,2}\otimes(\eta_{2,3})^{\dd}\otimes\ldots\otimes\xi_{n-1,n}$,
then
\begin{eqnarray}\label{star}
&&S_{W^*_{\lambda}\pi(\varphi)W_{\lambda}}(\zeta)=\\
&&=U_{\lambda_n}^*S_{\pi(\varphi)}(W_{\lambda_1,\lambda_2}\xi_{1,2}\otimes
(W_{\lambda_2,\lambda_3}\eta_{2,3})^{\dd}\otimes\ldots\otimes
W_{\lambda_{n-1},\lambda_n}\xi_{n-1,n})(U_{\lambda_1}^*)^{\dd}.\nonumber
\end{eqnarray}
Let $\Gamma_{\lambda}:\Gamma(H_1',\ldots, H_n')\to
\Gamma(H_1,\ldots, H_n)$ be the linear operator defined on
elementary tensors by
$$\Gamma_{\lambda}(\xi_{1,2}\otimes\eta_{2,3}^{\dd}\otimes\ldots\otimes\xi_{n-1,n})=
W_{\lambda_1,\lambda_2}\xi_{1,2}\otimes(W_{\lambda_2,\lambda_3}\eta_{2,3})^{\dd}\otimes\ldots\otimes
W_{\lambda_{n-1},\lambda_n}\xi_{n-1,n}.$$ It follows from
(\ref{star}) and Remark \ref{r_projnorm} that if $\varphi\in {\bf
M}_{\pi_1,\ldots,\pi_n}$ and $\zeta\in \Gamma(H_1',\ldots, H_n')$
then
\begin{eqnarray*}
\|S_{\pi'(\varphi)}(\zeta)\|_{\op}&\leq&\|S_{W_{\lambda}^*\pi(\varphi)
W_{\lambda}}(\zeta)\|_{\op} + \|S_{W_{\lambda}^*\pi(\varphi)
W_{\lambda}-\pi'(\varphi)}(\zeta)\|_{\op}\\
&\leq&\|S_{\pi(\varphi)}(\Gamma_{\lambda}\zeta)\|_{\op}+
\|S_{W_{\lambda}^*\pi(\varphi)
W_{\lambda}-\pi'(\varphi)}(\zeta)\|_{\op}\\ &\leq&
\|\varphi\|_{\pi_1,\ldots,\pi_n}\|\Gamma_{\lambda}\zeta\|_{\hh} +
\|W_{\lambda}^*\pi(\varphi)
W_{\lambda}-\pi'(\varphi)\|_{\op}\|\zeta\|_{2,\wedge}.
\end{eqnarray*}
Since $\|W_{\lambda}^*\pi(\varphi)
W_{\lambda}-\pi'(\varphi)\|_{\op}\rightarrow 0$, in order to prove
that $\nph\in {\bf M}_{\pi_1',\ldots,\pi_n'}$, it suffices to show
that $\|\Gamma_{\lambda}\zeta\|_{\hh}\leq \|\zeta\|_{\hh}$. If
$\xi_{i,i+1}\in H'_i\otimes H'_{i+1}$ then
$\theta(W_{\lambda_i,\lambda_{i+1}}\xi_{i,i+1})=U_{\lambda_{i+1}}\theta(\xi_{i,i+1})U_{\lambda_i}^{\dd}.$
Let $\zeta\in \Gamma(H_1',\ldots, H_n')$ be of the form
$$\zeta=\xi_{1,2}\otimes \eta_{2,3}^{\dd}\otimes\ldots\otimes \xi_{n-1,n}$$
where $\xi_{1,2}\in M_{1,k_2}(H_1'\otimes H_2')$,
$\eta_{2,3}^{\dd}\in M_{k_2,k_3}((H_2')^{\dd}\otimes
(H_3')^{\dd}),\ldots$, and $\xi_{n-1,n}\in
M_{k_{n-1},1}(H_{n-1}'\otimes H_n')$ are such that
$$\|\zeta\|_{\hh}=
\|\theta_{1,k_2}(\xi_{1,2})^{\ttt}\|_{\op}\|\theta_{k_2,k_3}(\eta_{2,3}^{\dd})^{\ttt}\|_{\op}\ldots
\|\theta_{k_{n-1},1}(\xi_{n-1,n})^{\ttt}\|_{\op}.$$ Then
$$\Gamma_{\lambda}\zeta=W_{\lambda_1,\lambda_2}\xi_{1,2}\odot
(W_{\lambda_2,\lambda_3}^{*,\dd}\otimes
I_{k_2})\eta_{2,3}^{\dd}\odot \ldots\odot
(W_{\lambda_{n-1},\lambda_n}\otimes I_{k_{n-1}})\xi_{n-1,n}$$ and
as
\begin{eqnarray*}
\displaystyle\theta_{1,k_2}(W_{\lambda_1,\lambda_2}\xi_{1,2})&=&U_{
\lambda_2}\theta_{1,k_2}(\xi_{1,2})(U_{\lambda_1}^{\dd}\otimes I_{k_2}),\\
\displaystyle\theta_{k_2,k_3}(((W_{\lambda_2,\lambda_3}^*)^{\dd}\otimes
I_{k_2})\eta_{2,3}^{\dd})&=&(U_{\lambda_3}^{\dd}\otimes
I_{k_2})\theta_{2,3}(\eta_{2,3}^{\dd})(U_{\lambda_2}\otimes I_{k_3}),\\
\dots\dots\dots\dots\dots\dots\dots\dots\dots\dots&\dots&\dots\dots\dots\dots\dots\dots\dots\dots\dots\dots\\
\theta_{k_{n-1}, 1}((W_{\lambda_{n-1},\lambda_n}\otimes
I_{k_{n-1}})\xi_{n-1,n})&=&(U_{\lambda_n}\otimes
I_{k_{n-1}})\theta_{k_{n-1},1}(\xi_{n-1,n})U_{\lambda_{n-1}}^{\dd},
\end{eqnarray*}
we get
\begin{eqnarray*}
\|\Gamma_{\lambda}\zeta\|_{\hh}&\leq &\|U_{\lambda_2}\otimes
I_{k_2}\|_{\op}
\|\theta_{1,k_2}(\xi_{1,2})^{\ttt}\|_{\op}\|U_{\lambda_1}^{\dd}\|_{\op}\ldots\\
& \dots &
\|\theta_{k_{n-1},1}(\xi_{n-1,n})^{\ttt}\|_{\op}\|U_{\lambda_{n-1}}^{\dd}\otimes
I_{k_{n-1}}\|_{\op}\\
&=&\|\theta_{1,k_2}(\xi_{1,2})^{\ttt}\|_{\op}\ldots\|\theta_{k_{n-1},1}(\xi_{n-1,n})^{\ttt}\|_{\op}=\|\zeta\|_{\hh}
\end{eqnarray*}

This completes the proof for the case where $n$ is even. Now
assume that $n$ is odd and let
$\Gamma_{\lambda}:\Gamma(H_1',\ldots, H_n')\to \Gamma(H_1,\ldots,
H_n)$ be the linear operator defined on elementary tensors by
$$\Gamma_{\lambda}(\xi_{1,2}^{\dd}\otimes\ldots
\otimes\eta_{n-1,n})=(W_{\lambda_1,\lambda_2}\xi_{1,2})^{\dd}\otimes
\otimes\ldots\otimes W_{\lambda_{n-1},\lambda_n}\eta_{n-1,n}.$$ An
estimate similar to the above shows again that
$\|\Gamma_{\lambda}\zeta\|_{\hh}\leq \|\zeta\|_{\hh}.$

By the definition of the map $S_{\pi'(\varphi)}$ and the arguments
above, we obtain
\begin{eqnarray*}
&&\|S_{\pi'(\varphi)}(\zeta)\|_{\op}  \leq
\|S_{W_{\lambda}^*\pi(\varphi) W_{\lambda}}(\zeta)\|_{\op} +\|S_{
(W_{\lambda}^*\pi(\varphi)
W_{\lambda}-\pi'(\varphi))}(\zeta)\|_{\op}\\
&& = \sup_{\xi_1\in H_1',\|\xi_1\|=1}\|S_{1\otimes
W_{\lambda}^*\pi(\varphi) W_{\lambda}}(\xi_1\otimes\zeta)\|_{H_n'}
+ \|S_{ (W_{\lambda}^*\pi(\varphi)
W_{\lambda}-\pi'(\varphi))}(\zeta)\|_{\op}\\
&& \leq \sup_{\xi_1\in
H_1',\|\xi_1\|=1}\|S_{1\otimes\pi(\varphi)}(U_{\lambda_1}\xi_1\otimes
\Gamma_{\lambda}\zeta)\|_{H_n} +\|S_{ (W_{\lambda}^*\pi(\varphi)
W_{\lambda}-\pi'(\varphi))}(\zeta)\|_{\op}\\
\end{eqnarray*}
\begin{eqnarray*}
&&\leq \sup_{\eta_1\in
H_1,\|\eta_1\|=1}\|S_{1\otimes\pi(\varphi)}(\eta_1\otimes\Gamma_{\lambda}\zeta)\|_{H_n}+\|W_{\lambda}^*\pi(\varphi)
W_{\lambda}-\pi'(\varphi)\|_{\op}\|\zeta\|_{2,\wedge}\\
&&=\|S_{\pi(\varphi)}(\Gamma_{\lambda}\zeta)\|_{\op}+\|W_{\lambda}^*\pi(\varphi)
W_{\lambda}-\pi'(\varphi)\|_{\op}\|\zeta\|_{2,\wedge}\\
&&\leq
\|\varphi\|_{\pi_1,\ldots,\pi_n}\|\|\Gamma_{\lambda}\zeta\|_{\hh}
+\|W_{\lambda}^*\pi(\varphi)
W_{\lambda}-\pi'(\varphi)\|_{\op}\|\zeta\|_{2,\wedge}\\
&&\leq\|\varphi\|_{\pi_1,\ldots,\pi_n}\|\|\zeta\|_{\hh}+\|W_{\lambda}^*\pi(\varphi)
W_{\lambda}-\pi'(\varphi)\|_{\op}\|\zeta\|_{2,\wedge}.
\end{eqnarray*}
As $\|W_{\lambda}^*\pi(\varphi)
W_{\lambda}-\pi'(\varphi)\|_{\op}\to 0$ we obtain the desired
statement.

(ii) is a direct consequence of (i). \prend

For $T\in B(H)$,  set $\text{rank}(T)=\overline{\text{dim}(TH)}$.
It was proved in \cite[Theorem~5.1]{hadwin}  that for
$*$-representations $\pi$ and $\pi'$ of a $C^*$-algebra $\cl A$
\begin{equation}\label{hadwin}
\pi'\stackrel{a}\ll\pi \ \Longleftrightarrow \
\text{rank}(\pi'(a))\leq\text{rank}(\pi(a))\text{ for each }
a\in{\cl A}.
\end{equation}

The next statement is a multidimensional version of
\cite[Corollory~5.3]{ks}. Its proof follows  the lines of the
proof of the corresponding statement in the two dimensional case
and uses  Theorem~\ref{ap} instead of \cite[Theorem~5.2]{ks}.
\begin{corollary}
Let $\pi_i$, $\pi_i'$ be representations of  separable
$C^*$-algebra $\cl A_i$, $i=1,\ldots,n$. Assume that
$$\min\{\aleph_0,\rank(\pi_i'(a_i))\}\leq \min\{\aleph_0,\rank(\pi_i(a_i))\}$$
for each $a_i\in\cl A_i$ and $i=1,\ldots,n$.


Then ${\bf M}_{\pi_1,\ldots,\pi_n}\subseteq{\bf
M}_{\pi_1',\ldots,\pi_n'}$ and
$\|\varphi\|_{\pi_1',\ldots,\pi_n'}\leq
\|\varphi\|_{\pi_1,\ldots,\pi_n}$ for $\varphi\in {\bf
M}_{\pi_1,\ldots,\pi_n}$.
\end{corollary}

Recall that a $*$-representation $\pi$ of a $C^*$-algebra $\cl A$
has a separating vector if there is a cyclic vector for the
commutant $\pi(\cl A)'$.

\begin{lemma}\label{smith+}
Let $\cl H, H_1,\ldots, H_n$ be Hilbert spaces,
$\pi_1,\ldots,\pi_n$ be representations of the $C^*$-algebras $\cl
A_1,\ldots,\cl A_n$ on $H_1,\ldots, H_n$ and $\pi_i\otimes 1$ be
the amplification of $\pi_i$ on $H_i\otimes \cl H$, respectively.
Assume that $\pi_1$ and $\pi_n$ have separating vectors. Then
$${\bf M}_{\pi_1,\ldots,\pi_n}={\bf M}_{\pi_1\otimes 1,\ldots,\pi_n\otimes
1}$$ and the multiplier norms on these spaces coincide.
\end{lemma}
\proof We use ideas from the proofs of \cite[Theorem~2.1]{smith}
and Lemma~\ref{l_gsmith}. For simplicity we assume that $n=3$ and
that $\cl H$ is separable. Let $\varphi\in {\bf
M}_{\pi_1,\pi_2,\pi_3}$ with $\|\varphi\|_{\pi_1,\pi_2,\pi_3}=1$
and set $S = S_{(\pi_1\otimes 1)\otimes(\pi_2\otimes
1)\otimes(\pi_3\otimes 1)(\varphi)}$. The mapping $S$ can be
regarded as a mapping on
\begin{equation}\label{space}
{\cl C}_2((H_2\otimes {\cl H})^{\dd}, H_3\otimes {\cl H})\odot{\cl
C}_2((H_1\otimes {\cl H}), (H_2\otimes {\cl H})^{\dd})
\end{equation}
by setting
$S(\theta(\xi_{2,3})\otimes\theta(\eta_{1,2}^{\dd}))=S(\eta_{1,2}^{\dd}\otimes\xi_{2,3})$
for $\zeta=\eta_{1,2}^{\dd}\otimes\xi_{2,3}\in \Gamma(H_1\otimes
{\cl H}, H_2\otimes {\cl H}, H_3\otimes {\cl H})$.  Similarly, the
mapping $S_{\pi_1\otimes\pi_2\otimes\pi_3(\varphi)}$ can be
regarded as a mapping on ${\cl C}_2(H_2^{\dd}, H_3)\odot{\cl
C}_2((H_1, H_2^{\dd})$. It follows from Lemma~\ref{opmul2} that
$S_{\pi_1\otimes\pi_2\otimes\pi_3(\varphi)}$ is $\pi_3(\cl
A_3)',(\pi_2(\cl A_2)')^{\dd},\pi_1(\cl A_1)'$-modular.

Assume that $\|\varphi\|_{\pi_1\otimes 1,\pi_2\otimes
1,\pi_3\otimes 1}>1$. Then there exists an element
$T=(T_1^2,\ldots,T_s^2)\odot (T_1^1,\ldots,T_s^1)^{\ttt}$ in the
space defined in  (\ref{space})  with
$$\|\sum (T_i^1)^*T_i^1\|\|\sum
T_i^2(T_i^2)^*\|=1$$ and vectors $\xi_0\in H_1\otimes {\cl H}$,
$\eta_0\in H_3\otimes {\cl H}$ of norm less than one such that
$$|(S(T)\xi_0,\eta_0)|>1.$$

Fix a basis $\{f_l\}$ of ${\cl H}$ and denote by $P_n$ the
projection onto the space generated by the first $n$ vectors in
this basis. Then, as
$$(1_{H_3}\otimes P_n) S(T)(1_{H_1}\otimes P_n)\to S(T),$$ weakly, there exists
$n\geq 1$ such that
$$|((1_{H_3}\otimes P_n) S(T)(1_{H_1}\otimes P_n)\xi_0,\eta_0)|>1.$$ Thus we
may assume that $\xi_0\in H_1\otimes P_n{\cl H}$ and $\eta_0\in
H_3\otimes P_n{\cl H}$, say $$\xi_0=(\xi_1,\ldots,\xi_n,0,\ldots),
\eta_0=(\eta_1,\ldots,\eta_n,0.\ldots).$$ As $\pi_1({\cl A_1})'$
and $\pi_3({\cl A_3})'$ have cyclic vectors, say $\xi$ and $\eta$
respectively, we may assume that $\xi_i=a_i\xi$, $\eta_i=b_i\eta$
for some $a_i\in \pi_1(\cl A_1)'$ and $b_i\in \pi_3(\cl A_3)'$.
Let $a=\sum a_i^*a_i$, $b=\sum b_i^*b_i$. Assuming first that $a$,
$b$ are invertible we set $\tilde a_i=a_ia^{-1/2}$, $\tilde
b_i=b_ib^{-1/2}$. Then for $\tilde\xi=a^{1/2}\xi$,
$\tilde\eta=b^{1/2}\eta$ we have $\xi_i=\tilde a_i\tilde\xi$ and
$\eta_i=\tilde b_i\tilde\eta$. We write $T_i^k=((T_i^k)_{lm})$,
where $(T_i^1)_{lm}=(1_{H_2^{\dd}}\otimes
P(f_l^{\dd}))T_i^1(1_{H_1}\otimes P(f_m))$,
$(T_i^2)_{lm}=(1_{H_3}\otimes P(f_l))T_i^2(1_{H_2^{\dd}}\otimes
P(f_m^{\dd}))$, where $P(f)$ is the projection onto the one
dimensional space generated by $f$. Using the modularity of
$S_{\pi_1\otimes\pi_2\otimes\pi_3(\varphi)}$, we obtain
\begin{eqnarray*}
|(S(T)\xi_0,\eta_0)|&=&\left|\sum_{i=1}^s (S(T_i^2\otimes
T_i^1)\xi_0,\eta_0)\right|\nonumber\\
\end{eqnarray*}
\begin{eqnarray}\label{ineq}
&=&\left|\sum_{i=1}^s\sum_{l,m=1}^n\sum_{k=1}^{\infty}
(S_{\pi_1\otimes\pi_2\otimes\pi_3(\varphi)}((T_i^2)_{lk}\otimes
(T_i^1)_{km})\tilde a_m\tilde\xi,\tilde
b_l\tilde\eta)\right|\\
&=&\left|\sum_{i=1}^s\sum_{l,m=1}^n\sum_{k=1}^{\infty}
(S_{\pi_1\otimes\pi_2\otimes\pi_3(\varphi)}(\tilde
b_l^*(T_i^2)_{lk}\otimes (T_i^1)_{km}\tilde a_m)\tilde\xi,
\tilde\eta)\right|.\nonumber
\end{eqnarray}
The next step is to prove that
$\displaystyle\sum_{i=1}^s\sum_{k=1}^{\infty}
\left(\sum_{l=1}^n\tilde b_l^*(T_i^2)_{lk}\right)\otimes
\left(\sum_{m=1}^n(T_i^1)_{km}\tilde a_m\right)$ belongs to $\cl
K(H_2^{\dd},H_3)\otimes_{\hh}\cl K(H_1,H_2^{\dd})$. Observe first
that the row operator
$$\displaystyle R_i=(\sum_{l=1}^n\tilde
b_l^*(T_i^2)_{l1},\ldots, \sum_{l=1}^n\tilde
b_l^*(T_i^2)_{lk},\ldots )$$ is equal to the product of the row
operator $\tilde B=(\tilde b_1,\ldots,\tilde b_n,0,\ldots )$ and
the Hilbert-Schmidt operator $T_i^2$. Set $R = (R_1,\ldots, R_s) =
(\tilde B T_1^2,\ldots,\tilde B T_s^2)$.

As each $T_i^2$ is the operator norm-limit of operators
$T_i^2(1_{H_2^{\dd}}\otimes P_k)$ as $k\to\infty$, the operator
$R_i$ is the uniform limit of the sequence of truncated operators
$R_i^k=(\sum_{l=1}^n\tilde
b_l^*(T_i^2)_{l1},\ldots,\sum_{l=1}^n\tilde
b_l^*(T_i^2)_{lk},0\ldots )$. Thus
$$RR^*=\sum_{i=1}^s\sum_{k=1}^{\infty}\left(\sum_{l=1}^n\tilde
b_l^*(T_i^2)_{lk}\right)\left(\sum_{l=1}^n\tilde
b_l^*(T_i^2)_{lk}\right)^*,$$ where the series converges uniformly
and
\begin{eqnarray*}
\|\sum_{i=1}^s\sum_{k=1}^{\infty}(\sum_{l=1}^n\tilde
b_l^*(T_i^2)_{lk})(\sum_{l=1}^n\tilde
b_l^*(T_i^2)_{lk})^*\|=\|RR^*\|=\|\sum_{i=1}^sR_iR_i^*\|\\
=\|\tilde B(\sum_{i=1}^sT_i^2(T_i^2)^*)\tilde B^*\|\leq \|\tilde
B\|^2\|\|\sum_{i=1}^sT_i^2(T_i^2)^*\|\leq 1.
\end{eqnarray*}
In the same way one shows that the series
$\displaystyle\sum_{k=1}^{\infty}(\sum_{m=1}^n(T_i^1)_{km}\tilde
a_m)(\sum_{m=1}^n(T_i^1)_{km}\tilde a_m)^*$ converges uniformly
and
$$\|\sum_{i=1}^s\sum_{k=1}^{\infty}(\sum_{m=1}^n(T_i^1)_{km}\tilde
a_m)(\sum_{m=1}^n(T_i^1)_{km}\tilde a_m)^*\|\leq 1.$$ Thus
$\displaystyle\sum_{i=1}^s\sum_{k=1}^{\infty}(\sum_{l=1}^n\tilde
b_l^*(T_i^2)_{lk})\otimes (\sum_{m=1}^n(T_i^1)_{km}\tilde a_m)\in
\cl K(H_1,H_2^{\dd})\otimes_{\hh}\cl K(H_2^{\dd},H_3)$ and
$$\|\sum_{i=1}^s\sum_{k=1}^{\infty}(\sum_{l=1}^n\tilde
b_l^*(T_i^2)_{lk})\otimes (\sum_{m=1}^n(T_i^1)_{km}\tilde
a_m)\|_{\hh}\leq 1.$$

Next
$\|\tilde\xi\|^2=(b^{1/2}\xi,b^{1/2}\xi)=(b\xi,\xi)=\sum_i(b_i\xi,
b_i\xi)= \|\xi_0\|^2<1$. Similarly, $\|\tilde\eta\|<1$. Since
$\|\nph\|_{\pi_1,\pi_2,\pi_3} = 1$, it now follows from
(\ref{ineq}) that
$$|(S(T)\xi_0,\eta_0)|\leq
\left\|\sum_{i=1}^s\sum_{k=1}^{\infty}\left(\sum_{l=1}^n\tilde
b_l^*(T_i^2)_{lk}\right)\otimes
\left(\sum_{m=1}^n(T_i^1)_{km}\tilde
a_m\right)\right\|_{\hh}\|\tilde\xi\|\|\tilde\eta\|\leq 1,$$ a
contradiction.

If $a$ or $b$ is not invertible, let $\epsilon > 0$ be such that
$\hat{\xi}_0 \stackrel{def}{=}
(\xi_1,\dots,\xi_n,\epsilon\xi,0,\dots)$ and $\hat{\eta}_0
\stackrel{def}{=} (\eta_1,\dots,\eta_n,\epsilon\eta,0,\dots)$ have
norm less than one and $|(S(T)\hat{\xi}_0,\hat{\eta}_0)| > 1$.
Choose $a_i$ and $b_i$ in the same way as before, and let $a_{n+1}
= \epsilon I$, $b_{n+1} = \epsilon I$, $a = \sum_{i=1}^{n+1}
a_i^*a_i$ and $b = \sum_{i=1}^{n+1} b_i^*b_i$. Then $a$ and $b$
are invertible and the proof proceeds in the same fashion.

We have proved that ${\bf M}_{\pi_1,\ldots,\pi_n}\subseteq{\bf
M}_{\pi_1\otimes 1,\ldots,\pi_n\otimes 1}$ and that
$\|\cdot\|_{\pi_1\otimes 1,\ldots,\pi_n\otimes 1}\leq
\|\cdot\|_{\pi_1,\ldots,\pi_n}$. The converse inequality is easy
to show, and thus the proof is complete. \prend

\begin{corollary}\label{kerin}
Let $\pi_i$ be a representation of the C*-algebra ${\mathcal
A}_i$, $i = 1,\dots,n$. Assume that $\pi_1$ and $\pi_n$ have
separating vectors. If
\begin{equation}\label{kernel}
\ker(\pi_i)\subseteq\ker(\pi_i'),\text{ for each } i=1,\ldots, n,
\end{equation}
then ${\bf M}_{\pi_1,\ldots,\pi_n}\subseteq {\bf
M}_{\pi_1',\ldots,\pi_n'}$ and
$\|\varphi\|_{\pi_1',\ldots,\pi_n'}\leq
\|\varphi\|_{\pi_1,\ldots,\pi_n}$, for each $\varphi\in {\bf
M}_{\pi_1,\ldots,\pi_n}$.
\end{corollary}
\proof The proof is similar to that  of \cite[Corollary~5.8]{ks};
we include it for completeness. Let ${\mathcal H}$ be an
infinite-dimensional Hilbert space of sufficiently large
dimension. Then (\ref{kernel}) implies
$$\text{rank}(\pi_i'(a_i))\leq \text{rank}(\pi_i(a_i)\otimes 1), \text{ for all
}a_i\in\cl A_i.
$$
By  (\ref{hadwin}), $\pi_i'\stackrel{a}\ll\pi_i\otimes 1$.
Applying now Theorem~\ref{ap} and then Lemma~\ref{smith+} we
obtain the statement. \prend

Using Corollary \ref{kerin} and results from \cite{ks} we will now
show that if the C*-algebras $\cl A_i$ are commutative then the
space ${\bf M}_{\pi_1,\ldots,\pi_n}(\cl A_1,\dots,\cl A_n)$ of
multipliers depends only on the supports of spectral measures
corresponding to the representations $\pi_i$.

Assume that ${\cl A_i}$ is commutative, $i = 1,\dots,n$ and let
$X_i$ be the maximal ideal spaces of $\cl A_i$; then $\cl A_i
\simeq C_0(X_i)$. Let $\pi_i$ be a representation of $\cl A_i$ and
$\cl E_{\pi_i}$ be the spectral measure on $X_i$ corresponding to
$\pi_i$.

It was proved in \cite[Lemma~7.2]{ks} that if $f\in C_0(X)$ and
the representation $\pi$ of $C_0(X)$ is such that $\text{rank
}(\pi(f))<\infty$ then $$\text{rank }(\pi(f))=\sum_{x\in S(f,{\cl
E}_{\pi})}\dim({\cl E}_{\pi}(\{x\})),$$ where $S(f,{\cl
E}_{\pi})=\{x\in\supp{\cl E}_{\pi}:f(x)\ne 0\}$. Thus the
condition
$$\supp{\cl E}_{\pi'}\subset\supp{\cl E}_{\pi}$$ implies
$\ker\pi(f)\subseteq\ker\pi'(f)$. As each representation $\pi$ of
a commutative algebra $C_0(X)$ has a separating vector we have the
following

\begin{corollary}\label{supp}
Let $\pi_i$, $\pi_i'$ be separable representations of the
C*-algebra $\cl A_i=C_0(X_i)$ and $\cl E_{\pi_i}$ and $\cl
E_{\pi_i'}$ be the corresponding spectral measures ($i =
1,\dots,n$). If
$$\supp{\cl E}_{\pi'_i}\subseteq\supp{\cl E}_{\pi_i},\text{ for each
}i=1,\ldots, n,$$ then ${\bf M}_{\pi_1,\ldots,\pi_n}\subseteq{\bf
M}_{\pi_1',\ldots,\pi_n'}$.
\end{corollary}

Let $\mu_i$ be measures on $X_i$. Let $\pi_i$ be a representation
of $C_0(X_i)$ on $L_2(X_i,\mu_i)$ defined by
$(\pi_i(f)h)(x_i)=f(x_i)h(x_i)$. We call $\varphi\in
C_0(X_1\times\ldots\times X_n)$ a {\it
$(\mu_1,\ldots,\mu_n)$-multiplier} if $\varphi\in {\bf
M}_{\pi_1,\ldots,\pi_n}$ and let $\|\nph\|_{\mu_1,\dots,\mu_n} =
\|\nph\|_{\pi_1,\dots,\pi_n}$.

By Corollary~\ref{supp}, the set of the all
$(\mu_1,\ldots,\mu_n)$-multipliers depends only on the supports of
measures $\mu_i$. The next statement shows the connection between
$(\mu_1,\ldots,\mu_n)$-multipliers and multidimensional Schur
multipliers (with respect to discrete measures).

\begin{corollary}
Let $X_i$ be locally compact spaces with countable bases and let
$\mu_i$ be Borel $\sigma$-finite measures on $X_i$ with
$\supp\mu_i=X_i$. Then $\varphi\in C_0(X_1\times\ldots\times X_n)$
is a $(\mu_1,\ldots,\mu_n)$-multiplier if and only if $\varphi$ is
a Schur multiplier on $X_1\times\ldots\times X_n$. Moreover, in
this case $\|\varphi\|_{\mu_1,\ldots,\mu_n}=\|S_{\varphi}\|$.
\end{corollary}

\proof The proof is similar to that of \cite[Theorem~7.5]{ks}.
\prend

\section{Universal multipliers}\label{s_universal}

The main goal of this section is to give a full description of the
multipliers which do not depend on the choice of the
representations of the C*-algebras $\mathcal{A}_1$,
$\mathcal{A}_2,\dots,\mathcal{A}_n$. Recall that an element
$\nph\in \mathcal{A}_1\otimes\ldots\otimes\mathcal{A}_n$ is called
a {\bf universal multiplier} if $\nph$ is a
$\pi_1,\pi_2,\ldots,\pi_n$-multiplier for all representations
$\pi_1$, $\pi_2$,..., $\pi_n$ of $\mathcal{A}_1$,
$\mathcal{A}_2,\dots,\mathcal{A}_n$, respectively. The set of all
universal multipliers in $\cl A_1\otimes\dots\otimes\cl A_n$ is
denoted by ${\bf M}(\mathcal{A}_1,\ldots,\mathcal{A}_n)$.

Along with the universal multipliers, we will describe another
class of multipliers which we call projective universal
multipliers and define as follows. Let $H_1,\dots,H_n$ be Hilbert
spaces. Equip $\Gamma(H_1,\dots,H_n)$ with the projective tensor
norm $\|\cdot\|_{\wedge}$, where each of the terms $H_i\otimes
H_{i+1}$ (resp. $H_{i-1}^{\dd}\otimes H_i^{\dd}$) is given its
operator norm. We call an element $\nph\in\cl
B(H_1\otimes\dots\otimes H_n)$ a concrete projective multiplier if
there exists $C > 0$ such that $\|S_{\nph}(\zeta)\|_{\op}\leq
C\|\zeta\|_{\wedge}$, for all $\zeta\in \Gamma(H_1,\dots,H_n)$. If
$\cl A_1,\dots\cl A_n$ are C*-algebras, an element $\nph\in \cl
A_1\otimes\dots\otimes\cl A_n$ will be called a {\bf projective
universal multiplier} if $(\pi_1\otimes\dots\otimes\pi_n)(\nph)$
is a concrete projective multiplier for all choices of
representations $\pi_1,\dots,\pi_n$ of $\cl A_1,\dots,\cl A_n$,
respectively. We denote by ${\bf
M}^{\wedge}(\mathcal{A}_1,\ldots,\mathcal{A}_n)$ the set of all
projective universal multipliers.

If $\nph\in {\bf M}(\mathcal{A}_1,\ldots,\mathcal{A}_n)$ let
$$\|\nph\|_{\univ}=\sup\limits_{\pi_1,\pi_2,\ldots,\pi_n}\|\nph\|_{\pi_1,\pi_2,\ldots,\pi_n}.$$
Note that $\|\nph\|_{\univ}$ is finite. In fact, assume that there
exist representations $\pi_{1,k}, \dots,\pi_{n,k}$, such that
$\|\nph\|_{\pi_{1,k},\pi_{2,k},\dots,\pi_{n,k}}\rightarrow_{k\rightarrow\infty}\infty$
and let $\pi_{1}=\bigoplus\limits_{k}\pi_{1,k}$,
$\pi_{2}=\bigoplus\limits_{k} \pi_{2,k}, \dots,$
$\pi_{n}=\bigoplus\limits_{k}\pi_{n,k}$. Then, by Theorem
\ref{ap},
$$\|\nph\|_{\pi_{1,k},\pi_{2,k},\ldots,\pi_{n,k}}\leq
\|\nph\|_{\pi_{1},\pi_{2},\ldots,\pi_{n}},$$ for all
$k\in\mathbb{N}$, which contradicts the fact that $\nph\in {\bf
M}(\mathcal{A}_1,\ldots,\mathcal{A}_n)$.

It is clear that ${\bf M}(\mathcal{A}_1,\ldots,\mathcal{A}_n)$ is
a linear subspace of $\cl A_1\otimes\dots\otimes\cl A_n$
containing $\mathcal{A}_1\odot\dots\odot \mathcal{A}_n$.

Recall that the Haagerup norm on $\mathcal{A}_{1}\odot
\mathcal{A}_2\odot\ldots\odot\mathcal{A}_{n}$ is
\begin{gather*}
\|\omega\|_{\hh}=\inf\{\|\omega_{1}\|\|\omega_{2}\|\ldots\|\omega_{n}\|:
 \omega=\omega_{1}\odot\omega_{2}\odot\ldots\odot\omega_{n},\\
\omega_{1}\in M_{1,i_1}(\mathcal{A}_1),\omega_2\in
M_{i_1,i_2}(\mathcal{A}_2),\ldots, \omega_{n}\in
M_{i_{n-1},1}(\mathcal{A}_{n}), i_1,i_2,\ldots,i_{n-1}\in
\mathbb{N}\}.
\end{gather*}

A modification of the Haagerup norm on the algebraic tensor
product of two $C^*$-algebras was considered in \cite{itoh, ks}.
We now introduce a natural generalisation of this norm for
arbitrary $n$. Recall the maps $\omega\mapsto\omega^{\ttt}$ and
$\omega\mapsto\omega^{\dd}$ on $M_n(\mathcal{A}) =
M_n(\bb{C})\otimes\cl A$ given on elementary tensors by $(a\otimes
b)^{\ttt}=a^{\ttt}\otimes b$ and $(a\otimes b)^{\dd}=a\otimes
b^{\dd}$ (here $\mathcal{A}$ is a $C^*$-subalgebra of $B(H)$ for
some Hilbert space $H$). We set
\begin{gather*}
\|\omega\|_{\ph}=\inf\{\prod\limits_{0\leq
i<\frac{n}{2}}\|\omega_{n-2i}^{\ttt}\|\|\omega_{n-2i-1}\|:
\omega=\omega_{1}\odot\omega_{2}\odot\ldots\odot\omega_{n}, \ \omega_0=I,\\
\omega_{1}\in M_{1,i_1}(\mathcal{A}_1),\omega_2\in
M_{i_1,i_2}(\mathcal{A}_2),\ldots, \omega_{n}\in
M_{i_{n-1},1}(\mathcal{A}_{n}), i_1,i_2,\ldots,i_{n-1}\in
\mathbb{N}\},
\end{gather*}
In the case $n = 2$, the above norm was denoted in \cite{itoh} by
$\|\cdot\|_{\hh'}$. Clearly, if the algebras ${\mathcal A}_i$, $i
= 1,\dots,n$, are commutative then the norms $\|\cdot\|_{\hh}$ and
$\|\cdot\|_{\ph}$ coincide. It was shown in \cite{itoh} that in
general they need not be even equivalent.

\begin{lemma}\label{K1}$\|\omega\|_{\univ}\leq\|\omega\|_{\ph}$ for all
$\omega\in \mathcal{A}_1\odot\ldots\odot\mathcal{A}_n$.
\end{lemma}
\proof Let $\pi_i$ be a representation of $\mathcal{A}_i$,
$i=1,\dots,n$, and let
$\omega=\omega_{1}\odot\omega_{2}\odot\ldots\odot\omega_{n}$,
where $\omega_{1}\in M_{1,k_1}(\mathcal{A}_1),\omega_2\in
M_{k_1,k_2}(\mathcal{A}_2),\ldots, \omega_{n}\in
M_{k_{n-1},1}(\mathcal{A}_{n})$ for some
$k_1,k_2,\ldots,k_{n-1}\in \mathbb{N}$.

Let $n$ be even, $\xi_{1,2}\in M_{1,l_1}(H_1\otimes H_2)$,
$\eta_{2,3}\in M_{l_1,l_2}(H_2^{d}\otimes H_3^{d}),\dots$,
$\xi_{n-1,n}\in M_{l_{n-2},1}(H_{n-1}\otimes H_n)$ and
$$\zeta =
\xi_{1,2}\odot\eta_{2,3}\odot\dots\odot\xi_{n-1,n}\in
\Gamma(H_1,\dots,H_n).$$ Letting $\pi=\pi_1\otimes\ldots\otimes
\pi_n$, by Lemma \ref{opmul2} we have
\begin{gather*}
S_{\pi(\omega)}(\zeta)=(\id\mbox{}_{1,k_{n-1}}\otimes\pi_n)(\omega_n^{\ttt})\dots
(\theta_{l_1,l_2}(\eta_{2,3})^{\ttt}\otimes I_{k_2})\\
\times ((\id\mbox{}_{k_1,k_{2}}\otimes
\pi_2)(\omega_2^{\ttt})\otimes
I_{l_1})(\theta_{1,l_1}(\xi_{1,2})^{\ttt}\otimes
I_{k_1})(\id\mbox{}_{k_1,1}\otimes\pi_1)(\omega_1^{\ttt})^{\dd}.\\
\end{gather*}
Since
$\|(\id_{k_{m-1},k_m}\otimes\pi_m)(\omega_{m}^{\ttt})^{\dd}\|=\|(\id_{k_{m-1},k_{m}}\otimes\pi_m)(\omega_{m})\|$,
we have
\begin{eqnarray*}\|S_{\pi(\omega)}(\zeta)\|_{\op} & \leq &
\|\theta_{1,l_1}(\xi_{1,2})^{\ttt}\|\dots
\|\theta_{l_{n-2},1}(\xi_{n-1,n})^{\ttt}\|\\ & \times &
\prod\limits_{0\leq
i<\frac{n}{2}}\|\omega_{n-2i}^{\ttt}\|\|\omega_{n-2i-1}\|
=\|\omega\|_{\ph} \|\zeta\|_{\hh}.
\end{eqnarray*}

Now let $n$ be odd and
$$\zeta = \eta_{1,2}\odot\xi_{2,3}\odot\dots\odot\xi_{n-1,n}\in
\Gamma(H_1,\dots,H_n),$$ where $\eta_{1,2}\in
M_{1,l_1}(H_1^{\dd}\otimes H_2^{\dd})$, $\xi_{2,3}\in
M_{l_1,l_2}(H_2\otimes H_3),\dots$, $\xi_{n-1,n}$ $\in$
$M_{l_{n-2},1}$ $(H_{n-1}\otimes H_n)$. Using the previously
obtained inequality, we have
\begin{eqnarray*}\|S_{\pi(\omega)}(\zeta)\|_{\op} & = &
\sup_{\|\xi\|\leq 1} \|S_{\pi(\omega)}(\zeta)(\xi)\|_{H_n}\\ & = &
\sup_{\|\xi\|\leq 1}
\|S_{\id\otimes\pi(1\otimes\omega)}((1\otimes\xi)\otimes\zeta)\|_{\cl
B(\bb{C}^{\dd},H_n)} \\ & \leq & \|\omega\|_{\ph}
\|\xi\|\|\zeta\|_{\hh}.
\end{eqnarray*}
The proof is complete. \prend

If $H_1,\dots,H_n$ are Hilbert spaces, we say that a net
$\{\nph_{\nu}\}\subseteq B(H_1\otimes\dots\otimes H_n)$ converges
semi-weakly to an operator $\nph\in B(H_1\otimes\dots\otimes H_n)$
if $(\nph_{\nu}\zeta_1,\zeta_2)\rightarrow (\nph\zeta_1,\zeta_2)$
for all $\zeta_1,\zeta_2\in H_1\odot\dots\odot H_n$. Note that if
the net $\{\nph_{\nu}\}$ is bounded then it converges semi-weakly
if and only if it converges weakly.

Let $\mathcal{A}_1\subseteq B(H_1)$, $\mathcal{A}_2\subseteq
B(H_2),\ldots$, $\mathcal{A}_n\subseteq B(H_n)$ be C*-algebras and
$(\mathcal{A}_1\odot \mathcal{A}_2\odot\ldots \odot
\mathcal{A}_n)^{\sharp}$ be the linear space of all $\nph\in
\mathcal{A}_1\otimes \mathcal{A}_2\otimes\ldots
\otimes\mathcal{A}_n$ for which there exists a net
$\{\nph_{\nu}\}\subseteq \mathcal{A}_1\odot
\mathcal{A}_2\odot\ldots \odot \mathcal{A}_n$ converging to $\nph$
semi-weakly (as a net of operators in $B(H_1\otimes
H_2\otimes\ldots \otimes H_n)$) and such that
$\sup\limits_{\nu}\|\nph_{\nu}\|_{\ph} < \infty$.

\begin{proposition}\label{subs}
Let $\mathcal{A}_i\subseteq \cl B(H_i)$, $i = 1,\dots,n$, be
C*-algebras. Then $(\mathcal{A}_1\odot\dots\odot
\mathcal{A}_n)^{\sharp}\subseteq {\bf
M}(\mathcal{A}_1,\dots,\mathcal{A}_n)\subseteq {\bf
M}^{\wedge}(\mathcal{A}_1,\dots,\mathcal{A}_n)$.
\end{proposition}
\proof Since $\|\zeta\|_{h}\leq\|\zeta\|_{\wedge} $ for all
$\zeta\in\Gamma(H_1,\ldots, H_n)$ we have ${\bf
M}(\mathcal{A}_1,\dots,\mathcal{A}_n)\subseteq {\bf
M}^{\wedge}(\mathcal{A}_1,\dots,\mathcal{A}_n)$.

Let us first prove that
$$(\mathcal{A}_1\odot\dots\odot \mathcal{A}_n)^{\sharp}\subseteq
{\bf M}_{\pi_1,\ldots,\pi_n}(\mathcal{A}_1,\dots,\mathcal{A}_n),$$
in the case where $\pi_i=\bigoplus\limits_{\lambda_i}\id$ is the
sum of $\lambda_i$ copies of the identity representation. Let
$\{\nph_{\nu}\}\subseteq \mathcal{A}_1\odot \ldots \odot
\mathcal{A}_n$ be a net converging semi-weakly to $\nph$ and such
that $D=\sup\limits_{\nu}\|\nph_{\nu}\|_{\ph}<\infty$ and
$\pi=\pi_1\otimes\ldots\otimes \pi_n$. By Lemma \ref{K1},
$$\|S_{\pi(\nph_{\nu})}(\zeta)\|_{\op}\leq D
\|\zeta\|_{\hh}$$ for all $\nu$ and $\zeta\in
\Gamma(H_1,\dots,H_n)$.

Suppose first that $n$ is even. To prove that
$\|S_{\pi(\nph)}(\zeta)\|_{\op}\leq D \|\zeta\|_{\hh}$, it
suffices to show that the net $\{S_{\pi(\nph_{\nu})}(\zeta)\}$ of
operators in $B(\widetilde{H}_1^{\dd},\widetilde{H}_n)$ converges
weakly to the operator $S_{\pi(\nph)}(\zeta)$ (here and in the
sequel we set $\widetilde{H}_i=\bigoplus\limits_{\lambda_i}H_i$,
$i = 1,\dots,n$). By linearity and the uniform boundedness of the
net $\{S_{\pi(\nph_{\nu})}(\zeta)\}$, it is suffices to prove that
$$(S_{\pi(\nph_{\nu})}(\zeta)x^{\dd},y)\rightarrow
(S_{\pi(\nph)}(\zeta)x^{\dd},y)$$ for all $x^{\dd}$ and $y$ which
have only one non-zero entry in the corresponding direct sums of
$H_1^{\dd}$ and $H_n$, respectively.

Fix such $x^{\dd}$ and $y$, and let
$\zeta=\xi_{1,2}\otimes\eta_{2,3}^{\dd} \otimes\ldots
\otimes\xi_{n-1,n}\in
\Gamma(\widetilde{H}_1,\ldots,\widetilde{H}_n)$. Then
$$(S_{\pi(\nph_{\nu})}(\zeta)x^{\dd},y)= (\pi(\nph_{\nu})(\xi_{1,2}\otimes\ldots
\otimes\xi_{n-1,n}),
x\otimes\eta_{2,3}\otimes\eta_{4,5}\otimes\ldots
\otimes\eta_{n-2,n-1}\otimes y)$$ Indeed, assuming $n=4$ for the
simplicity we get
\begin{eqnarray*}
(S_{\pi(\nph_{\nu})}(\zeta)x^{\dd},y)&=&
(\sigma_{\pi}(\nph_{\nu})\theta(\xi_{1,2}\otimes\xi_{3,4})(\theta(\eta_{2,3}^{\dd})),\theta(x\otimes y))_2\\
&=&(\sigma_{\pi}(\nph_{\nu})\theta(\xi_{1,2}\otimes\xi_{3,4}),\theta(\theta(\eta_{2,3})\otimes\theta(x\otimes
y)))_2\\
&=&(\sigma_{\pi}(\nph_{\nu})\theta(\xi_{1,2}\otimes\xi_{3,4}),\theta(x\otimes\eta_{2,3}\otimes y))_2\\
& = &
(\pi(\nph_{\nu})(\xi_{1,2}\otimes\xi_{3,4}),x\otimes\eta_{2,3}\otimes
y).
\end{eqnarray*}

Fix $\epsilon > 0$ and let $\tilde{\zeta} =
\tilde{\xi}_{1,2}\otimes\tilde{\eta}_{2,3}^{\dd}\otimes\dots\otimes\tilde{\xi}_{n-1,n}$
be such that all norms $\|\xi_{1,2} - \tilde{\xi}_{1,2}\|$,
$\|\eta_{2,3} - \tilde{\eta}_{2,3}\|$,$\dots$,$\|\xi_{n-1,n} -
\tilde{\xi}_{n-1,n}\|$ are smaller than $\epsilon$ and all vectors
$\tilde{\xi}_{1,2}$,
$\tilde{\eta}_{2,3}^{\dd}$,$\dots$,$\tilde{\xi}_{n-1,n}$ are
finite sums of elementary tensors which have only finitely many
non-zero entries in the direct sums of the corresponding Hilbert
spaces. Thus, we may assume that $\tilde{\xi}_{1,2}\in
H_1^{(k)}\odot H_2^{(k)}$,$\tilde{\eta}_{2,3}\in H_2^{(k)}\odot
H_3^{(k)}$ $\dots$,$\tilde{\xi}_{n-1,n}\in H_{n-1}^{(k)}\odot
H_n^{(k)}$, $x^{\dd}\in H_1^{(k)}$ and $y\in H_n^{(k)}$ for some
$k\in\bb{N}$.

It follows from the formula above that there exists $\nu_0$ such
that if $\nu\geq\nu_0$ then
$$|(S_{\pi(\nph_{\nu})}(\tilde{\zeta})x^{\dd},y) -
(S_{\pi(\nph)}(\tilde{\zeta})x^{\dd},y)| < \epsilon.$$ On the
other hand,
\begin{eqnarray*}
& & |(S_{\pi(\nph_{\nu})}(\zeta)x^{\dd},y) -
(S_{\pi(\nph_{\nu})}(\tilde{\zeta})x^{\dd},y)|\\ & \leq &
D\|x\|\|y\|\|\tilde{\zeta}-\zeta\|_{\hh} \leq
(C+\epsilon)^{n-2}D(n-1)\|x\|\|y\|\epsilon,
\end{eqnarray*}
for every $\nu$, where $C =
\max\{\|\xi_{1,2}\|,\|\eta_{2,3}\|\dots,\|\xi_{n-1,n}\|\}$. Using
Remark \ref{r_projnorm}, we have
\begin{eqnarray*}
& & |(S_{\pi(\nph)}(\zeta)x^{\dd},y) -
(S_{\pi(\nph)}(\tilde{\zeta})x^{\dd},y)|\\ & \leq &
\|\nph\|\|x\|\|y\|\|\zeta-\tilde\zeta\|_{2,\wedge}\leq
\|\nph\|(C+\epsilon)^{n-2}(n-1)\|x\|\|y\|\epsilon.
\end{eqnarray*}
Thus,
\begin{eqnarray*}
& & |(S_{\pi(\nph_{\nu})}(\zeta)x^{\dd},y) -
(S_{\pi(\nph)}(\zeta)x^{\dd},y)|\\ & \leq & \epsilon(1 +
(C+\epsilon)^{n-2}D(n-1)\|x\|\|y\| +
\|\nph\|(C+\epsilon)^{n-2}(n-1)\|x\|\|y\|)
\end{eqnarray*}
whenever $\nu\geq\nu_0$. It follows that the net
$\{S_{\pi(\nph_{\nu})}(\zeta)\}$ converges weakly to
$S_{\pi(\nph)}(\zeta)$ and hence $\nph\in$ ${\bf
M}_{\pi_1,\dots,\pi_n}$ $(\mathcal{A}_1,\dots,\mathcal{A}_n)$.

In the case $n$ is odd, a calculation similar to the one above
shows that $(S_{\pi(\nph_{\nu})}(\zeta)x,y)$ is equal to
$$(\pi(\nph_{\nu})(x\otimes\xi_{2,3}\otimes\ldots
\otimes\xi_{n-1,n}),\eta_{1,2}\otimes\ldots
\otimes\eta_{n-2,n-1}\otimes y),$$ whenever $x$ $\in$
$\widetilde{H}_1$, $y$ $\in$ $\widetilde{H}_n$,
$\zeta=\eta_{1,2}^{\dd}\otimes\xi_{2,3} \otimes\ldots
\otimes\xi_{n-1,n}\in\Gamma(\widetilde{H}_1,\ldots,\widetilde{H}_n)$,
and the proof proceeds in a similar fashion.

Now let $\pi_1$, $\dots$, $\pi_n$ be representations of
$\mathcal{A}_1$, $\dots$, $\mathcal{A}_n$ on $H_{\pi_1}$, $\dots$,
$H_{\pi_n}$ and $\pi=\pi_1\otimes\ldots\otimes \pi_n$. Then
$$\text{rank}(\pi_i(a_i))\leq
\text{rank}\left(\bigoplus\limits_{dim(H_{\pi_i})}\id(a_i)\right)$$
for all $a_i\in \mathcal{A}_i$ and $i=1,..,n$. By Theorem \ref{ap}
(i),
$${\bf M}_{\bigoplus\limits_{\lambda_1} \id,\bigoplus\limits_{\lambda_2}
\id,\ldots ,\bigoplus\limits_{\lambda_k}
\id}(\mathcal{A}_1,\ldots,\mathcal{A}_n)\subseteq {\bf
M}_{\pi_1,\pi_2,\ldots,\pi_k}(\mathcal{A}_1,\mathcal{A}_2,\ldots,\mathcal{A}_n).$$
The proof is complete. \prend

Assume that $n$ is even. Then the mapping $S_{\id(\nph)}$ acting on
$\Gamma(H_1$, $\dots$, $H_n)$ $=$ $(H_1\otimes
H_{2})\odot(H_2^{\dd}\otimes H_{3}^{\dd})\odot\ldots
\odot(H_{n-1}\otimes H_{n})$ can be regarded as a mapping on the
algebraic tensor product
\begin{gather}\label{HSG}
HS(H_{n-1},H_{n})\odot HS(H_{n-2},H_{n-1})^{\dd}\odot\ldots \odot
HS(H_{1},H_{2})
\end{gather}
of the corresponding spaces of Hilbert-Schmidt operators by
letting
$$S_{\nph}(\theta(\xi_{n-1,n})\otimes\theta(\eta_{n-2,n-1})^{\dd}\otimes\theta(\xi_{n-3,n-2})
\otimes\ldots \otimes\theta(\xi_{1,2}))=S_{\nph}(\zeta),$$ where
$\zeta=\xi_{1,2}\otimes\eta_{2,3}^{\dd}\otimes\xi_{3,4}\otimes\ldots
\otimes\xi_{n-1,n}$. Denote the space (\ref{HSG}) by $HS\Gamma$
$(H_1,\ldots,H_n)$. If $\nph$ is an elementary tensor then Lemma
\ref{opmul2} (i) shows that $S_{\id(\nph)}$ is $\mathcal{A}_n',
(\mathcal{A}_{n-1}^{\dd})',\ldots,\mathcal{A}_2',({\mathcal{A}_1}^{\dd})'$-modular.
It follows by continuity that $S_{\id(\nph)}$ is $\mathcal{A}_n',
(\mathcal{A}_{n-1}^{\dd})',\dots,\mathcal{A}_2',({\mathcal{A}_1}^{\dd})'$-modular
for every $\nph\in \cl A_1\otimes\dots\otimes\cl A_n$. If moreover
$\nph\in {\bf M}_{\id,\dots,\id}(\mathcal{A}_1,\ldots,
\mathcal{A}_n)$ then $S_{\id(\nph)}$ can be extended to a bounded
mapping (denoted in the same way) from the algebraic tensor product
$$\cl K(H_{n-1}^{\dd},H_{n})\odot \cl K(H_{n-2}^{\dd},H_{n-1})^{\dd}\odot\dots\odot
\cl K(H_{1}^{\dd},H_{2})$$ into $\cl K(H_1^{\dd},H_n)$. By
continuity, this extension is $\mathcal{A}_n',
(\mathcal{A}_{n-1}^{\dd})',\dots,\mathcal{A}_2',({\mathcal{A}_1}^{\dd})'$-modular.

Similarly, if $n$ is odd and $\nph\in {\bf
M}_{\id,\dots,\id}(\mathcal{A}_1,\ldots, \mathcal{A}_n)$ then
$S_{\id(\nph)}$ can be regarded as a multilinear $\mathcal{A}_n',
({\mathcal{A}_{n-1}}^{\dd})',\dots,({\mathcal{A}_2}^{\dd})',\mathcal{A}_1'$-modular
map from
$$\cl K(H_{n-1}^{\dd},H_{n})\odot \cl K(H_{n-2}^{\dd},H_{n-1})^{\dd}\odot\dots\odot
\cl K(H_{1}^{\dd},H_{2})$$ into $\cl B(H_1,H_n)$. Denote by ${\bf
M}_{\id,\dots,\id}^{cb}(\mathcal{A}_1,\ldots , \mathcal{A}_n)$ the
set of all $(\id,\ldots ,\id)$-multipliers for which the mapping
$S_{\id(\nph)}$ is completely bounded.

\begin{proposition}\label{MIdId}
Let $\mathcal{A}_i\subseteq \cl B(H_i)$, $i = 1,\dots,n$, be von
Neumann algebras. Then ${\bf
M}_{\id,\dots,\id}^{cb}(\mathcal{A}_1, \dots,\mathcal{A}_n)
\subseteq(\mathcal{A}_1\odot\dots\odot \mathcal{A}_n)^{\sharp}$.
\end{proposition}
\proof
Assume first that $n$ is even. For notational simplicity we assume
that $H_i$ is separable, $i=1,\dots,n$. Let $\id : \cl
A_1\otimes\dots\otimes\cl A_n\rightarrow\cl
B(H_1\otimes\dots\otimes H_n)$ be the identity representation.

Let $\nph\in {\bf M}_{\id,\ldots ,\id}^{cb}(\mathcal{A}_1,
\dots,\mathcal{A}_n)$. Then $S_{\id(\nph)}$ is a multilinear
$\mathcal{A}_n'$, $(\mathcal{A}_{n-1}^{\dd})'$, $\ldots$,
$\mathcal{A}_2'$, $({\mathcal{A}_1}^{\dd})'$-modular mapping on
$$\cl K(H_{n-1}^{\dd},H_{n})\odot
\cl K(H_{n-2},H_{n-1}^{\dd})\odot\dots \odot \cl
K(H_{1}^{\dd},H_{2})$$ taking values in $\cl K(H_1^{\dd},H_n)$.
Let $H^{\infty}=H\otimes l^2$ and $I_{\infty}$ be the identity
operator on $l^2$.

Since $S_{\id(\nph)}$ is completely bounded, it extends to a
completely bounded mapping, denoted in the same way, from $$\cl
K(H_{n-1}^{\dd},H_{n})\otimes_{\hh} \cl
K(H_{n-2},H_{n-1}^{\dd})\otimes_{\hh}\dots \otimes_{\hh} \cl
K(H_{1}^{\dd},H_{2})$$ into $\cl K(H_1^{\dd},H_n)$. Then the
second dual $S_{\id(\nph)}^{**}$ is a weak* continuous completely
bounded mapping from $\cl B(H_{n-1}^{\dd},H_{n})$ $\otimes_{\sigma
\hh}$ $\dots$ $\otimes_{\sigma \hh}$ $\cl B(H_{1}^{\dd},H_{2})$
into $\cl B(H_1^{\dd},H_n)$ and hence gives rise to a weak*
continuous completely bounded $\mathcal{A}_n'$,
$(\mathcal{A}_{n-1}^{\dd})'$, $\ldots$, $\mathcal{A}_2'$,
$({\mathcal{A}_1}^{\dd})'$-modular multilinear map, denoted in the
same way, from
$$\cl B(H_{n-1}^{\dd},H_{n})\times \cl
B(H_{n-2},H_{n-1}^{\dd})\times\dots \times \cl
B(H_{1}^{\dd},H_{2})$$ into $\cl B(H_1^{\dd},H_n)$.

It follows from Corollary 5.9 of \cite{cs} that there exist
bounded linear operators $A_{1}:H_1^{\dd}\rightarrow
(H_1^{\dd})^{\infty}$, $A_j:H_j^{\infty}\rightarrow H_j^{\infty}$,
if $j$ is even, $A_j:(H_j^{\dd})^{\infty}\rightarrow
(H_j^{\dd})^{\infty}$ if $j$ is odd ($j = 2,\dots,n-1$) and
$A_n:H_n^{\infty}\rightarrow H_n$ such that the entries of $A_j$
with respect to the corresponding direct sum decomposition belong
to $\mathcal{A}_j''=\mathcal{A}_j$ for even $j$ and to
$(\mathcal{A}_j^{\dd})''=\mathcal{A}_j^{\dd}$ for odd $j$,
$$S_{\id(\nph)}(\zeta)=A_n(\theta(\xi_{n-1,n})\otimes
I_{\infty})A_{n-1}(\theta(\eta_{n-2,n-1})^{\dd}\otimes
I_{\infty})A_{n-2}\ldots A_1,$$ for all
$$\zeta=\theta(\xi_{n-1,n})\otimes\theta(\eta_{n-2,n-1})^{\dd}\otimes\ldots
\otimes\theta(\xi_{1,2})\in HS\Gamma(H_1,\ldots ,H_n),$$ and
$$\|S_{\id(\nph)}\|_{cb}=\prod\limits_{1\leq i\leq
n}\|A_i\|.$$

Let $P_{m,\nu}=(p^{m}_{ij})_{i,j=1}^{\infty}$ be the projection
with $p^{m}_{ij}\in B(H_m)$ (resp. $p^{m}_{ij}\in B(H_m^{\dd})$),
$p^{m}_{ii}=I_{H_m}$ (resp. $p^{m}_{ii}=I_{H_m^{\dd}}$) if $m$ is
even (resp. if $m$ is odd) and $1\leq i\leq \nu$, and
$p^{m}_{ij}=0$ otherwise.

Set $\nph_{\nu}=A_1^{\dd,\ttt}P_{1,\nu}^{\dd}\odot P_{2,\nu} A_2
P_{2,\nu}\odot P_{3,\nu} A_3^{\dd} P_{3,\nu}\ldots\odot
P_{n,\nu}A_n$. Clearly, $\|\nph_{\nu}\|_{\ph}\leq
\prod\limits_{1\leq i\leq n}\|A_i\|$ for each $\nu$; it hence
suffices to prove that $\{\nph_{\nu}\}$ converges semi-weakly to
$\nph$.

As $S_{\id(\nph_{\nu})}(\zeta)$ equals
$$A_nP_{n,\nu}(\theta(\xi_{n-1,n})\otimes
I_{\infty})P_{n-1,\nu}A_{n-1}P_{n-1,\nu}(\theta(\eta_{n-2,n-1})^{\dd}\otimes
I_{\infty})\ldots P_{1,\nu}A_1$$ and $P_{l,\nu}$ converges
strongly to $I_{H_{l}}$, we have that $S_{\id(\nph_{\nu})}(\zeta)$
converges weakly to $S_{\id(\nph)}(\zeta)$. By the proof of
Proposition \ref{subs}, if $x^{\dd}\in H_1^{\dd}$, $y\in H_n$ and
$\psi\in {\mathcal A}_1\otimes\ldots\otimes{\mathcal A}_n$ then
$(S_{\id(\psi)}(\zeta)x^{\dd},y)$ equals
\begin{eqnarray*}
& & (\sigma_{\id}(\psi) \theta(\xi_{1,2}\otimes\ldots
\otimes\xi_{k-1,k}), \theta(x\otimes\eta_{2,3}\otimes\ldots
\otimes\eta_{k-2,k-1}\otimes
y))_2\\
& = & (\psi(\xi_{1,2}\otimes\ldots \otimes\xi_{k-1,k}),
x\otimes\eta_{2,3}\otimes\ldots \otimes\eta_{k-2,k-1}\otimes y).
\end{eqnarray*}
Thus $\nph_{\nu}$ converges semi-weakly to $\nph$ and therefore
$\nph\in (\mathcal{A}_1\odot \ldots \odot
\mathcal{A}_n)^{\sharp}$, giving the inclusion ${\bf
M}_{\id,\dots,\id}^{cb}(\mathcal{A}_1,$ $\dots,\mathcal{A}_n)$
$\subseteq$ $(\mathcal{A}_1\odot$ $\ldots$ $\odot$
$\mathcal{A}_n)^{\sharp}$.


Now assume that $n$ is odd. In this case $S_{\id(\nph)}^{**}$ is a
weak* continuous completely bounded multilinear $\mathcal{A}_n'$,
$(\mathcal{A}_{n-1}^{\dd})'$, $\ldots$, $(\mathcal{A}_2^{\dd})'$,
${\mathcal{A}_1}'$-modular mapping on
$$\cl B(H_{n-1}^{\dd},H_{n})\times
\cl B(H_{n-2},H_{n-1}^{\dd})\times\dots \times \cl
B(H_{1},H_{2}^{\dd})$$ taking values in $\cl B(H_1,H_n)^{**}$. Let
$Q$ be the weak* continuous projection from $\cl B(H_1,H_n)^{**}$
onto $\cl B(H_1,H_n)$. Then $Q\circ S_{\id(\nph)}^{**}$ takes
values in $\cl B(H_1,H_n)$, and coincides with $S_{\id(\nph)}$ on
$HS\Gamma(H_1,\ldots ,H_n)$. The proof now proceeds as above.
\prend

\begin{proposition}\label{MinMcb}
Let $\mathcal{A}_i\subseteq \cl B(H_i)$, $i = 1,\dots,n$, be
C*-algebras. Then ${\bf M}^{\wedge}(\mathcal{A}_1,$ $\dots,$
$\mathcal{A}_n)$ $\subseteq$ ${\bf
M}_{\id,\dots,\id}^{cb}(\mathcal{A}_1,\dots,\mathcal{A}_n)$.
\end{proposition}
\proof Let $\nph\in {\bf
M}^{\wedge}(\mathcal{A}_1,\ldots,\mathcal{A}_n)$. Then there
exists a constant $D>0$ such that
$$\|\sigma_{\pi_{1}\otimes\ldots\otimes\pi_{n}}(\nph)(\zeta)\|_{\op}\leq
D\|\zeta\|_{\wedge}$$ for all $\zeta\in \Gamma(H_1,\ldots ,H_n)$
and all representations $\pi_{1}$, $\dots$, $\pi_n$ of $\cl A_1,
\dots, \cl A_n$, respectively.

Let $k\in\bb{N}$. The space $HS\Gamma(H_1^k,\ldots ,H_n^k)$ is
naturally isomorphic to
\begin{gather}\label{mn}M_k(HS(H_{n-1},H_n))\odot
M_k(HS(H_{n-2},H_{n-1})^{\dd})\odot\ldots \odot
M_k(HS(H_1,H_2)),\end{gather} and thus the mapping $S_{(\id\otimes
1_k)\otimes\ldots \otimes(\id\otimes 1_k)(\nph)}$ is well-defined
on the space (\ref{mn}). One can easily check that
\begin{equation}\label{dotim}
S^{(k)}_{\id\otimes\ldots \otimes \id(\nph)}(\Xi_{n-1}\odot\ldots
\odot\Xi_1)= S_{(\id\otimes 1_k)\otimes\ldots \otimes(\id\otimes
1_k)(\nph)}(\Xi_{n-1}\otimes\ldots \otimes\Xi_1),
\end{equation}
where $\Xi_i\in M_k(HS(H_{i},H_{i+1}))$ (resp. $\Xi_i\in
M_k(HS(H_{i},H_{i+1})^{\dd})$) if $i$ is even (resp, if $i$ is
odd) and $\Xi_i\in M_k(HS(H_{i},H_{i+1})^{\dd})$ (resp. $\Xi_i\in
M_k(HS(H_{i},H_{i+1}))$) if $i$ is odd (resp, if $i$ is even). If
the matrices $\Xi_i$ are of arbitrary sizes such that the product
$\Xi_{n-1}\odot\ldots \odot\Xi_1$ is well defined then they may be
considered as square matrices, all of the same size, by
complementing with zeros, and identity (\ref{dotim}) will still
hold. It follows that
\begin{gather*}
\|S^{(k)}_{\id \otimes\ldots  \otimes \id(\nph)}(\Xi_1\odot\ldots
\odot\Xi_{n-1})\|_{\op}\leq D\prod\limits_{1\leq i\leq
n-1}\|\Xi_i\|_{\op},\text{ for all }\Xi_1,\ldots \Xi_{n-1},
\end{gather*}
and hence the mapping $S_{\id\otimes \ldots\otimes \id(\nph)}$ is
completely bounded and $\varphi$ is an $(\text{id},\ldots,$
$\text{id})$-multiplier. \prend

\begin{theorem}\label{sista}
Let $\mathcal{A}_i\subseteq \cl B(H_i)$, $i = 1,\dots,n$, be
C*-algebras. Then ${\bf M}(\mathcal{A}_1,$ $\dots,$
$\mathcal{A}_n)={\bf M}^{\wedge}(\mathcal{A}_1,$ $\dots,$
$\mathcal{A}_n)=(\mathcal{A}_1\odot\dots\odot\mathcal{A}_n)^{\sharp}$.
\end{theorem}
\proof By Propositions \ref{subs}, \ref{MIdId} and \ref{MinMcb},
$${\bf M}^{cb}_{\id,\ldots ,\id}(\mathcal{A}_1'',\ldots ,
\mathcal{A}_n'')=(\mathcal{A}_1''\odot \ldots \odot
\mathcal{A}_n'')^{\sharp}.$$ Evidently,
$${\bf M}^{cb}_{\id,\ldots ,\id}(\mathcal{A}_1,\dots,\mathcal{A}_n)\subseteq {\bf
M}^{cb}_{\id,\ldots ,\id}(\mathcal{A}_1'',
\dots,\mathcal{A}_n'')\cap (\mathcal{A}_1\otimes \ldots\otimes
\mathcal{A}_n).$$ Applying Propositions  \ref{subs}, \ref{MIdId}
and \ref{MinMcb}, we obtain
\begin{eqnarray*}
(\mathcal{A}_1\odot \ldots \odot
\mathcal{A}_n)^{\sharp}&\subseteq& {\bf M}(\mathcal{A}_1, \ldots ,
\mathcal{A}_n)\\&\subseteq& {\bf M}^{\wedge}(\mathcal{A}_1, \ldots
, \mathcal{A}_n)\\&\subseteq& {\bf M}^{cb}_{\id,\ldots
,\id}(\mathcal{A}_1,\ldots,\mathcal{A}_n)\\&\subseteq& {\bf
M}^{cb}_{\id,\dots,\id}(\mathcal{A}_1'',\dots,\mathcal{A}_n'')\cap
(\mathcal{A}_1\otimes \ldots\otimes
\mathcal{A}_n)\\&=&(\mathcal{A}_1''\odot\ldots \odot
\mathcal{A}_n'')^{\sharp}\cap (\mathcal{A}_1\otimes \ldots\otimes
\mathcal{A}_n).
\end{eqnarray*}
It hence suffices to show that
$$(\mathcal{A}_1''\odot \ldots \odot
\mathcal{A}_n'')^{\sharp}\cap \mathcal{A}_1\otimes \ldots\otimes
\mathcal{A}_n\subseteq (\mathcal{A}_1\odot\ldots
\odot\mathcal{A}_n)^{\sharp}.$$

Let $\nph\in (\mathcal{A}_1''\odot\ldots \odot
\mathcal{A}_n'')^{\sharp}\cap (\mathcal{A}_1\otimes \ldots\otimes
\mathcal{A}_n)$. Then there exists a net $\{\nph_{\nu}\}_{\nu\in
J}\subseteq \mathcal{A}_1''\odot\ldots \odot \mathcal{A}_n''$ with
$\sup\limits_{\nu}\|\nph_{\nu}\|_{ph}<\infty$ which converges
semi-weakly to $\nph$. Write $\nph_{\nu}=A_{1,\nu}\odot\ldots\odot
A_{n,\nu}$, where $A_{1,\nu}\in
M_{1,i_1}(\mathcal{A}_1''),A_{2,\nu}\in
M_{i_1,i_2}(\mathcal{A}_2''),\ldots, A_{n,\nu}\in
M_{i_{n},1}(\mathcal{A}_n'')$.

By Kaplansky's Density Theorem for TRO's \cite{harris}, for each
pair $(m,\nu)$ there exists a net
$\{A_{m,\nu,\tau(m)}\}_{\tau(m)}\subset
 M_{i_{m-1},i_m}(\mathcal{A}_m)$ converging strongly to $A_{m,\nu}$ and such that
$||A_{m,\nu,\tau(m)}||\leq||A_{m,\nu}||$ for all $\tau(m)$. Thus
if $A_{\nu,\tau}=A_{1,\nu,\tau(1)}\odot
A_{2,\nu,\tau(2)}\odot\ldots\odot A_{n,\nu,\tau(n)}$, where
$\tau=(\tau(1),\ldots,\tau(n))$, then the net
$\{A_{\nu,\tau}\}_{\tau}$ converges strongly to $\varphi_{\nu}$
and $||A_{\nu,\tau}||_{ph}\leq||\varphi_\nu||_{ph}$.



Let $\cl U$ be the collection of all weak neighbourhoods of $0$ of
the form $\{S\in \cl B(H_1\otimes\dots\otimes H_n) :
|(S(\zeta_1^j),\zeta_2^j)| < \epsilon_j, j = 1,\dots,k\}$, where
$\zeta_1^j,\zeta_2^j\in H_1\odot\dots\odot H_n$ and $\epsilon_j >
0$, $j = 1,\dots,k$. Note that $\cl U$ is directed with respect to
reverse inclusion. The convergence of the net
$\{\nph_{\nu}\}_{\nu\in J}$ semi-weakly to $\nph$ implies that for
every $U\in\cl U$ there exists $\nu(U)$ such that for every
$\lambda \in J$ with $\lambda\geq \nu(U)$, we have that
$\nph_{\lambda}-\nph\in U$. The convergence of
$\{A_{\nu,\tau}\}_{\tau}$ to $\nph_{\nu}$ implies the existence of
$T(\nu(U),U)$ such that for every $\tau\geq T(\nu(U),U)$, we have
that $A_{\nu(U),\tau}-\nph_{\nu(U)}\in U$. Consider the net
$A_{U}=A_{\nu(U),T(\nu(U),U)}$ indexed by $\cl U$. It is easy to
check that $A_{U}$ converges semi-weakly to $\nph$. The proof is
complete. \prend

Note that in Theorem \ref{sista} we actually proved that if $n$ is
even, $\nph\in {\bf M}(\mathcal{A}_1,\dots,\mathcal{A}_n)$, $\zeta
= \xi_{1,2}\otimes\dots\otimes\xi_{n-1,n}
\in\Gamma(H_1,\dots,H_n)$ and
$$S_{\id\otimes\dots\otimes\id(\nph)}(\zeta) =
A_n(\theta(\xi_{n-1,n})\otimes I)\dots(\theta(\xi_{1,2})\otimes
I)A_1^{\dd},$$ where $A_i$ for $i$ even (resp. $A_i^{\dd}$ for $i$
odd) is a bounded block operator matrix with entries in $\cl
A_i''$ (resp. $(\cl A_i^{\dd})''$), then there exists a net
$\nph_{\nu} = A_1^{\nu}\odot A_2^{\nu}\odot\dots\odot A_n^{\nu}$,
where $A_i^{\nu}$ is a finite block operator matrix with entries
in $\cl A_i$ such that $\nph_{\nu}\rightarrow\nph$ semi-weakly,
$A_i^{\nu}\rightarrow A_i$ (resp. $A_i^{\nu \dd}\rightarrow
A_i^{\dd}$) strongly for $i$ even (resp. for $i$ odd) and all
operator norms $\|A_i^{\nu}\|$, $\|A_i\|$ are bounded by a
constant depending only on $n$. A similar statement holds in the
case $n$ is odd.

\medskip

Denote by $(\mathcal{A}_1\odot\ldots \odot\mathcal{A}_n)^{\sim}$
the set of all $\nph\in\mathcal{A}_1\otimes \ldots\otimes
\mathcal{A}_n$ for which there exists a net
$\{\nph_{\nu}\}\subseteq
\mathcal{A}_1\odot\dots\odot\mathcal{A}_n$, such that
$\sup\limits_{\nu}\|\nph_{\nu}\|_{\ph}<\infty$ and if $\pi_i$ is
an irreducible representation of $\cl A_i$, $i = 1,\dots,n$, then
$\{(\pi_1\otimes\ldots\otimes\pi_n)(\nph_{\nu})\}$ converges
semi-weakly to $(\pi_1\otimes\ldots\otimes\pi_n)(\nph)$. Note that
if $\sup\limits_{\nu}\|\nph_{\nu}\|_{min}<\infty$, which holds for
example when the norms $\|\cdot\|_{\ph}$ and $\|\cdot\|_{\hh}$ are
equivalent (see \cite{itoh}), then in the definition of the space
$(\mathcal{A}_1\odot\ldots \odot\mathcal{A}_n)^{\sim}$ the
semi-weak convergence can be replaced by the convergence in the
weak operator topology.

It follows from \cite{ks} that if $\cl A$ and $\cl B$ are
commutative C*-algebras then ${\bf
M}(\mathcal{A},\mathcal{B})=(\mathcal{A}\odot\mathcal{B})^{\sim}$.
As a corollary of Theorem \ref{sista}, we show that the same
equality holds for an arbitrary number of arbitrary C*-algebras,
giving an answer to a problem posed in \cite{ks}.

\begin{theorem}\label{desc}
Let $\mathcal{A}_i$, $i = 1,\dots,n$, be C*-algebras. Then
$${\bf
M}(\mathcal{A}_1, \dots, \mathcal{A}_n)={\bf
M}^{\wedge}(\mathcal{A}_1,\dots,
\mathcal{A}_n)=(\mathcal{A}_1\odot \ldots \odot
\mathcal{A}_n)^{\sim}.$$
\end{theorem}
\proof Let $\pi_1=\bigoplus\limits_{\pi\in IrrRep
(\mathcal{A}_1)}\pi,\ldots , \pi_n = \bigoplus\limits_{\pi\in
IrrRep (\mathcal{A}_n)}\pi$, where $IrrRep(\mathcal{A}_i)$ is a
set whose elements are all inequivalent irreducible
representations of $\mathcal{A}_i$. Then
\begin{eqnarray*}
{\bf M}(\mathcal{A}_1,\dots,
\mathcal{A}_n)&=&(\pi_1\otimes\ldots\otimes
\pi_n)^{-1}(\pi_1(\mathcal{A}_1)\odot\ldots \odot
\pi_n(\mathcal{A}_n))^{\sharp}\\&\subseteq&(\mathcal{A}_1\odot\ldots
\odot \mathcal{A}_n)^{\sim}.
\end{eqnarray*}
Using arguments similar to the ones from the proof of Proposition
\ref{subs}, one can show that
$$(\mathcal{A}_1\odot \ldots \odot
\mathcal{A}_n)^{\sim}\subseteq {\bf M}(\mathcal{A}_1, \ldots ,
\mathcal{A}_n),$$ which together with Theorem~\ref{sista} gives
the statement of the theorem. \prend

\end{document}